\documentclass[a4paper,reqno]{amsart}
\usepackage{amsfonts}
\usepackage[dvips]{graphicx}

\catcode`@=11 
\def\blfootnote{\xdef\@thefnmark{}\@footnotetext}
\catcode`@=12 

\theoremstyle{plain}
\begingroup
\newtheorem{theorem}{Theorem}[section]
\newtheorem{lemma}[theorem]{Lemma}
\newtheorem{proposition}[theorem]{Proposition}

\endgroup

\theoremstyle{definition}
\begingroup
\newtheorem{definition}[theorem]{Definition}
\newtheorem{remark}[theorem]{Remark}

\endgroup

\theoremstyle{remark}
\begingroup

\endgroup

\mathchardef\emptyset="001F

\numberwithin{equation}{section}

\newcommand{\thone}{\theta_{\e,\eta}^{1}}
\newcommand{\tone}{t_{\e,\eta}^{1}}
\newcommand{\ttwo}{t_{\e,\eta}^{2}}

\newcommand{\Ttwo}{t_{\eta}^{2}}
\newcommand{\Ttwota}{t_{\delta}^{2}}
\newcommand{\tthree}{t_{\e}^{3}}
\newcommand{\tfive}{t_{\e,\eta}^{5}}
\newcommand{\thfour}{\theta_{\e,\rho}^{4}}
\newcommand{\tauf}{\tau_{\e,\rho}^{4}}
\newcommand{\tfour}{t_{\e,\eta}^{4}}
\newcommand{\tfourr}{t_{\e,\eta^2}^{4}}
\newcommand{\thfive}{\theta_{\e,\eta}^{5}}
\newcommand{\thstar}{\theta_{\e,\eta}^*}
\newcommand{\e}{\varepsilon}
\newcommand{\Om}{\Omega}
\newcommand{\Ga}{\Gamma}
\newcommand{\C}{{\mathbb C}}

\newcommand{\R}{{\mathbb R}}
\newcommand{\Rn}{{\R}^d}

\newcommand{\Mnn}{{\mathbb M}^{d{\times}d}_{sym}}
\newcommand{\MD}{{\mathbb M}^{d{\times}d}_D}

\newcommand{\wto}{\rightharpoonup}

\renewcommand{\div}{{\rm div}}
\renewcommand{\hom}{{hom}}
\newcommand{\supp}{{\rm supp}}
\newcommand{\tr}{{\rm tr}}
\newcommand{\hn}{{\mathcal H}^{d-1}}
\newcommand{\Ln}{{\mathcal L}^d}
\newcommand{\LL}{{\mathcal L}}

\newcommand{\QQ}{{\mathcal Q}}
\newcommand{\D}{{\mathcal D}}

\newcommand{\HH}{{\mathcal H}}
\newcommand{\K}{{\mathcal K}}

\newcommand{\bary}{{\rm bar}}
\newcommand{\smtr}{{_\triangle}}
\newcommand{\V}{{\mathcal V}}
\newcommand{\ee}{{\boldsymbol e}}
\newcommand{\pp}{{\boldsymbol p}}
\newcommand{\uu}{{\boldsymbol u}}
\newcommand{\vv}{{\boldsymbol v}}
\newcommand{\ww}{{\boldsymbol w}}
\newcommand{\zz}{{\boldsymbol z}}
\newcommand{\deltaa}{{\boldsymbol\delta}}
\newcommand{\zetaa}{{\boldsymbol\zeta}}

\newcommand{\muu}{{\boldsymbol\mu}}

\newcommand{\sigmaa}{{\boldsymbol\sigma}}
\newcommand{\psii}{{\boldsymbol\psi}}

\newcommand{\ol}{\overline}
\newcommand{\tki}{t_k^i}
\newcommand{\tkim}{t_k^{i-1}}
\newcommand{\wki}{w_k^i}
\newcommand{\uki}{u_k^i}
\newcommand{\eki}{e_k^i}
\newcommand{\pki}{p_k^i}
\newcommand{\zki}{z_k^i}

\title[Quasistatic evolution in plasticity with softening]
{A vanishing viscosity approach to quasistatic evolution in plasticity with softening}
\author{G.\ Dal Maso}
\author{A.\ DeSimone}
\author{M.G.\ Mora}
\author{M.\ Morini}
\address[G.~Dal Maso, A.~DeSimone, M.G.~Mora, and M.\ Morini]{SISSA, Via Beirut 2-4, 
34014 Trieste, Italy}
\email[Gianni Dal Maso]{dalmaso@sissa.it}
\email[Antonio DeSimone]{desimone@sissa.it}
\email[Maria Giovanna Mora]{mora@sissa.it}
\email[Massimiliano Morini]{morini@sissa.it}

\begin{document}
\begin{abstract}
We deal with quasistatic evolution problems in plasticity with softening, in the framework of
small strain associative elastoplasticity. The presence of a nonconvex term due to the
softening phenomenon requires a nontrivial  extension of the variational framework for
rate-independent problems to the case of a nonconvex energy functional.
We argue that, in this case, the use of global minimizers in the corresponding
incremental problems
is not justified from the mechanical point of view. Thus,  we analize a different selection
criterion
for the solutions of the quasistatic evolution problem, based on a viscous approximation.
This leads to a generalized formulation in terms of Young measures, developed in the
first part of the paper. In the second part we apply our approach to some concrete
examples.
\end{abstract}
\maketitle

{\small
\keywords{\noindent {\bf Keywords:} 
plasticity with softening, quasistatic evolution, rate-independent processes, shear bands, incremental problems, viscous approximation, Young measures. \blfootnote{Preprint SISSA 38/2006/M (June 2006)}}

}
\tableofcontents

\begin{section}{Introduction}

In plasticity theory the term {\em softening\/} refers to the reduction of the yield stress as
plastic deformation proceeds. Classically this is described by a family of yield surfaces
depending on a parameter $\zeta$. The evolution laws are formulated in such a way that the
yield surface shrinks when the time derivative of the plastic deformation is not zero
(see  \cite{Han-Red}, \cite{Hill}, \cite{Lub}, and \cite{Mar}).

We deal with this problem in the {\em quasistatic\/} case, in the framework of {\em small
strain associative elastoplasticity\/} in a bounded and Lipschitz domain $\Om\subset\Rn$, $d\ge 2$.
For simplicity we consider only the case of {\em no applied forces\/} and of {\em prescribed
boundary displacements\/} on a closed subset $\Ga_0$ of the boundary
$\partial\Om$ with positive ${(d-1)}$-dimensional measure.
The linearized strain $Eu$, defined as the symmetric part of the spatial
gradient of the displacement $u$, is decomposed as the sum $Eu=e+p$, 
where $e$ and $p$ are the elastic and plastic strains.
The stress $\sigma$ is determined only by $e$, through the formula $\sigma=\C e$, where
$\C$ is the elasticity tensor. We assume that, for every value of the parameter 
$\zeta$, the {\it elastic domain\/} -- the set of admissible stresses enclosed by the yield surface -- has the form
$\{\sigma\in\Mnn:\sigma_D\in K(\zeta)\}$, where $\Mnn$  is the space of symmetric
$d{\times}d$ matrices, $\sigma_D$ denotes the deviatoric part of $\sigma$, and
$K(\zeta)$ is a subset of the subspace $\MD$ of trace-free symmetric matrices.
To simplify the mathematics of the problem, we assume that the set
$$
K:=\{(\sigma,\zeta)\in \MD{\times}\R :\sigma\in K(\zeta)\}
$$ 
is a {\em compact convex neighbourhood\/} of $(0,0)$ in $\MD{\times}\R$.

To express the evolution laws, it is convenient to introduce an internal variable $z$,
related to $\zeta$ by the equation $\zeta=-V'(z)$,
where $V\colon \R\to\R$ is a given function of class $C^2$ with bounded second derivatives,
called the {\em softening potential\/}.

The strong formulation of the quasistatic evolution problem consists in finding functions
$\uu(t,x)$, $\ee(t,x)$, $\pp(t,x)$, $\sigmaa(t,x)$, $\zz(t,x)$, and $\zetaa(t,x)$ satisfying  
the following conditions for every $t\in[0,+\infty)$ and every $x\in\Om$:
\begin{itemize}
\smallskip
\item[(sf1)] additive decomposition: $E\uu(t,x)=\ee(t,x)+\pp(t,x)$,
\item[(sf2)] constitutive equations: $\sigmaa(t,x)=\C \ee(t,x)$ and $\zetaa(t,x)=-V'(\zz(t,x))$,
\item[(sf3)] equilibrium: $\div\,\sigmaa(t,x)=0$,
\item[(sf4)] stress constraint: $\sigmaa(t,x)_D\in K(\zetaa(t,x))$,
\item[(sf5)] associative flow rule: 
$(\dot \pp(t,x),\dot \zz(t,x))\in N_K(\sigmaa(t,x)_D,\zetaa(t,x))$,
\smallskip
\end{itemize}
where dots denote time derivatives and $N_K(\sigma,\zeta)$ is the normal cone to $K$ at
$(\sigma,\zeta)$ in $\MD{\times}\R$. The evolution is driven by a 
prescribed time-dependent boundary condition
$$
\uu(t,x)=\ww(t,x)\quad \hbox{for every } t\in [0,+\infty) \hbox{ and every } x\in\Ga_0\,.
$$
It is supplemented by initial conditions at $t=0$ and by the traction-free boundary condition 
$$
\sigmaa(t,x) n(x)=0\quad \hbox{for every } t\in [0,+\infty) \hbox{ and every } x\in\Ga_1\,,
$$
where $n(x)$ is the normal to $\partial\Omega$ at $x$ and 
$\Ga_1=:\partial\Om\setminus\Ga_0$.

Introducing  the support function
$$
H(\xi,\theta):=\sup_{(\sigma,\zeta)\in K} (\xi{\,:\,}\sigma + \theta\,\zeta) \,,
$$
where the colon denotes the scalar product between matrices,
the flow rule (sf5) can be written in the equivalent form
\begin{itemize}
\smallskip
\item[$\rm(sf5')$] dissipation pseudo-potential formulation:
$(\sigmaa(t,x)_D,\zetaa(t,x))\in \partial H(\dot \pp(t,x),\dot \zz(t,x))$,
\smallskip
\end{itemize}
where $\partial H(\xi,\theta)$ denotes the subdifferential of $H$ at $(\xi,\theta)$.

If $V$ is strictly convex, this model describes plasticity with hardening, where the
yield surface expands when the time derivative of the plastic strain is not zero. In this case 
it is possible to give a variational formulation of the problem, which does not assume 
the existence of time derivatives, and is based on the energies
\begin{eqnarray*}
&\displaystyle \QQ(e):={\textstyle\frac12} \int_\Om \C e(x){\,:\,}e(x)\,dx \,,
\qquad
\HH(p,z):= \int_\Om H(p(x),z(x))\,dx\,,\nonumber
\\
&
\displaystyle \V(z) :=\int_\Om V(z(x))\,dx  \,.
\end{eqnarray*}
The term $\HH(p,z)$ is used to introduce the notion of dissipation of a function 
$t\mapsto (\pp(t),\zz(t))$ on an interval $[a,b]\subset[0,+\infty)$, defined by
$$
\D_{\!H}(\pp,\zz;a,b):= 
 \sup \sum_{j=1}^k \HH(\pp(t_j)-\pp(t_{j-1}),\zz(t_j)-\zz(t_{j-1}))\,,
$$
where the supremum is taken over all finite sequences $(t_j)$ such that 
$a=t_0<t_1<\dots<t_{k-1}<t_k=b$.

According to the energetic approach to rate-independent processes developed in \cite{Mie-The-Lev}, \cite{Mai-Mie}, \cite{Mie-review}, the variational formulation of the quasistatic evolution problem consists in finding functions
$\uu(t,x)$, $\ee(t,x)$, $\pp(t,x)$, and $\zz(t,x)$ satisfying the following conditions:
\begin{itemize}
\smallskip
\item[(vf1)] global stability: for every $t\in [0,+\infty)$ we have $E\uu(t)=\ee(t)+\pp(t)$ on 
$\Om$, 
$\uu(t)=\ww(t)$ on $\Ga_0$, and
$$
\QQ(\ee(t))+\V(\zz(t))\le \QQ(\hat e)+ \HH(\hat p-\pp(t),\hat z-\zz(t))+ \V(\hat z)
$$
for every $\hat u$, $\hat e$, $\hat p$, $\hat z$ such that $E\hat u=\hat e+\hat p$ on $\Om$,
$\hat u=\ww(t)$ on $\Ga_0$;
\smallskip
\item[(vf2)] energy inequality: for every $T\in[0,+\infty)$ we have
$$
\QQ(\ee(T))+\D_{\!H}(\pp,\zz;0,T)+ \V(\zz(T))\le \QQ(\ee(0))+ \V(\zz(0))+
\int_{0}^{T}\langle\sigmaa(t),E\dot \ww(t)\rangle\,dt\,.
$$
\end{itemize}
In the convex case, thanks to the Euler conditions, (vf1) is equivalent to the following property:
\begin{itemize}
\smallskip
\item[$\rm(vf1')$] stability: for every $t\in [0,+\infty)$ we have $E\uu(t)=\ee(t)+\pp(t)$ on
$\Om$, $\uu(t)=\ww(t)$ on $\Ga_0$, and
$$
\begin{array}{c}
\div\,\sigmaa(t,x)=0\quad\hbox{for }x\in\Om\,,\quad 
\sigmaa(t,x)n(x)=0\quad\hbox{for }x\in \Ga_1\,,
\\
(\sigmaa(t,x)_D,\zetaa(t,x))\in K\quad\hbox{for }x\in\Om\,,
\end{array}
$$
where $\sigmaa$ and $\zetaa$ are defined by (sf2).
\smallskip
\end{itemize}

There is a vast literature on variational methods in the study of evolution problems in
elasto-plasticity. Among the papers which are closest in spirit to our approach we quote
 \cite{Ort-Mar}, \cite{Ort-Sta},
\cite{Miehe}, and \cite{Mie}.

In this paper we assume that $V$ is {\em concave\/}, which reflects the fact that the yield surface shrinks as $\dot \pp(t,x)\neq 0$. To simplify the mathematics of the problem, we also assume that the image of $-V'$ is contained in the interior of the projection of $K$ onto the $\zeta$-axis. By lack of convexity, condition $\rm(vf1')$ is no longer equivalent to (vf1), and we have only $\rm(vf1)\Longrightarrow\rm(vf1')$. In this case the selection criterion provided by global minimality is not justified from the mechanical point of view. Indeed, as we shall
show in \cite{DM-DeS-Mor-Mor-2}, global minimality leads to missing the softening
phenomenon altogether.

We explore a different selection criterion,  based on the approximation by solutions of some
regularized evolution problems, depending on a small {\em ``viscosity''\/}  parameter $\e>0$.
In its strong formulation this regularized problem consists in finding functions
$\uu_\e(t,x)$, $\ee_\e(t,x)$, $\pp_\e(t,x)$, $\sigmaa_\e(t,x)$, $\zz_\e(t,x)$, and $\zetaa_\e(t,x)$ satisfying (sf1), (sf2), (sf3), and
\begin{itemize}
\smallskip
\item[(sf4)$\!_\e$] regularized flow rule: 
$(\dot \pp_\e(t,x),\dot \zz_\e(t,x))= N_K^\e(\sigmaa_\e(t,x)_D,\zetaa_\e(t,x))$,
\smallskip
\end{itemize}
where 
$$
N_K^\e(\sigma,\zeta):=\tfrac{1}{\e}\big((\sigma,\zeta)-P_K(\sigma,\zeta)\big)\,,
$$
$P_K$ being the projection onto $K$. If $\zz$ and $\zetaa$ were not present, this condition would coincide with the flow rule of Perzyna viscoplasticity.
We observe that condition (sf4)$\!_\e$ is closely related to (sf5).
It  also acts as a penalization which leads to the stress constraint (sf4) as $\e\to 0$.

The existence of a solution to the $\e$-regularized evolution problem is proved by a
variational method based on time discretization and on the solution of
suitable incremental minimum problems (see Theorem~\ref{young-main2}).
Some parts of the proof are inspired by \cite{Suq}. The uniqueness is based on Gronwall's lemma. 
In Theorem~\ref{equivalence} we also prove that this solution is characterized by the following conditions:
\begin{itemize}
\smallskip
\item[(re1)$\!_\e$] equilibrium condition: for every $t\in [0,+\infty)$ we have
$E\uu_\e(t)=\ee_\e(t)+\pp_\e(t)$ in $\Om$, $\uu_\e(t)=\ww(t)$ on $\Ga_0$, and
$$
\begin{array}{c}
\div\,\sigmaa_\e(t,x)=0\quad\hbox{for }x\in\Om\,,\quad \sigmaa_\e(t,x)n(x)=0
\quad\hbox{for }x\in \Ga_1\,,
\\
(\sigmaa_\e(t,x)_D-\e\dot \pp_\e(t,x),\zetaa_\e(t,x)-\e \dot \zz_\e(t,x))\in K
\quad\hbox{for }x\in\Om\,,
\end{array}
$$
where $\sigmaa_\e$ and $\zetaa_\e$ are defined by $\sigmaa_\e(t,x):=\C \ee_\e(t,x)$ and
$\zetaa_\e(t,x):=-V'(\zz_\e(t,x))$;
\smallskip
\item[(re2)$\!_\e$] energy equality: for every $T\in[0,+\infty)$ we have
$$
\begin{array}{c}
\displaystyle
\QQ(\ee_\e(T))+\D_{\!H}(\pp_\e,\zz_\e;0,T)+ \V(\zz_\e(T))
+ \e\int_{0}^{T}\|\dot \pp_\e(t)\|_2^2\,dt
+\e\int_{0}^{T}\|\dot \zz_\e(t)\|_2^2\,dt
=
\smallskip
\\
\displaystyle
= \QQ(\ee_\e(0))+ \V(\zz_\e(0))
+\int_{0}^{T}\langle \sigmaa_\e(t),E\dot \ww(t)\rangle\,dt\,.
\end{array}
$$
\end{itemize}

By accepting only those solutions of $\rm(vf1')$ and (vf2) which can be approximated by
solutions of (sf1), (sf2), (sf3), (sf4)$\!_\e$ (or, equivalently, by solutions of (re1)$\!_\e$ and (re2)$\!_\e$), we regard quasistatic evolution as the limiting case of a viscosity-driven dynamics (Definition~\ref{def:qsym}). A similar approach in finite dimension was
used in \cite{Efe-Mie}. Other rate-independent problems with nonconvex energy have been studied  in \cite{Car-Hac-Mie},\cite{Kru-Mie-Rou}, and \cite{Fra-Mie}. The last paper considers a different regularizing term based on the space gradient of the internal variable.

The main difficulty in our approach is due to the fact that, 
by the nonconvexity of the energy, the components $\pp_\e$ and $\zz_\e$ of the
solutions of the $\e$-regularized problem may develop stronger and stronger space
oscillations as $\e\to 0$. As a consequence of this fact, their weak limits do not
satisfy, in general,  $\rm(vf1')$ and (vf2)  (see Section~\ref{conc-osc}). To overcome
this difficulty, we propose a weaker formulation in terms of Young measures. 

Since the functionals $\HH$ and $\V$ have linear growth, the classical notion of
Young measure is not enough. To take into account possible concentrations at
infinity, we use the notion of generalized Young measure introduced in 
\cite{Dip-Maj}, \cite{Ali-Bou}, and \cite{Fon-Mue-Ped}, following the presentation 
of~\cite{DM-DeS-Mor-Mor-1}.  

In addition, to write the Young measure version of (vf2) we need to introduce a notion of dissipation for a time-dependent family of generalized Young measures. A natural definition can be given by taking the limit of the dissipations of suitable time-dependent generating functions. Unfortunately, this limit does not depend only on the values of the generalized Young measures at each time, but it also involves the mutual correlations between oscillations at different times. We solve this problem by using the notion of system of generalized Young measures introduced in \cite{DM-DeS-Mor-Mor-1}. This allows us to write a Young measure formulation of  problem $\rm(vf1')$, (vf2) 
(Theorem~\ref{Thm55}) and to prove an existence result (Theorem~\ref{young-main3}).

The second part of the paper is devoted to some examples.
The first example, developed in Section~\ref{examples}, deals with a spatially homogeneous
case, where the $\e$-regularized evolution is described by a system of ordinary differential
equations, and the study of the limit as $\e\to 0$ reduces to the analysis of a singular perturbation problem for ordinary differential equations. In this example the stress
$\sigmaa(t)$ converges, as $t\to\infty$, to a constant matrix $\sigma_\infty$, 
which coincides with
the yield stress that would be obtained in the perfectly plastic case with elastic domain
$$
K_\infty:=\bigcap\nolimits_{\zeta} K(\zeta)\,.
$$
In other words, as time tends to infinity, the material behaves in the weakest
way permitted by its internal variable. Moreover, after a critical time, $|\sigmaa(t)|$ is decreasing with respect to $t$, reflecting the fact that the material softens as  plastic deformation proceeds (see Figure~1). Finally, for some values of the parameters, this example exhibits
a jump discontinuity at a certain time (see Figures~2 and~3). This leads to a strict inequality in the energy balance (vf2), which shows that an instantaneous dissipation occurs at the discontinuity time.

A second group of examples, developed in Section~\ref{conc-osc}, shows that 
strain localization, in the form of a shear band, may occur in this model even if the initial 
and boundary data are sufficiently regular.
One of the examples exhibits also a strong oscillation of the internal variable $\zz_\e$ localized near the shear band. As $\e\to0$, this leads to a Young measure solution of the quasistatic evolution problem. In this example the usual weak$^*$ limit of $\zz_\e$ (always given by the barycentre of the Young measure solution) still satisfies the equilibrium condition $\rm( vf1')$, but does not satisfy the energy inequality~(vf2). This shows that, because of the lack of convexity, some important terms generated by the space oscillations of the approximate solutions can be captured only by the Young measure formulation.

\end{section}

\begin{section}{Notation and preliminary results}

\subsection{Mathematical preliminaries}
We begin with a quick presentation of the mathematical tools used in the paper.

\medskip
\noindent
{\bf Measures.} 
The Lebesgue measure on $\Rn$, $d\ge 1$, is denoted by ${\mathcal L}^d$, and the 
$k$-dimen\-sional Hausdorff measure by ${\mathcal H}^k$. 
Given a Borel set $B\subset\Rn$ and a finite dimensional Hilbert space 
$\Xi$, $M_b(B;\Xi)$ denotes the space of 
bounded Borel measures on $B$  with values in $\Xi$, endowed with the norm 
$\|\mu\|_1:=|\mu|(B)$, where 
$|\mu|\in M_b(B):=M_b(B;\R)$ is the variation of the measure $\mu$. 
The space of nonnegative bounded Borel measures on $B$ is denoted by $M^+_b(B)$.
For every $\mu\in M_b(B;\Xi)$ we consider the Lebesgue decomposition 
$\mu=\mu^a+\mu^s$, where $\mu^a$ is 
absolutely continuous and $\mu^s$ is singular with respect to Lebesgue 
measure $\Ln$. 

If $\mu^s=0$, we always 
identify $\mu$ with its density with respect to Lebesgue measure ${\mathcal L}^d$. In 
this way
$L^1(B;\Xi)$ is regarded as a subspace of $M_b(B;\Xi)$, with the induced norm. 
In particular $\mu^a\in L^1(B;\Xi)$ 
for every $\mu\in M_b(B;\Xi)$.
The $L^r$ norm, $1\le r\le\infty$, is denoted 
by $\|\cdot\|_r$. The brackets $\langle \cdot,\cdot\rangle$ denote the duality product between conjugate $L^r$ spaces, as well as between other pairs of spaces, according to the context. The symbols $\lor$ and $\land$ denote the maximum and minimum of two numbers or functions, while the symbol $(\cdot)^+$ denotes the positive part.

If  $B$ is locally compact (in the relative topology), by the Riesz representation theorem 
(see, e.g., \cite[Theorem~6.19]{Rud})
$M_b(B;\Xi)$ can be identified with the dual of 
$C_0(B;\Xi)$, the space of continuous functions $\varphi\colon B\to \Xi$
such that $\{|\varphi|\ge \e\}$ is compact for every $\e>0$. The 
weak$^*$ topology of $M_b(B;\Xi)$ is defined using this duality.

\medskip

\noindent {\bf Matrices.}
The space of {\it symmetric $d{\times}d$ matrices\/} is denoted by $\Mnn$; it is endowed with the euclidean scalar product 
$\xi{\,:\,}\zeta:=\tr(\xi\zeta)=\sum_{ij}\xi_{ij}\zeta_{ij}$ and with the 
corresponding euclidean norm 
$|\xi|:= (\xi{\,:\,}\xi)^{1/2}$. The {\it symmetrized tensor product\/} $a{\,\odot\,}b$ of two vectors $a$, 
$b\in\Rn$ is the symmetric matrix with 
entries $(a_ib_j+a_jb_i)/2$. 
It is easy to see that $\tr(a{\,\odot\,}b)=a{\,\cdot\,}b$, the scalar 
product of $a$ and $b$.

When $d\ge 2$ we define $\MD$ as the space of all matrices of 
$\Mnn$ with trace zero. It turns out that $\MD$ is 
the orthogonal complement of the subspace $\R I$ spanned by the identity 
matrix $I$. For every $\xi\in\Mnn$ the orthogonal projection 
of $\xi$ on $\R I$ is $\frac1d \tr(\xi)I$, while the 
orthogonal projection on $\MD$ is the {\it deviator\/} $\xi_D$ of $\xi$, so that 
we have the orthogonal decomposition
$$
\textstyle
\xi=\xi_D+\frac1d (\tr\,\xi) I\,.
$$
The case $d=1$ is special: we have ${\mathbb M}^{1{\times}1}_{sym}=\R$ and we do 
not need any
orthogonal decomposition of the space. For the purposes of this paper 
(see Section~\ref{simple-shear}) it is convenient 
to define ${\mathbb M}^{1{\times}1}_D:=\R$ and 
$\xi_D:=\xi$ for every $\xi\in\R$, although this does not agree with the definition given for
$d\ge 2$.

\medskip

\noindent {\bf Functions with bounded deformation.}
Let $U$ be an open set in $\Rn$, $d\ge 1$. For every $u\in L^1(U;\Rn)$ let $Eu$ be the 
$\Mnn$-valued distribution on $U$, whose components are defined by $E_{ij}u= (D_j u_i +D_i u_j)/2$. The space $BD(U)$ of functions with {\it bounded deformation\/} is the space of all $u\in L^1(U;\Rn)$ such that $Eu\in M_b(U;\Mnn)$. It is easy to see that $BD(U)$ is a Banach space with the norm
$$
\|u\|_1+\|Eu\|_1\,.
$$
It is possible to prove that $BD(U)$ is the dual of a normed space
(see \cite{Mat} and \cite{Tem-Stra}). The weak$^*$ topology of $BD(U)$ is defined
using this duality. A sequence $u_k$ converges to $u$ weakly$^*$ in $BD(U)$ if and
only if $u_k\wto u$ weakly in $L^1(U;\Rn)$ and $Eu_k\wto Eu$ weakly$^*$ in 
$M_b(U;\Mnn)$. Every bounded sequence in $BD(U)$ has a weakly$^*$ convergent 
subsequence. Moreover, if $U$ is bounded and has Lipschitz boundary, every bounded sequence in $BD(U)$ has a subsequence which converges weakly in
$L^{d/(d-1)}(U;\Rn)$ and strongly in $L^r(U;\Rn)$ for every ${r<d/(d-1)}$.
For the general properties of $BD(U)$ we refer to \cite{Tem}.

In our problem $u\in BD(U)$ represents the {\it displacement\/} of an elasto-plastic body and $Eu$ is the corresponding linearized {\it strain\/}.

\medskip

\noindent {\bf Generalized Young measures.}
As mentioned in the introduction, the results of this paper are based on the notion of systems of generalized Young measures developed in~\cite{DM-DeS-Mor-Mor-1}.  For the reader's convenience we collect here the main definitions and our notational conventions, while we refer to \cite{DM-DeS-Mor-Mor-1} for the motivations behind them and for the main properties.

Let $U$ be a bounded open set in $\Rn$, $d \ge 1$. We shall apply the notions introduced in \cite{DM-DeS-Mor-Mor-1} with $X:=\ol U$ and $\lambda={\mathcal L}^d$.
To define the space of generalized Young measures we introduce the space
$C^\hom(\ol U{\times}\Xi)$ of all
$f\in C(\ol U{\times}\Xi)$ such that  $f(x,\cdot)$ is positively homogeneous of degree
one on $\Xi$ for every $x\in\ol U$. This space is endowed with the norm
$$
\|f\|_{\hom}:= \max\{ |f(x,\xi)|: x\in \ol U,\ \xi\in\Sigma_\Xi\} \,,
$$
where $\Sigma_\Xi:=\{\xi\in\Xi: |\xi|=1\}$. The dual of the Banach space 
$C^\hom(\ol U{\times}\Xi)$ is denoted by $M_*(\ol U{\times}\Xi)$, and the corresponding dual norm by $\|\cdot\|_*$; the weak$^*$ topology of 
$M_*(\ol U{\times}\Xi)$ is defined by using this duality. As in the case of $M_b(\ol U{\times} \Xi)$, it is sometimes convenient to write the dummy variables explicitly also for this duality product and to use the notation $\langle f(x,\xi),\mu(x,\xi)\rangle$ instead of $\langle f,\mu\rangle$.

The space $\Xi{\times} \R$ is endowed with the product Hilbert structure. The corresponding spaces $C^\hom(\ol U{\times}(\Xi{\times} \R))$ and $ M_*(\ol U{\times}(\Xi{\times} \R))$ are denoted by $C^\hom(\ol U{\times}\Xi{\times} \R)$ and $ M_*(\ol U{\times}\Xi{\times} \R)$.
The space $GY(\ol U;\Xi)$ of {\it generalized Young measures\/} on $\ol U$ with values in $\Xi$ is defined as the set of all $\mu\in M_*(\ol U{\times}\Xi{\times} \R)$ satisfying the following properties:
\begin{itemize}
\item[(a)] positivity property: 
\begin{equation}\label{mupos}
\langle f,\mu\rangle\ge 0\quad\text{for every }f\in
C^\hom(\ol U{\times}\Xi{\times} \R)\text{ with }f\ge 0\,;
\end{equation}
\item[(b)] support property:
$$
\langle f,\mu\rangle= 0 \quad\text{for every }f\in
C^\hom(\ol U{\times}\Xi{\times} \R) \text{ vanishing on }
\ol U{\times}\Xi{\times} {[0,{+\infty})}\,;
$$
\item[(c)] projection property: 
\begin{equation}\label{muproj}
\displaystyle \langle\varphi(x)\eta, \mu(x,\xi,\eta)\rangle=\int_U \varphi(x) \, dx \qquad\text{for every }\varphi\in
C(\ol U)\,.
\end{equation}
\end{itemize}
It is easy to see that for every $\mu\in GY(\ol U;\Xi)$ we have
\begin{equation}\label{normpositive}
\|\mu\|_*=\langle \sqrt{|\xi|^2 + |\eta|^2},\mu(x,\xi,\eta)\rangle\,.
\end{equation}

We recall that for any $\mu\in GY(\ol U;\Xi)$ the {\em barycentre\/} of $\mu$, denoted by 
$\bary(\mu)$, is defined as the element of $M_b(\ol U;\Xi)$ such that 
\begin{equation}\label{2bary}
\langle\varphi,\bary(\mu)\rangle=\langle\varphi(x) {\,\cdot\,}\xi,
\mu(x,\xi,\eta)\rangle
\end{equation}
for every $\varphi\in C(\ol U;\Xi)$. 
If $\mu_k\wto\mu$ weakly$^*$ in $GY(\ol U;\Xi)$, then 
\begin{equation}\label{weakbar}
\bary(\mu_k)\wto \bary(\mu)\quad \hbox{weakly}^*\hbox{ in }M_b(\ol U;\Xi)\,. 
\end{equation}

Given a measure $p\in M_b(\ol U;\Xi)$, the {\it generalized Young measure
associated with $p$\/} is the element $\delta_p$ of $GY(\ol U;\Xi)$ defined for
every $f\in C^\hom(\ol U{\times}\Xi{\times}\R)$ by 
\begin{equation}\label{mup}
\langle f, \delta_p \rangle := \int_ {\ol U}
\textstyle f(x, \frac{dp}{d\lambda}(x), \frac{d{\mathcal L}^d}{d\lambda}(x))\,d\lambda(x)\,,
\end{equation}
where $\lambda\in M_b^+(\ol U)$ is an arbitrary measure such that 
${\mathcal L}^d<<\lambda$ and
$p<<\lambda$ (the homogeneity of $f$ implies that the integral does not depend on 
$\lambda$). 



Given another finite dimensional Hilbert space $\Xi'$, let 
$\psi\colon \ol U{\times}\Xi{\times}\R \to \ol U{\times}\Xi'{\times}\R$ be a continuous 
map of the form $\psi(x,\xi,\eta)=(x,\phi(x,\xi,\eta),\eta)$, with $\phi(x,\xi,\eta)$
positively one-homogeneous in $(\xi,\eta)$. The image $\psi(\mu)$ of any
$\mu\in GY(\ol U;\Xi)$ is defined as the element of $GY(\ol U;\Xi')$ such that
\begin{equation}\label{image2}
\langle f,\psi(\mu) \rangle :=\langle f\circ\psi  ,\mu\rangle
=\langle f(x,\phi(x,\xi,\eta),\eta),\mu(x,\xi,\eta)\rangle
\end{equation}
for every $f\in C^{hom}(\ol U{\times}\Xi'{\times}\R)$.
In \cite{DM-DeS-Mor-Mor-1} it is proved that formula (\ref{image2}) makes sense also if $\phi$ is a Borel map, positively one-homogeneous in $(\xi,\eta)$, and satisfying $|\phi(x,\xi,\eta)|\le a|\xi| + b(x)|\eta|$
with $a\in\R$ and $b\in L^1(U)$.  

For additional properties of the space $GY(\ol U;\Xi)$ we refer to 
\cite{DM-DeS-Mor-Mor-1}. 

\medskip

\noindent {\bf Systems of generalized Young measures.}
Let $U$ be a bounded open set in $\R^n$, $n\ge 1$, and let $\Xi$ be a finite dimensional 
Hilbert space.
If $\{s_1,s_2,\dots,s_n\}\subset\{t_1,t_2,\dots,t_m\}\subset \R$, with 
$s_1<s_2<\dots<s_n$ and $t_1 < t_2 < \dots < t_m$, the projection 
$\pi_{s_1\dots s_n}^{t_1\dots t_m}\colon 
\ol U{\times}\Xi^m{\times}\R\to \ol U{\times}\Xi^n{\times}\R$ is defined by 
$$
\pi_{s_1\dots s_n}^{t_1\dots t_m}(x,\xi_{t_1},\dots,\xi_{t_m},\eta)=
(x,\xi_{s_1},\dots,\xi_{s_n},\eta)\,.
$$
A {\it compatible system of generalized Young measures\/} on $\ol U$ with values in
$\Xi$ with time set $\Theta\subset\R$ is a family $\muu=(\muu_{t_1\dots t_m})$ of
generalized Young measures $\muu_{t_1\dots t_m}\in GY(\ol U;\Xi^m)$, with 
$t_1, \dots,t_m$ running over all finite sequences of elements of $\Theta$ with
$t_1 < t_2 < \dots < t_m$, such that the following compatibility condition holds:
\begin{equation}\label{compatib}
\muu_{s_1\dots s_n}=\pi_{s_1\dots s_n}^{t_1\dots t_m}(\muu_{t_1\dots t_m})
\quad \text{whenever } \{s_1,s_2,\dots,s_n\}\subset\{t_1,t_2,\dots,t_m\}\,.
\end{equation}
The space of all such systems is denoted by $SGY(\Theta,\ol U;\Xi)$.

Given a function $\pp\colon \Theta\to M_b(\ol U;\Xi)$,  the {\em compatible system of generalized Young measures $\deltaa_\pp$ associated with\/} $\pp$ is the element of $SGY(\Theta,\ol U;\Xi)$ defined by
\begin{equation}\label{delta1m}
(\deltaa_\pp)_{t_1\dots t_m}:=\delta_{(\pp(t_1),\dots,\pp(t_m))}\,,
\end{equation}
where the right-hand side is defined by \eqref{mup}. 
 
Finally, we recall that the {\em variation\/} of $\muu\in SGY(\Theta,\ol U;\Xi)$ on a time interval $[a,b]$, with $a,\, b\in \Theta$, is defined as 
\begin{equation}\label{nuvar}
{\rm Var}(\muu;a,b):=
\sup \sum_{i=1}^k 
\langle |\xi_i-\xi_{i-1}|, \muu_{t_0t_1\dots t_k}(x,\xi_0,\dots,\xi_k,\eta) \rangle\,,
\end{equation}
where the supremum is taken over all finite families $t_0,t_1,\dots,t_k$ in $\Theta$ such that $a=t_0<t_1<\dots<t_k=b$.
If $\muu=\deltaa_\pp$ for some $\pp\colon \Theta\to M_b(\ol U;\Xi)$, then ${\rm Var}(\muu;a,b)$ coincides with 
\begin{equation}\label{pvar}
{\rm Var}(\pp;a,b):=
\sup \sum_{i=1}^k 
\|\pp(t_i)-\pp(t_{i\!-\!1})\|_1\,,
\end{equation}
where the supremum is taken over all finite families $t_0,t_1,\dots,t_k$ in $\Theta$ such that $a=t_0<t_1<\dots<t_k=b$.

For the main properties of the space $SGY(\Theta,\ol U;\Xi)$ we refer to \cite[Section~7]{DM-DeS-Mor-Mor-1}. 

\subsection{Mechanical preliminaries}
We now introduce the mechanical notions used in the paper.

\medskip
\noindent
{\bf The reference configuration.}
Throughout the paper the {\em reference configuration\/} $\Om$ is a {\it bounded connected open set\/} in $\Rn$, 
$d\ge 1$, with 
{\it Lipschitz boundary\/} $\partial\Om=\Gamma_0\cup\Gamma_1$. We assume that 
$\Gamma_0$ is closed, $\hn(\Gamma_0)>0$, and  $\Gamma_0\cap\Gamma_1=\emptyset$.

On $\Ga_0$ we will prescribe a Dirichlet boundary condition.
This will be done by assigning a sufficiently regular function  $w \colon \Ga_0 \to \Rn$, or, equivalently, a function $w\in H^1(\Om;\Rn)$, whose trace on $\Ga_0$
(also denoted by $w$) is the prescribed boundary value. 

Every function $u\in BD(\Om)$ has a {\it trace\/} on $\partial\Om$, still denoted by $u$, which belongs to $L^1(\partial\Om;\Rn)$. Moreover, there exists a constant $C>0$, 
depending on $\Om$ and $\Ga_0$, such that
\begin{equation}
\label{seminorm}
\|u\|_{1,\Om} \le C\, \|u\|_{1,\Ga_0}+ 
C\, \|Eu\|_{1,\Om}
\end{equation}
(see \cite[Proposition~2.4 and Remark~2.5]{Tem}).

We will study two alternative situations.
\begin{itemize}
\item[(DN)] {\it Dirichlet-Neumann\/}. In this case we will consider 
traction-free boundary conditions on $\Ga_1$.

\item[(DP)] {\it Dirichlet-Periodic\/}. In this case $\Om$ is the cube $Q:=(-\frac12,\frac12)^d$, $d \ge 2$, and
\begin{equation}\label{Gper}
\Gamma_0:=\{ x\in \partial Q : |x_1| =\tfrac12\}\,, \quad 
\Gamma_1:=\{ x\in \partial Q : |x_1| <\tfrac12\}\,,
\end{equation}
where $x=(x_1,\ol x)$,  $\ol x=(x_2,\dots,x_d)$, and we consider only functions $u$ which are $\ol x$-periodic, in the sense that $u(x+e_i)=u(x)$ for every $x\in \Rn$ and every element $e_i$ of the canonical basis of $\Rn$, with $i=2,\dots, d$.
\end{itemize}

\medskip

\noindent {\bf Admissible stresses and dissipation.}
Let $K$ be a closed convex set in $\MD{\times}\R$, which will play the role of a constraint on the deviatoric part of the stress and on the internal variable~$\zeta$. For every $\zeta\in \R$
the set 
\begin{equation}\label{Kzeta}
K(\zeta):=\{\sigma\in \MD:(\sigma,\zeta)\in K\}
\end{equation}
 is interpreted as the {\it elastic domain\/} and its boundary as the {\it yield surface\/} corresponding to~$\zeta$. We assume that there exist two constants $A$ and $B$, with $0<A\le B<\infty$, such that
\begin{equation}\label{rk}
\{(\sigma,\zeta)\in\MD{\times}\R: |\sigma|^2+|\zeta|^2\le A^2\}\subset K\subset \{(\sigma,\zeta)\in\MD{\times}\R: |\sigma|^2+|\zeta|^2\le B^2\}\,.
\end{equation}
We assume in addition that 
\begin{equation}\label{0zeta}
(\sigma,\zeta)\in K \quad \Longrightarrow \quad (0,\zeta)\in K\,.
\end{equation}
Let $\pi_\R \colon \MD {\times} \R \to \R$ be the projection onto $\R$.
The hypotheses on $K$ imply that there exist two constants $a_K>0$ and $b_K>0$ such that 
\begin{equation}\label{prK}
\pi_\R(K) = [-a_K,b_K]\,.
\end{equation}

The {\it support function\/} 
$H\colon\MD{\times}\R\to {[0,+\infty)}$ of $K$, defined by
\begin{equation}\label{HD}
H(\xi,\theta):=\sup_{(\sigma,\zeta)\in K} \{\sigma{\,:\,}\xi+\zeta\,\theta\} \,,
\end{equation}
will play the role of the {\it dissipation density\/}.
It turns out that $H$ is convex and positively homogeneous of degree one on
$\MD{\times}\R$. In particular it satisfies the triangle inequality
$$
H(\xi_1+\xi_2,\theta_1+\theta_2)\le H(\xi_1,\theta_1)+H(\xi_2,\theta_2)\,.
$$
{}From (\ref{rk}) it follows that
\begin{equation}\label{boundsH}
A\sqrt{|\xi|^2+\theta^2}\le H(\xi,\theta)\le B\sqrt{|\xi|^2+\theta^2}\,,
\end{equation}
for every $(\xi,\theta)\in \MD{\times}\R$, and from (\ref{0zeta}) we obtain
\begin{equation}\label{H0zeta}
H(\xi,\theta)\geq H(0,\theta) = \begin{cases}
-a_K \theta & \text{if } \theta\le 0\,, \\
b_K \theta & \text{if } \theta\geq 0
\end{cases}
\end{equation}
for every $(\xi,\theta)\in \MD{\times}\R$.

Using the theory of convex functions of measures developed in \cite{Gof-Ser},
we introduce the functional $\HH\colon M_b(\ol\Om;\MD){\times}M_b(\ol\Om)\to\R$
defined by
$$
\HH(p,z):=\int_{\ol\Om} 
H({\textstyle\frac{dp}{d\lambda}(x),\frac{dz}{d\lambda}}(x))\,d\lambda(x) \,,
$$
where $\lambda\in M_b^+(\ol\Om)$ is any measure such that 
$p<<\lambda$ and
$z<<\lambda$ (the homogeneity of $H$ implies that the integral does not depend
on~$\lambda$).
Using \cite[Theorem~4]{Gof-Ser} and  \cite[Chapter~II, Lemma~5.2]{Tem} we can see that $\HH(p,z)$ coincides with the integral over $\ol\Om$ of the measure studied in 
\cite[Chapter~II, Section~4]{Tem}, hence
$\HH$ is lower semicontinuous on $M_b(\ol\Om;\MD){\times}M_b(\ol\Om)$ with respect to weak$^*$ convergence of measures.
It follows from the properties of $H$ that $\HH$ satisfies the triangle inequality, i.e., 
\begin{equation}\label{triangle}
\HH(p_1+p_2,z_1+z_2)\le \HH(p_1,z_1)+\HH(p_2,z_2)
\end{equation}
for every $p_1,p_2\in M_b(\ol\Om;\MD)$ and every $z_1,z_2\in M_b(\ol\Om)$.

For every $\e>0$ we introduce the function $H_\e\colon\MD{\times}\R\to\R$ defined as
\begin{equation}\label{He}
H_\e(\xi,\theta):=H(\xi,\theta)+\tfrac{\e}{2}|\xi|^2+\tfrac{\e}{2}|\theta|^2\,,
\end{equation}
and the corresponding integral functional
$\HH_\e\colon L^2(\Om;\MD){\times}L^2(\Om)\to\R$ defined by
$$
\HH_\e(p,z):=\int_{\Om}H_\e(p(x),z(x))\, dx\,.
$$
The convex conjugate $H_\e^*\colon\MD{\times}\R\to\R$ of $H_\e$ is defined by
$$
H_\e^*(\sigma,\zeta):=\sup_{(\xi,\theta)\in \MD{\times}\R}
\{\sigma{\,:\,}\xi+\zeta\,\theta-H_\e(\xi,\theta)\}\,.
$$
Since the convex conjugate $H^*$ of $H$ satisfies $H^*(\sigma,\zeta)=0$ 
for $(\sigma,\zeta)\in K$ and $H^*(\sigma,\zeta)=+\infty$ for 
$(\sigma,\zeta)\not\in K$ (see \cite[Theorem~13.2]{Roc}),
using  \cite[Theorem~16.4]{Roc} one can prove that
\begin{equation}\label{He*}
H_\e^*(\sigma,\zeta)=\tfrac{1}{2\e}|(\sigma,\zeta)-P_K(\sigma,\zeta)|^2\,,
\end{equation}
where $P_K\colon \MD{\times}\R\to K$ is the projection onto $K$. This implies that
$H_\e^*$ is differentiable, and that its gradient is given by
\begin{equation}\label{partialHe*}
N_K^\e(\sigma,\zeta):=\tfrac{1}{\e}\big((\sigma,\zeta)-P_K(\sigma,\zeta)\big)\,.
\end{equation}
Note that $N_K^\e$ is Lipschitz continuous.

Let $\HH_\e^*\colon L^2(\Om;\MD){\times}L^2(\Om)\to\R$ be the convex conjugate of $\HH_\e$.  By a general property of integral functionals 
(see, e.g., \cite[Proposition~IX.2.1]{Eke-Tem}) we have
$$
\HH_\e^*(\sigma,\zeta)=\int_\Om H_\e^*(\sigma(x),\zeta(x))\, dx\,,
$$
so that, by the Dominated Convergence Theorem, its gradient $\partial \HH_\e^*$ is given by
\begin{equation}\label{dHeNe}
\partial \HH_\e^*(\sigma,\zeta)(x)=N_K^\e(\sigma(x),\zeta(x))
\,, \qquad \hbox{for a.e. } x\in \Om\,.
\end{equation}
Therefore $\partial \HH_\e^*$ is Lipschitz continuous.

\medskip

\noindent {\bf The elasticity tensor.}
Let $\C$ be the {\it elasticity tensor\/}, considered as a symmetric 
positive definite linear operator $\C\colon\Mnn\to\Mnn$. 
We assume that the orthogonal subspaces $\MD$ and $\R I$ are invariant under~$\C$.
This is equivalent to saying that there exist a 
symmetric positive definite linear 
operator $\C_D\colon\MD\to\MD$ and a constant
$\kappa>0$ such that
\begin{equation}
\C\xi:= \C_D\xi_D+\kappa (\tr\,\xi) I 
\end{equation}
for every $\xi\in\Mnn$. Note that when $\C$ is isotropic, we have
\begin{equation}\label{iso}
\C\xi=2\mu \xi_D + \kappa (\tr\xi)I\,,
\end{equation}
where $\mu>0$ is the {\it shear modulus\/} and $\kappa>0$ is the 
{\it modulus of compression\/}, so that our assumptions are satisfied. 

Let $Q\colon\Mnn\to{[0,+\infty)}$ be the quadratic form associated with $\C$, defined by
\begin{equation}\label{W}
\textstyle Q(\xi):=\frac12 \C\xi{\,:\,}\xi=
\frac12\C_D\xi_D{\,:\,}\xi_D+\frac{\kappa}2(\tr\,\xi)^2\,.
\end{equation}
It turns out that there exist two constants $\alpha_\C$ and $\beta_\C$, 
with $0<\alpha_\C\le\beta_\C<+\infty$, such that
\begin{equation}\label{boundsC}
\alpha_\C|\xi|^2\le Q(\xi)\le 
\beta_\C|\xi|^2
\end{equation}
for every $\xi\in\Mnn$.
These inequalities imply 
\begin{equation}\label{normC}
|\C\xi|\le 2\beta_\C| \xi| \,.
\end{equation}

\medskip

\noindent {\bf The softening potential.}
Let $V\colon\R\to\R$ be a function of class $C^2$, which will control the evolution of the internal variable $\zeta$, and consequently of the set $K(\zeta)$ of admissible stresses. We assume that there exists a constant $M\ge 0$ such that 
\begin{eqnarray}
& -M\le V''(\theta)\le 0 \qquad\hbox{for every } \theta\,, 
\label{2der}
\\
& -b_K<V'(+\infty)\le V'(-\infty)<a_K\,,
\label{V'infty}
\end{eqnarray}
where $V'(\pm\infty)$ denote the limits of $V'(\theta)$ as $\theta\to\pm\infty$, 
while $a_K$ and $b_K$ are the constants in~(\ref{prK}).
We denote the {\it recession function\/} of $V$  by $V^\infty$; it is defined as
\begin{equation}\label{Vinf}
V^\infty(\theta):=\lim_{t\to+\infty}\frac{V(t\theta)}{t}=
\begin{cases}
V'(-\infty)\theta & \text{if } \theta\le 0\,, \\
V'(+\infty)\theta & \text{if } \theta\ge 0\,.
\end{cases}
\end{equation}
Let us fix $\e\in(0,1)$ such that
\begin{equation}\label{lambdaV'}
V'(-\infty)<(1-\e) a_K \qquad \hbox{and}\qquad
V'(+\infty)> -(1-\e) b_K \,.
\end{equation}
Then, by (\ref{H0zeta}),
$$
\begin{array}{c}
H(\xi_2-\xi_1,\theta_2-\theta_1)+V(\theta_2)-V(\theta_1)\ge
\\ \ge \e H(\xi_2-\xi_1,\theta_2-\theta_1)+(1-\e)H(0,\theta_2-\theta_1) +V'(\theta_3)(\theta_2 -\theta_1)
\end{array}
$$
for a suitable $\theta_3$. By (\ref{H0zeta}) and (\ref{lambdaV'}), using the monotonicity of $V'$, we obtain
$(1-\e)H(0,\theta_2-\theta_1) +V'(\theta_3)(\theta_2 -\theta_1)\ge 0$, hence
$$
H(\xi_2-\xi_1,\theta_2-\theta_1)+V(\theta_2)-V(\theta_1)\ge  \e H(\xi_2-\xi_1,\theta_2-\theta_1)\,.
$$
Therefore (\ref{boundsH}) implies that there exists a constant $  C^K_V>0$ such that
\begin{equation}\label{gammaM}
H(\xi_2-\xi_1,\theta_2-\theta_1)+V(\theta_2)-V(\theta_1)\ge  C^K_V |\xi_2-\xi_1| +  C^K_V |\theta_2-\theta_1|
\end{equation}
for every $\xi_1,\xi_2\in \MD$ and every $\theta_1,\theta_2\in\R$.

It is convenient to introduce the function $\V\colon L^1(\Om)\to\R$ defined by
$$
\V(z):=\int_\Om V(z(x))\,dx
$$
for every $z\in L^1(\Om)$.

\medskip

\noindent {\bf The prescribed boundary displacements.}
For every $t\in[0,+\infty)$ we prescribe a {\it boundary displacement\/} $\ww(t)$ in the
space $H^1(\Om;\Rn)$. This choice is motivated by the fact that we do not want to
impose ``discontinuous'' boundary data, so that, if the displacement develops sharp
discontinuities, this is due to energy minimization.

We assume also that $\ww\in H^1_{loc}([0,+\infty); H^1(\Om;\Rn))$, which means, by definition, that $\ww\in H^1([0,T]; H^1(\Om;\Rn))$ for every $T>0$, so that the time derivative $\dot \ww$ belongs to 
$L^2([0,T]; H^1(\Om;\Rn))$ and its strain $E\dot \ww$ belongs to $L^2([0,T];L^2(\Om;\Mnn))$. In the Dirichlet-Periodic 
case (DP) we assume also that $\ww(t)$ is (the restriction to $Q$ of) an $\ol x$-periodic function.
For the main properties of Sobolev functions with values in 
reflexive Banach spaces we refer to \cite[Appendix]{Bre}.

\medskip
\noindent
{\bf Elastic and plastic strains.}
Given a displacement $u\in BD(\Om)$ and a 
boundary datum $w\in H^1(\Om;\Rn)$, the {\it elastic strain\/}
$e\in L^2(\Om;\Mnn)$ and the {\it plastic strain\/} $p\in M_b(\ol\Om;\MD)$
satisfy the {\it weak kinematic admissibility conditions\/}
\begin{eqnarray}
& Eu=e+p \quad \hbox{in }\Om\,, \label{900}
\\
& p=(w-u){\,\odot\,}n\,\hn \quad \hbox{on } \Ga_0\,,  \label{901}
\end{eqnarray}
where $n$ denotes the outward unit normal. The condition on $\Ga_0$ shows, in particular, that the prescribed boundary condition $w$ is not attained on $\Ga_0$ whenever a plastic slip occurs at the boundary.
It follows from \eqref{900} and \eqref{901} that $e=E^au-p^a$ a.e.\ in $\Om$ and $p^s=E^su$ in $\Om$. Since $\tr\,p=0$, 
it follows from (\ref{900}) that $\div\,u=\tr\, e\in L^2(\Om)$
and from (\ref{901}) that $(w-u){\,\cdot\,}n=0$ $\hn$-a.e.\ on~$\Ga_0$. 
This shows that, if $u\in H^1(\Om;\Rn)$ and  $p\in L^2(\Om;\MD)$, a stronger kinematic admissibility condition is satisfied.

\begin{definition}\label{ADw}
Given $w\in H^1(\Om;\Rn)$, the set $A(w)$ of {\it admissible displacements and strains\/} for the boundary datum $w$ on $\Ga_0$ is defined as the set of all triples $(u,e,p)$, with $u\in H^1(\Om;\Rn)$, $e\in L^2(\Om;\Mnn)$, $p\in L^2(\Om;\MD)$, which satisfy the 
 {\it kinematic admissibility conditions\/}
\begin{eqnarray}
& Eu=e+p \quad \hbox{a.e.\ in }\Om\,, \label{900R}
\\
& u=w\quad \hn\hbox{-a.e.\ on } \Ga_0\,.  \label{901R}
\end{eqnarray}
In the Dirichlet-Periodic case (DP) we always assume that $w$ is the restriction to 
$Q$ of an 
$\ol x$-periodic function belonging to $H^1_{loc}(\Rn;\Rn)$ and we add the requirement that  the same property holds for $u$.
\end{definition}

\medskip

\noindent{\bf The stress.}
The {\it stress\/} $\sigma\in L^2(\Om;\Mnn)$ is given by
\begin{equation}\label{sigma}
\sigma:=\C e= \C_D e_D+\kappa\,\tr\, e\,,
\end{equation}
and the {\it stored elastic energy\/} by
\begin{equation}\label{elenergy}
\QQ(e)=\int_{\Om}Q(e(x))\, dx={\textstyle\frac12}\langle \sigma, e\rangle \,.
\end{equation}
It is well known that $\QQ$ is lower semicontinuous on 
$L^2(\Om;\Mnn)$ with respect to weak convergence.

Let $\QQ^*\colon L^2(\Om;\Mnn)\to[0,+\infty)$ be the convex conjugate of $\QQ$. It is well known that 
$$
\QQ^*(\sigma)=\int_\Om Q^*(\sigma(x))\, dx\,,
$$
where $Q^*\colon\Mnn\to\R$ is the convex conjugate of $Q$ given by
\begin{equation}\label{Q*}
Q^*(\sigma)=\tfrac12 \sigma{\,:\,}\C^{-1}\sigma\,.
\end{equation}

If $\sigma\in L^2(\Om;\Mnn)$ and $\div\,\sigma\in L^2(\Om;\Rn)$, then the trace of the normal component of $\sigma$ on $\partial\Om$, denoted by $[\sigma n]$, is  defined as the distribution 
 on $\partial\Om$ such that
\begin{equation}\label{sigmanu}
\langle [\sigma n],\psi\rangle_{\partial\Om}:=\langle \div\, \sigma,\psi 
\rangle + \langle \sigma,E\psi\rangle \qquad 
\end{equation}
for every $\psi\in H^1(\Om;\Rn)$.
It turns out that $[\sigma n]\in H^{-1/2}(\partial\Om;\Rn)$ (see, e.g., 
\cite[Theorem~1.2, Chapter~I]{Tem}). We say that $[\sigma n]=0$ on $\Ga_1$ if 
$\langle [\sigma n],\psi\rangle_{\partial\Om}=0$ for every  $\psi\in H^1(\Om;\Rn)$ with $\psi=0$ $\hn$-a.e. on $\Ga_0$. 

\end{section}

\begin{section}{The minimum problem}

In this section we study in detail the incremental minimum problems used in the discrete-time formulation of the regularized evolution. 
The data are the current values $p_0\in L^2(\Om;\MD)$ and $z_0\in L^2(\Om)$ of the plastic strain
and of the internal variable and the updated value $w_1\in  H^1(\Om;\Rn)$ of the boundary displacement. By solving the minimum problem
\begin{equation}\label{minp0}
\min \{\QQ(e)+\HH(p-p_0,z-z_0)+\V(z)
+ \tfrac{\e}{2\tau}\|p-p_0\|_2^2 + \tfrac{\e}{2\tau}\|z-z_0\|_2^2  \}
\end{equation}
over the set of all $(u,e,p)\in A(w_1)$ and all $z\in L^2(\Om)$,
we get the updated values $u$, $e$, $p$, and $z$ of displacement, elastic strain, plastic strain,
and internal variable. In (\ref{minp0}) $\e> 0$ is a prescribed ``viscosity'' parameter, while $\tau>0$ will play the role of the time discretization step.


We are now in a position to prove the existence of a solution to~(\ref{minp0}), provided $\tau<\e/M$, where $M$ is the constant appearing in (\ref{2der}).

\begin{theorem}\label{existencemin}
Let $\e> 0$, $\tau\in {(0,\e/M)}$, $w_1\in H^1(\Om;\Rn)$, $p_0\in L^2(\Om;\MD)$, and $z_0\in L^2(\Om)$.
Then the minimum problem (\ref{minp0}) has a unique solution.
\end{theorem}

\begin{proof} 
Let $(u_k,e_k,p_k,z_k)$ be a minimizing sequence. 
By (\ref{gammaM}) we have
$$
\HH(p_k-p_0,z_k-z_0)+\V(z_k) \ge   C^K_V \|p_k-p_0\|_1
+  C^K_V \|z_k-z_0\|_1+\V(z_0)\,;
$$
therefore, \eqref{minp0} implies that
$\QQ(e_k)
+ \tfrac{\e}{2\tau}\|p_k-p_0\|_2^2 + \tfrac{\e}{2\tau}\|z_k-z_0\|_2^2$
is bounded. Consequently, the sequences $p_k$ and $z_k$ are bounded in
$L^2(\Om;\MD)$ and
in $L^2(\Om)$, respectively, and $\QQ(e_k)$ is bounded. Using (\ref{boundsC}) we deduce that the sequence $e_k$ is bounded in $L^2(\Om;\Mnn)$.
Since $Eu_k=e_k+p_k$ a.e.\ in $\Om$, it follows that $Eu_k$ is bounded in $L^2(\Om;\Mnn)$. 
Since $u_k=w_1$ $\hn$-a.e\ on $\Ga_0$, $u_k$ is bounded in $H^1(\Om;\Rn)$ by the Korn-Poincar\'e inequality (see e.g. \cite[Theorem 6.3-4]{Cia}). Since $\tau<\e/M$, the functional in (\ref{minp0}) is convex, hence weakly lower semicontinuous. The existence of a minimizer follows now from the direct methods of the calculus of variations. The uniqueness is a consequence of the strict convexity of the functional.
\end{proof}

We now derive the Euler conditions for a minimizer of (\ref{minp0}). 
In what follows, the symbol $\partial$ denotes the subdifferential of convex analysis
( see, e.g., \cite[Definition I.5.1]{Eke-Tem}).

\begin{theorem}\label{Euler}
Let $\e>0$, $\tau\in (0, \e/M )$, $p_0\in L^2(\Om; \MD)$, $\zeta_0 \in L^2(\Om)$, $w_1\in H^1(\Om;\Rn)$, $(u_1,e_1,p_1)\in A(w_1)$, $z_1\in L^2(\Om)$,
$\sigma_1:=\C e_1\in L^2(\Om;\Mnn)$, and $\zeta_1:= -V'(z_1)\in L^\infty(\Om)$.
If $(u_1,e_1,p_1,z_1)$ is the solution of (\ref{minp0}), then
\begin{equation}\label{EulereqZ} 
\begin{array}{c}
\HH(\tilde p+ p_1-p_0,\tilde z+z_1-z_0) - \HH( p_1-p_0,z_1-z_0)\ge
\smallskip
\\
\ge -\langle \sigma_1,\tilde e\rangle + \langle \zeta_1, \tilde z\rangle -
\tfrac{\e}{\tau}\langle p_1-p_0, \tilde p\rangle-
\tfrac{\e}{\tau}\langle z_1-z_0, \tilde z\rangle
\end{array}
\end{equation}
for every $(\tilde u,\tilde e,\tilde p)\in A(0)$ and every $\tilde z\in L^2(\Om)$. 
Condition \eqref{EulereqZ} implies that
\begin{eqnarray}
& \div\,\sigma_1=0  \text{ in } \Om\,, \quad [\sigma_1 n]=0 \text { on } \Ga_1\,, 
\label{3.5}
\\
& ((\sigma_1)_D-\frac{\e}{\tau}(p_1-p_0),\zeta_1 -\frac{\e}{\tau}(z_1-z_0))\in
\partial \HH(p_1-p_0, z_1-z_0)\,.
\label{3.6}
\end{eqnarray}
In the Dirichlet-Periodic case (DP) the condition $[\sigma n]=0$ on $\Gamma_1$ is replaced by the requirement that $\sigma$ is the restriction to $Q$ of an $\ol x$-periodic function satisfying $\div\,\sigma=0$ in~$\Rn$. 

Conversely, conditions \eqref{3.5} and \eqref{3.6} imply that $(u_1,e_1,p_1,z_1)$ is a solution of (\ref{minp0}).
\end{theorem}

\begin{proof}

To prove \eqref{EulereqZ} we fix $(\tilde u,\tilde e,\tilde p)\in A(0)$ and
$\tilde z\in L^2(\Om)$. 
For every $s>0$ the triple $({u_1+s \tilde u},\allowbreak{e_1+s \tilde e}, \allowbreak{p_1+s\tilde p})$ belongs to $A(w_1)$, and hence by minimality,
\begin{eqnarray*}
& \QQ(e_1+s\tilde e)+\HH(s\tilde p+p_1-p_0,s\tilde z+z_1-z_0)+\V(z_1+s\tilde z)+
\\
& + \tfrac{\e}{2\tau}\|s\tilde p+p_1-p_0\|_2^2
+ \tfrac{\e}{2\tau}\|s\tilde z+z_1-z_0\|_2^2
\ge\\
&\ge\QQ(e_1)+\HH(p_1-p_0,z_1-z_0)
+\V(z_1)+\tfrac{\e}{2\tau}\|p_1-p_0\|_2^2+ \tfrac{\e}{2\tau}\|z_1-z_0\|_2^2\,.
\end{eqnarray*}
Using the convexity of $H$ we obtain
\begin{eqnarray*}
& s\{ \HH(\tilde p+p_1-p_0,\tilde z+z_1-z_0)- \HH(p_1-p_0,z_1-z_0) \} \geq
\\
& \geq
 - \QQ(e_1+ s\tilde e) +\QQ(e_1)
-\V(z_1+ s\tilde z)+\V(z_1) - {}
\\
& {} - \tfrac{\e}{2\tau} \|p_1+ s\tilde p-p_0\|_2^2
+ \tfrac{\e}{2\tau}\|p_1-p_0\|_2^2
-\tfrac{\e}{2\tau} \|z_1+ s\tilde z-z_0\|_2^2
+ \tfrac{\e}{2\tau}\|z_1-z_0\|_2^2
\end{eqnarray*}
for every $s>0$.
Taking the derivative with respect to $s$ at $s=0$, we obtain (\ref{EulereqZ}). 

Taking $\tilde e=E\tilde u$, $\tilde p=0$, and $\tilde z=0$ in   (\ref{EulereqZ}) we deduce that $\langle\sigma_1,E\tilde u\rangle=0$ for every $\tilde u\in H^1(\Om)$ with $\tilde u=0$
$\hn$-a.e\ on $\Gamma_0$, which gives~\eqref{3.5} thanks to \eqref{sigmanu}. The changes in the Dirichlet-Periodic case are obvious.

Taking $\tilde u=0$ and $\tilde e=-\tilde p$, we obtain from  (\ref{EulereqZ}) 
$$
\begin{array}{c}
\HH(\tilde p+ p_1-p_0,\tilde z+z_1-z_0) - \HH( p_1-p_0,z_1-z_0)\ge
\smallskip
\\
\ge \langle \sigma_1,\tilde p\rangle + \langle \zeta_1, \tilde z\rangle -
\tfrac{\e}{\tau}\langle p_1-p_0, \tilde p\rangle-
\tfrac{\e}{\tau}\langle z_1-z_0, \tilde z\rangle=
\\
=\langle (\sigma_1)_D,\tilde p\rangle + \langle \zeta_1, \tilde z\rangle -
\tfrac{\e}{\tau}\langle p_1-p_0, \tilde p\rangle-
\tfrac{\e}{\tau}\langle z_1-z_0, \tilde z\rangle
\end{array}
$$
for every $\tilde p\in L^2(\Om;\MD)$ and every $\tilde z\in L^2(\Om)$,  which 
implies~\eqref{3.6}. 

The fact that \eqref{3.5} and \eqref{3.6} imply minimality follows from the strict convexity of the functional.
\end{proof}

The following theorem provides a dual formulation of problem~(\ref{minp0}).

\begin{theorem}\label{mindual}
Let $\e>0$, $\tau\in(0,\e/M)$, $w_0\in H^1(\Om;\Rn)$, $(u_0,e_0,p_0)\in A(w_0)$,
$z_0\in L^2(\Om)$, $\sigma_0:=\C e_0$, and $\zeta_0:=-V'(z_0)$. 
Let $w_1\in H^1(\Om;\Rn)$, let $(u_1,e_1,p_1,z_1)$
be the solution of (\ref{minp0}), and let $\sigma_1:=\C e_1$ and $\zeta_1:=-V'(z_1)$. 
Then 
\begin{equation}\label{subPar}
\big(\tfrac{1}{\tau}(p_1-p_0), \tfrac{1}{\tau}(z_1-z_0)\big)=
\partial \HH_\e^*((\sigma_1)_D,\zeta_1 )\,.
\end{equation}
Moreover $(\sigma_1,\zeta_1)$  is a solution to
\begin{equation}\label{minSigma}
\min  \{ \tfrac{1}{\tau}\QQ^*(\sigma-\sigma_0)
+ \HH_\e^*(\sigma_D,\zeta)
-\tfrac{1}{\tau}\langle  \zeta, z_1-z_0\rangle
 -\tfrac{1}{\tau}\langle \sigma , Ew_1-Ew_0 \rangle  \}
\end{equation}
over the set of all $(\sigma,\zeta)\in L^2(\Om;\Mnn){\times}L^2(\Om)$ with 
$\div\,\sigma=0$ in $\Om$ and $[\sigma n]=0$ on $\Gamma_1$. In the Dirichlet-Periodic case (DP) the condition $[\sigma n]=0$ on $\Gamma_1$ is replaced by the requirement that $\sigma$ is the restriction to $Q$ of an $\ol x$-periodic function satisfying $\div\,\sigma=0$ in~$\Rn$.
\end{theorem}

\begin{proof}
As $\partial H$ is positively homogeneous of degree $0$, formula \eqref{3.6} can be written as
$$
\big((\sigma_1)_D-\tfrac{\e}{\tau}(p_1-p_0),\,\zeta_1 -\tfrac{\e}{\tau}(z_1-z_0)\big)\in
\partial \HH\big(\tfrac{1}{\tau}(p_1-p_0), \tfrac{1}{\tau}(z_1-z_0)\big)\,,
$$
which is equivalent to
$$
((\sigma_1)_D,\zeta_1 )\in
\partial \HH_\e\big(\tfrac{1}{\tau}(p_1-p_0), \tfrac{1}{\tau}(z_1-z_0)\big)\,.
$$
By a general duality formula (see, e.g., \cite[Corollary~I.5.2]{Eke-Tem}) this is equivalent to \eqref{subPar}.

Using \eqref{Q*}, we obtain
\begin{equation}\label{451}
\tfrac1\tau\QQ^*(\sigma-\sigma_0)\geq \tfrac1\tau\QQ^*(\sigma_1-\sigma_0)+
\tfrac1\tau \langle \sigma-\sigma_1 , e_1-e_0\rangle
\end{equation}
for every $\sigma\in L^2(\Om;\Mnn)$. By \eqref{subPar} we have
\begin{eqnarray*}
&
\HH_\e^*(\sigma_D,\zeta)\geq \HH_\e^*((\sigma_1)_D,\zeta_1 )+
 \tfrac{1}{\tau} \langle \sigma_D-(\sigma_1)_D , p_1-p_0 \rangle
 +\tfrac{1}{\tau} \langle \zeta-\zeta_1  , z_1-z_0 \rangle=
 \\
 &
 =\HH_\e^*((\sigma_1)_D,\zeta_1 )+
 \tfrac{1}{\tau} \langle \sigma-\sigma_1  , p_1-p_0\rangle
 +\tfrac{1}{\tau} \langle \zeta-\zeta_1 , z_1-z_0 \rangle\,.
\end{eqnarray*}
Since $p_0=Eu_0-e_0$ and $p_1=Eu_1-e_1$, we obtain
\begin{equation}\label{452}
\begin{array}{c}
\HH_\e^*(\sigma_D,\zeta)\geq \HH_\e^*((\sigma_1)_D,\zeta_1 )
+ \tfrac{1}{\tau} \langle \sigma-\sigma_1  , Eu_1-Eu_0 \rangle-{}
\\
{}- \tfrac{1}{\tau} \langle \sigma-\sigma_1 , e_1-e_0 \rangle
+\tfrac{1}{\tau} \langle \zeta-\zeta_1 , z_1-z_0 \rangle
\end{array}
\end{equation}
for every $(\sigma,\zeta)\in L^2(\Om;\Mnn){\times}L^2(\Om)$. 
Adding \eqref{451} and \eqref{452} we deduce
\begin{equation}\label{452bis}
\begin{array}{c}
\tfrac{1}{\tau}\QQ^*(\sigma-\sigma_0) + \HH_\e^*(\sigma_D,\zeta)
- \tfrac{1}{\tau} \langle  \zeta , z_1-z_0\rangle
- \tfrac{1}{\tau} \langle  \sigma , Eu_1-Eu_0\rangle \geq
\\
\geq \tfrac{1}{\tau}\QQ^*(\sigma_1-\sigma_0) + \HH_\e^*((\sigma_1)_D,\zeta_1 )
- \tfrac{1}{\tau} \langle \zeta_1 , z_1-z_0\rangle
- \tfrac{1}{\tau} \langle \sigma_1  , Eu_1-Eu_0 \rangle\,.
\end{array}
\end{equation}
If $\div\,\sigma=0$ in $\Om$ and $[\sigma n]=0$ on $\Gamma_1$, integrating by parts and using the equality $u_1-u_0 = w_1 - w_0$ $\hn$-a.e. on $\Ga_0$, we obtain $\langle Eu_1-Eu_0, \sigma\rangle=\langle Ew_1-Ew_0, \sigma\rangle$. Together with \eqref{452bis}, this proves that $(\sigma_1,\zeta_1)$ is a solution to \eqref{minSigma}.
\end{proof}
\end{section}

\begin{section}{Regularized evolution}\label{reg-evol}

In this section we study the notion of regularized evolution.

\subsection{Definition and properties}
We begin with the definition.

\begin{definition}\label{def:reym}
Let   $\ww\in H^1_{loc}([0,+\infty);H^1(\Om;\Rn))$, $u_0\in H^1(\Om; \Rn)$, $e_0\in L^2(\Om;\Mnn)$, $p_0\in L^2(\Om;\MD)$, $z_0\in L^2(\Om)$, and let $\e>0$.
A {\em solution of the $\e$-regularized evolution problem\/}
with boundary datum $\ww$ and initial condition $(u_0, e_0, p_0,z_0)$ is a function
$(\uu_\e,\ee_\e,\pp_\e,\zz_\e)$, with 
\begin{equation}\label{5567}
\begin{array}{c}
\uu_\e\in H^1_{loc}([0,+\infty);H^1(\Om;\Rn))\,,\qquad
\ee_\e\in H^1_{loc}([0,+\infty);L^2(\Om;\Mnn))\,,
\\
\pp_\e\in H^1_{loc}([0,+\infty);L^2(\Om;\MD))\,,
\qquad
\zz_\e\in H^1_{loc}([0,+\infty);L^2(\Om))\,,
\end{array}
\end{equation}
such that, setting 
$$ 
\sigmaa\!_\e(t):=\C \ee_\e(t)\,, \qquad \zetaa_\e(t):=-V'(\zz_\e(t))\,,
$$ 
the following conditions are satisfied:
\begin{itemize}
\smallskip
\item[(ev0)$\!_\e$] {\it initial condition:\/} $(\uu_\e(0),\ee_\e(0),\pp_\e(0),\zz_\e(0))=
(u_0, e_0, p_0,z_0)$;
\smallskip
\item[(ev1)$\!_\e$] {\it kinematic admissibility:\/} for every $t\in(0,+\infty)$
$$
\begin{array}{c}
E\uu_\e(t)=\ee_\e(t)+\pp_\e(t) \quad \text{a.e.\ in }\Om\,,
\smallskip
\\
\uu_\e(t)=\ww(t)
\quad \hn\text{-a.e.\ in } \Gamma_0\,;
\end{array}
$$
\item[(ev2)$\!_\e$] {\it equilibrium condition:\/} for every $t\in(0,+\infty)$
$$
\div\,\sigmaa_\e(t)=0 \text{ in } \Om\,, \qquad
[\sigmaa_\e(t) n]=0 \text{ on } \Ga_1\,;
$$
\item[(ev\^3)$\!_\e$] {\it regularized flow rule:\/} for a.e.\ $t\in(0,+\infty)$
$$
(\dot\pp_\e(t),\dot\zz_\e(t))=N^\e_K(\sigmaa_\e(t)_D,\zetaa_\e(t))
\quad \text{a.e.\ in }\Om\,,
$$
where $N^\e_K$ is defined by~\eqref{partialHe*}.
\smallskip
\end{itemize}

In the Dirichlet-Periodic case (DP) we assume that for every $t\in [0,+\infty)$ the functions $\ww(t)$, $u_0$, and  $\uu_\e(t)$ are restrictions to $Q$ of 
$\ol x$-periodic functions of class 
$H^1_{loc}(\Rn;\Rn)$, and in the equilibrium condition (ev2)$\!_\e$ the equality 
$[\sigmaa_\e(t) n]=0$ on $\Ga_1$ is replaced by the requirement that $\sigmaa_\e(t)$ is the restriction of an $\ol x$-periodic function with 
$\div\,\sigmaa_\e(t)=0$ in~$\Rn$.
\end{definition}

\begin{remark}\label{rmk:eq-def}
Let us fix $t>0$ such that the derivatives $\dot \pp_\e(t)$ and $\dot \zz_\e(t)$ exist.
Then the following conditions are equivalent:
\begin{eqnarray}
&\label{dualdef}
(\dot\pp_\e(t),\dot\zz_\e(t))=N^\e_K(\sigmaa_\e(t)_D,\zetaa_\e(t))
\quad \text{a.e.\ in }\Om\,,
\\
&\label{primdef}
(\sigmaa_\e(t)_D,\zetaa_\e(t))\in\partial\HH_\e(\dot\pp_\e(t),\dot\zz_\e(t))
\quad \text{a.e.\ in }\Om\,,
\\
&\label{primdef2}
(\sigmaa_\e(t)_D-\e \dot\pp_\e(t),\zetaa_\e(t)-\e\dot\zz_\e(t))\in\partial\HH(\dot\pp_\e(t),\dot\zz_\e(t))
\quad \text{a.e.\ in }\Om\,.
\end{eqnarray}
Indeed, by \eqref{dHeNe} we have $\partial \HH^*_\e (\sigmaa_\e (t), \zetaa_|e (t))
= N^\e_K(\sigmaa_\e (t), \zetaa_\e (t))$, so that \eqref{dualdef} and \eqref{primdef} are equivalent by a standard property of conjugate functions (see, e.g., \cite[Corollary~I.5.2]{Eke-Tem}).
The equivalence between \eqref{primdef} and \eqref{primdef2} follows immediately from the definition of $\HH_\e$ (see~\eqref{He}).
\end{remark}

The following theorem shows that the modified flow rule (ev\^3)$\!_\e$ can be replaced
by a suitable stress constraint and an energy equality. To formulate these conditions
it is convenient to introduce the convex set 
\begin{equation}\label{KOm}
\K(\Om):=\{(\sigma,\zeta)\in L^2(\Om;\MD){\times}L^\infty(\Om): 
(\sigma(x),\zeta(x))\in K\, \hbox{ for a.e. }\,x\in\Om\}\,.
\end{equation}

\begin{theorem}\label{equivalence}
Let  $\ww$, $u_0$, $e_0$, $p_0$, $z_0$, and $\e$ be as in Definition~\ref{def:reym}, and let $(\uu_\e,\ee_\e,\pp_\e,\zz_\e)$ be a function satisfying \eqref{5567}.
Then $(\uu_\e,\ee_\e,\pp_\e,\zz_\e)$ is a solution of the $\e$-regularized evolution problem with boundary datum $\ww$ and initial condition $(u_0, e_0, p_0,z_0)$ if and only if it satisfies the initial condition (ev0)$\!_\e$, the kinematical admissibility 
(ev1)$\!_\e$, the equilibrium condition (ev2)$\!_\e$, and the following properties:
\begin{itemize}
\smallskip
\item[(ev3)$\!_\e$] {\it modified stress constraint:\/} for a.e.\ $t\in(0,+\infty)$
$$
(\sigmaa_\e(t)_D- \e\dot\pp_\e(t),\zetaa_\e(t)- \e\dot \zz_\e(t))\in \K(\Om)\,;
$$
\smallskip 
\item[(ev4)$\!_\e$] {\it energy equality:\/} for every $T>0$ we have 
\begin{eqnarray*}
&\displaystyle\QQ(\ee_\e(T))+\int_{0}^{T}\HH(\dot\pp_\e(t),\dot\zz_\e(t))\,dt
+\V(\zz_\e(T))
+\e\int_{0}^{T}\|\dot \pp_\e(t)\|_2^2\, dt + \e
\int_{0}^{T} \|\dot \zz_\e(t)\|_2^2\, dt
=
 \\
&\displaystyle= \QQ(e_0)+\V(z_0)+
\int_{0}^{T} \langle\sigmaa_\e(t), E\dot \ww(t)\rangle\, dt\,.
\end{eqnarray*}
\end{itemize}
\end{theorem}

\begin{proof} 
Suppose that $(\uu_\e, \ee_\e,\pp_\e,\zz_\e)$ satisfies (ev1)$\!_\e$, (ev2)$\!_\e$, and
(ev\^3)$\!_\e$.
As $\HH$ is positively homogeneous of degree one and $\partial H(0,0)=K$ (see, e.g., \cite[Corollary~23.5.3]{Roc}), by a general property of integral functionals (see, e.g., \cite[Proposition IX.2.1]{Eke-Tem}), we have
$$
\partial\HH(\dot\pp^\e(t),\dot\zz^\e(t))\subset \partial\HH(0,0)=\K(\Om)\,.
$$
Therefore \eqref{primdef2} implies (ev3)$\!_\e$.

Since $\HH$ is positively homogeneous of degree one, the Euler relation gives
$\langle\sigma,p\rangle+\langle\zeta, z\rangle=\HH(p,z)$ whenever 
$(\sigma,\zeta)\in\partial \HH(p,z)$. Therefore, \eqref{primdef2} implies 
\begin{equation}\label{z0}
\HH(\dot\pp_\e(t),\dot\zz_\e(t))=\langle \sigmaa_\e(t)_D-\e\dot\pp_\e(t), \dot\pp_\e(t)\rangle +\langle\zetaa_\e(t)-\e\dot\zz_\e(t),\dot\zz_\e(t)\rangle\,,
\end{equation}
which is equivalent to
\begin{equation}\label{z1}
\HH(\dot\pp_\e(t),\dot\zz_\e(t)) +\langle V'(\zz_\e(t)),\dot\zz_\e(t)\rangle
+\e\|\dot\pp_\e(t)\|_2^2+\e\|\dot\zz_\e(t)\|_2^2 =
\langle  \sigmaa_\e(t),\dot\pp_\e(t)\rangle\,.
\end{equation}
By (ev1)$\!_\e$ we have
\begin{equation}\label{z2}
\langle  \sigmaa_\e(t),\dot\pp_\e(t)\rangle=
\langle \sigmaa_\e(t),E\dot\uu_\e(t)\rangle-
\langle  \sigmaa_\e(t),\dot\ee_\e(t)\rangle\,.
\end{equation}
Integrating by parts and using (ev2)$\!_\e$, we obtain
\begin{equation}\label{z3}
\langle \sigmaa_\e(t),E\dot\uu_\e(t)\rangle=
\langle \sigmaa_\e(t),E\dot\ww(t)\rangle\,.
\end{equation}
Combining  \eqref{z1}, \eqref{z2}, and \eqref{z3}, we deduce that
\begin{equation}\label{z4}
\begin{array}{c}
\langle  \sigmaa_\e(t),\dot\ee_\e(t)\rangle+
\HH(\dot\pp_\e(t),\dot\zz_\e(t)) +\langle V'(\zz_\e(t)),\dot\zz_\e(t)\rangle+{}
\smallskip
\\
{}+\e\|\dot\pp_\e(t)\|_2^2+\e\|\dot\zz_\e(t)\|_2^2 =
\langle \sigmaa_\e(t),E\dot\ww(t)\rangle\,.
\end{array}
\end{equation}
The energy equality (ev4)$\!_\e$ can be obtained from \eqref{z4} by integration.

Conversely, assume that $(\uu_\e, \ee_\e,\pp_\e,\zz_\e)$ satisfies conditions (ev1)$\!_\e$, (ev2)$\!_\e$, (ev3)$\!_\e$, and (ev4)$\!_\e$. By differentiating (ev4)$\!_\e$ we obtain \eqref{z4}. Thanks to \eqref{z2} and \eqref{z3}, from \eqref{z4} we deduce \eqref{z1}, which is equivalent to \eqref{z0}. 
By  (ev3)$\!_\e$ for a.e. $t\in(0,+\infty)$ we have
\begin{equation}\label{primdef3}
(\sigmaa_\e(t)_D-\e \dot\pp_\e(t),\zetaa_\e(t)-\e\dot\zz_\e(t))\in \K(\Om)=\partial\HH(0,0)\,.
\end{equation}
Using the homogeneity of $\HH$, condition \eqref{primdef2} follows easily from \eqref{z0} and \eqref{primdef3}.
\end{proof}

\begin{theorem}\label{young-main2} 
Let  $\ww$, $u_0$, $e_0$, $p_0$, $z_0$, and $\e$ be as in Definition~\ref{def:reym}, and let $\sigma_0:= \C e_0$. Suppose that the following conditions are satisfied:
\begin{itemize}
\smallskip
\item[(ev1)$\!_0$] {\it kinematic admissibility:\/}
$$
\begin{array}{c}
Eu_0=e_0+p_0 \quad \text{a.e.\ in }\Om\,,
\smallskip
\\
u_0=\ww(0)
\quad \hn\text{-a.e.\ in } \Gamma_0\,.
\end{array}
$$
\smallskip
\item[(ev2)$\!_0$] {\it equilibrium condition:\/} 
$$
\div \, \sigma_0=0 \hbox{ in } \Om\,, \qquad [\sigma_0 n]=0 \hbox { on } \Ga_1\,.
$$
\end{itemize}
Then there exists a unique solution of the
$\e$-regularized evolution problem with boundary datum $\ww$ and initial condition 
$(u_0, e_0, p_0,z_0)$.
\end{theorem}

The proof is postponed to the next subsection. 

Let ${\mathbb M}^{d{\times}d}$ be the space of all $d{\times}d$-matrices. The uniqueness result of the previous theorem allows to simplify the problem in the spatially homogeneous case, when the boundary condition has the form
\begin{equation}\label{bd-hom}
\ww(t,x):=\xi(t)x\,, \qquad \Ga_0:= \partial \Om\,,
\end{equation}
with
\begin{equation}\label{xit}
\xi\in H^1_{loc}([0,+\infty);{\mathbb M}^{d{\times}d} )\,,
\end{equation}
while the initial condition has the form
\begin{equation}\label{initial-hom}
u_0(x):=\xi_0x\,, \quad e_0(x):=\xi_0^e\,,\quad p_0(x):=\xi_0^p\,,\quad z_0(x):=\theta_0\,,
\end{equation}
with 
\begin{equation}
\xi_0\in {\mathbb M}^{d{\times}d}\,,\quad \xi_0^e\in\Mnn\,,\quad \xi_0^p\in \MD\,,\quad\theta_0\in \R\,.\label{initial-hom2}
\end{equation}
We assume the following compatibility condition
\begin{equation}
\xi_0^s= \xi_0^e+ \xi_0^p\,, \quad \xi(0)=\xi_0\,,\label{initial-hom3}
\end{equation}
where $\xi_0^s$ denotes the symmetric part of $\xi_0$.
\begin{proposition}\label{prop:sp-hom}
Assume that $\ww$, $\Ga_0$,  $u_0$, $e_0$, $p_0$,  $z_0$ satisfy \eqref{bd-hom}--\eqref{initial-hom3}, and let $\e>0$. Then $(\uu_\e,\ee_\e,\pp_\e,\zz_\e)$ is the solution of the
$\e$-regularized evolution problem with boundary datum $\ww$ and initial condition 
$(u_0, e_0, p_0,z_0)$  if and only if 
\begin{equation}\label{sol-hom}
\uu_\e(t,x)=\xi(t)x\,,\quad \ee_\e(t,x)=\xi^e_\e(t)\,,\quad\pp_\e(t,x)=\xi^p_\e(t)\,,\quad \zz_\e(t,x)=\theta_\e(t)\,,
\end{equation}
with
\begin{equation}\label{xip}
\xi_\e^p(t):=\xi^s(t)-\xi_\e^e(t)\,,
\end{equation}
where $\xi^s(t)$ denotes the symmetric part of $\xi(t)$ and $(\xi_\e^e, \theta_\e)$ is the unique solution in $H^1_{loc}([0, +\infty);\Mnn){\times}H^1_{loc}([0, +\infty))$ of the Cauchy problem
\begin{equation}\label{equation-diff}
(\dot\xi^s(t)-\dot\xi_\e^e(t),\dot\theta_\e(t))=N_K^\e(\C_D\xi^e_\e(t)_D,-V'(\theta_\e(t)))\,,
\end{equation}
\begin{equation}\label{initial-Cauchy}
\xi_\e^e(0)=\xi_0^e\,,\quad \theta_\e(0)=\theta_0\,.
\end{equation}
\end{proposition}
\begin{proof}
Let $(\xi_\e^e, \theta_\e)$ be the solution of the Cauchy problem \eqref{equation-diff}, \eqref{initial-Cauchy}, whose existence follows from the Lipschitz-continuity of $N_\K^\e$ and $V'$. By the definition of $N_K^\e$ (see \eqref{partialHe*}) it follows that 
$\dot\xi^s(t)-\dot\xi_\e^e(t)\in \MD$ for a.e.\ $t\in [0,+\infty)$. Moreover, by \eqref{initial-hom2}, \eqref{initial-hom3}, and \eqref{initial-Cauchy} we have also $\xi^s(0)-\xi_\e^e(0)\in \MD$. It follows that $\xi^s(t)-\xi_\e^e(t)\in \MD$ for every $t\in [0,+\infty)$,
hence $\xi_\e^p(t)\in \MD$ by \eqref{xip}. 
It is then easy to see that the function $(\uu_\e,\ee_\e,\pp_\e,\zz_\e)$ defined by \eqref{sol-hom} satisfies all  conditions of Definition~\ref{def:reym}. 
\end{proof}

\subsection{Proof of the existence theorem}
To prove Theorem~\ref{young-main2} it is useful to introduce the notion of $H$-dissipation of a possibly discontinuous pair of functions $\pp\colon[0,+\infty)\to M_b(\ol\Om;\MD)$ and $\zz\colon[0,+\infty)\to M_b(\ol\Om)$.
For every pair of times $a$, $b\in [0,+\infty)$, with $a<b$, the $H$-dissipation of $(\pp,\zz)$ on the interval $[a,b]$ is defined by
\begin{equation}\label{diss}
\D_{\!H}(\pp,\zz;a,b):=\sup \sum_{j=1}^k \HH(\pp(t_j)-\pp(t_{j-1}), \zz(t_j)-\zz(t_{j-1}))\,,
\end{equation}
where the supremum is taken over all finite families $t_0, t_1, \dots, t_k$ with
$ a=t_0< t_1< \dots< t_k=b$.
For every $\pp\in H^{1,1}_{loc}([0,+\infty);L^2(\Om;\MD))$ and every 
$\zz\in H^{1,1}_{loc}([0,+\infty);L^2(\Om))$ we have
\begin{equation}\label{diss-reg}
\D_{\!H}(\pp,\zz;a,b)=
\int_{a}^{b}\HH(\dot\pp(t),\dot\zz(t)) \, dt
\end{equation}
(see, e.g., \cite[Theorem~7.1]{DM-DeS-Mor}).

\begin{proof}[Proof of Theorem~\ref{young-main2}] We prove the existence of a solution to the $\e$-regularized evolution problem by time discretization, using an implicit Euler scheme, which leads to an incremental minimization problem. The use of the dual formulation of the problem is inspired by \cite{Suq}. The proof is divided into several steps.

\medskip

\noindent
{\it Step 1. The incremental problems.\/}
Let us fix a sequence of subdivisions $(\tki)_{i\geq 0}$ of the half-line $[0,+\infty)$, with
\begin{eqnarray}
& 0=t_k^0<t_k^1<\dots<t_k^{i-1}<t_k^{i}\to+\infty\,,
\label{subdivX}\\
&\displaystyle
\tau_k:=\sup_i (\tki-\tkim)\to 0\quad\hbox{as } k\to\infty\,.
\label{fineX}
\end{eqnarray}
{}It is not restrictive to assume also that $\tau_k<\e/M$ for every $k$, where $M$ is the constant introduced in \eqref{2der}. 

For every $i$ we set $\wki:=\ww(\tki)$. We define $\uki$, $\eki$, $\pki$, and $\zki$ by induction on
$i$.
We set $(u_k^0,e_k^0, p_k^0, z_k^0):=(u_0,e_0,p_0,z_0)$, and for $i\geq1$ we 
define $(\uki,\eki,\pki,\zki)$ as the solution to the incremental problem
\begin{equation}\label{min-incX}
\min \{ \QQ(e)+\HH(p-p_k^{i-1},z-z_k^{i-1}) +\V(z)
+\tfrac{\e}{2\tau_k^i}\|p-p_k^{i-1}\|^2_2
+\tfrac{\e}{2\tau_k^i}\|z-z_k^{i-1}\|^2_2 \}\,,
\end{equation}
where $\tau_k^i:={\tki-t_k^{i-1}}$ and the minimum is taken under the conditions
$(u,e,p)\in A(\wki)$ and $z\in L^2(\Om)$.
The existence and uniqueness of the solution to this problem is proved in 
Theorem~\ref{existencemin}. 

{}For every $i\geq0$ we set $\sigma_k^i:=\C e_k^i$ and $\zeta_k^i:=-V'(z_k^i)$. We consider the piecewise constant interpolations defined by
\begin{equation}\label{uktX}
\begin{array}{c} 
\uu_k(t):=\uki\,, \quad \ee_k(t):=\eki\,,  \quad   \pp_k(t):=\pki\,,  \quad \zz_k(t):=\zki\,,
\vspace{.1cm}\\
\sigmaa_k(t):=\sigma_k^i\,, \quad   \zetaa_k(t) := \zeta_k^i \,, \quad\ww_k(t):=w_k^i \,, \quad [t]_k:=\tki
\end{array}
\end{equation}
for $t\in [\tki,t^{i+1}_k)$. 
By definition 
\begin{equation}\label{admk}
(\uu_k(t),\ee_k(t),\pp_k(t))\in A(\ww_k(t))
\end{equation}
and
\begin{equation}\label{initk}
\uu_k(0)=u_0\,, \quad \ee_k(0)=e_0\,, \quad \pp_k(0)=p_0\,,  \quad \zz_k(0)=z_0\,.
\end{equation}
A solution of the $\e$-regularized evolution problem will be obtained by taking the limit, as $k\to\infty$, of these piecewise constant interpolations. We will also consider the piecewise affine interpolations 
\begin{equation}\label{affinterp}
\begin{array}{c}
\uu^\smtr_k(t):= u^i_k + (t-t^i_k)(u_k^{i+1} - u_k^{i})/\tau_k^{i+1} \,, \\
\ee^\smtr_k(t):= e^i_k + (t-t^i_k)(e_k^{i+1} - e_k^{i})/\tau_k^{i+1} \,, \\
\pp^\smtr_k(t):= p^i_k + (t-t^i_k)(p_k^{i+1} - p_k^{i})/\tau_k^{i+1} \,, \\
\zz^\smtr_k(t):= z^i_k + (t-t^i_k)(z_k^{i+1} - z_k^{i})/\tau_k^{i+1} \,, \\
\sigmaa^\smtr_k(t):= \sigma^i_k + (t-t^i_k)(\sigma_k^{i+1} - \sigma_k^{i})/\tau_k^{i+1} \,,
\end{array}
\end{equation}
for every $t\in [t^i_k, t^{i+1}_k)$.

We derive now an energy estimate for the solutions of the incremental problems. 

\begin{lemma}\label{lm:41X}
For every $T>0$ there exists a sequence $\omega^T_k\to 0^+$ such that 
\begin{equation}\label{inc-ineqX}
\begin{array}{c}
\displaystyle
\vphantom{\int_{[t_1]_k}^{[t_2]_k}}
\QQ(\ee_k(t_2))+\D_{\!H}(\pp_k,\zz_k;t_1,t_2)+\V(\zz_k(t_2))+{}
\\
\displaystyle
{}+
\frac{\e}{2} \int_{[t_1]_k}^{[t_2]_k}\|\dot \pp^\smtr_k(t)\|_2^2\, dt
+\frac{\e}{2} \int_{[t_1]_k}^{[t_2]_k}\|\dot \zz^\smtr_k(t)\|_2^2\, dt
\le
\medskip
\\
\displaystyle
\le \QQ(\ee_k(t_1))+\V(\zz_k(t_1))
+\int_{[t_1]_k}^{[t_2]_k}\langle\sigmaa_k(t),E\dot \ww(t)\rangle\, dt +\omega^T_k\,.
\end{array}
\end{equation}
for every $k$ and every
$t_1,t_2\in[0,T]$ with $t_1\le t_2$.
\end{lemma}

\begin{proof}
Let us fix $T>0$. Since $(\pp_k,\zz_k)$ is piecewise constant on $[0,T]$,
it is easy to see that
\begin{equation}\label{pc-diss}
\D_{\!H}(\pp_k,\zz_k;t_1,t_2)= \sum_{r=i+1}^j \HH(p_k^r-p_k^{r-1},z_k^r-z_k^{r-1})\,,
\end{equation}
where $i$ and $j$ are such that $t_k^i=[t_1]_k$ and $t_k^j=[t_2]_k$.
Therefore, we have to prove that there exists a sequence $\omega^T_k\to 0^+$ such that
\begin{equation}\label{inc-ineq-2X}
\begin{array}{c}
\displaystyle
\vphantom{\int_{\tki}^{t_k^j}}
\QQ(e_k^j)+
\sum_{r=i+1}^j \HH(p_k^r-p_k^{r-1},z_k^r-z_k^{r-1}) +\V(z_k^j)+{}
\\
\displaystyle
{}+\frac{\e}{2} \int_{\tki}^{t_k^j}\|\dot \pp^\smtr_k(t)\|_2^2\, dt
+\frac{\e}{2} \int_{\tki}^{t_k^j}\|\dot \zz^\smtr_k(t)\|_2^2\, dt
\le
\\
\displaystyle
\le \QQ(\eki)+\V(\zki)
+\int_{\tki}^{t_k^j}\langle\sigmaa_k(t),E\dot \ww(t)\rangle\, dt +\omega^T_k\,.
\end{array}
\end{equation}
for every $i,j,k$ with $0\le i\le j$ and $t^j_k\le T$. 

Let us fix an integer $r$  with $i < r \leq j$ and let $\hat u:=u_k^{r-1}-w_k^{r-1}+w_k^r$
and $\hat e:=e_k^{r-1}-Ew_k^{r-1}+Ew_k^r$; then, 
$(\hat u,\hat e, p_k^{r-1})\in A(w_k^r)$. Testing with 
$(\hat u,\hat e, p_k^{r-1},z_k^{r-1})$, by the minimality condition
(\ref{min-incX}) we have
\begin{equation}\label{b01X}
\begin{array}{c}
\QQ(e_k^r)+\HH(p_k^r-p_k^{r-1},z_k^r-z_k^{r-1})+
\V(z_k^r) +{} \smallskip
\\
{}+\tfrac{\e}{2\tau_k^r}\|p_k^r-p_k^{r-1}\|_2^2
+\tfrac{\e}{2\tau_k^r}\|z_k^r-z_k^{r-1}\|_2^2
 \le  \smallskip
 \\
\le 
\QQ(e_k^{r-1}+Ew_k^r -Ew_k^{r-1})+\V(z_k^{r-1})\,,
\end{array}
\end{equation}
where the quadratic form in the right-hand side can be developed as
\begin{equation}\label{b02X}
\QQ(e_k^{r-1}+Ew_k^r-Ew_k^{r-1})
=\QQ(e_k^{r-1})+\langle \sigma_k^{r-1}, Ew_k^r- 
Ew_k^{r-1} \rangle 
+ \QQ(Ew_k^r-Ew_k^{r-1})\,.
\end{equation}
{}From the absolute continuity of $\ww$ with respect to $t$ we obtain
$$
w_k^r-w_k^{r-1}=\int_{t_k^{r-1}}^{t_k^r}\dot \ww(t)\, dt\,,
$$
where the right-hand side is a Bochner integral of a function with values in 
$H^1(\Om;\Rn)$.
This implies that
\begin{equation}\label{b03X}
Ew_k^r- Ew_k^{r-1}=\int_{t_k^{r-1}}^{t_k^r} E\dot \ww(t)\, dt\, ,
\end{equation}
where the right-hand side is a Bochner integral of a function with values in 
$L^2(\Om;\Mnn)$.
By (\ref{boundsC}) and (\ref{b03X}) we get
\begin{equation}\label{b03.5}
\QQ(Ew_k^r-Ew_k^{r-1})\le \beta_{\C}\Big(\int_{t_k^{r-1}}^{t_k^r} 
\|E\dot \ww(t)\|_2\, dt\Big)^2. 
\end{equation}
Note that we can also write
\begin{equation}\label{b03.55X}
\begin{array}{c}
\displaystyle
\|p_k^r-p_k^{r-1}\|_2^2=(t_k^r-t_k^{r-1})
\int_{t_k^{r-1}}^{t_k^r}\|\dot \pp^\smtr_k(t)\|_2^2\, dt\,,\smallskip
\\
\displaystyle
\|z_k^r-z_k^{r-1}\|_2^2=(t_k^r-t_k^{r-1})
\int_{t_k^{r-1}}^{t_k^r}\|\dot \zz^\smtr_k(t)\|_2^2\, dt\,.
\end{array}
\end{equation}
By (\ref{b01X})--(\ref{b03.55X}) we obtain
\begin{equation}\label{b04X}
\begin{array}{c}
\displaystyle
\vphantom{\int_{t_k^{r-1}}^{t_k^r}}
\QQ(e_k^r)+\HH(p_k^r-p_k^{r-1},z_k^r-z_k^{r-1})+
\V(z_k^r)+\frac{\e}{2}\int_{t_k^{r-1}}^{t_k^r}
\big( \|\dot \pp^\smtr_k(t)\|_2^2+
\|\dot \zz^\smtr_k(t)\|_2^2\big)\, dt \le \vspace{.1cm} 
\\
\le \displaystyle
\QQ(e_k^{r-1})+ \V(z_k^{r-1})
+\int_{t_k^{r-1}}^{t_k^r} \langle \sigma_k^{r-1}, E\dot \ww(t)\rangle \, dt
 +\beta_{\C}\Big(\int_{t_k^{r-1}}^{t_k^r}\|E\dot \ww(t)\|_2\, dt 
\Big)^2
\le \vspace{.1cm} 
\\
\le \displaystyle
\QQ(e_k^{r-1})+\V(z_k^{r-1})
+\int_{t_k^{r-1}}^{t_k^r} \langle \sigma_k^{r-1}, E\dot \ww(t)\rangle \, dt
+ \rho^T_k\int_{t_k^{r-1}}^{t_k^r}\|E\dot \ww(t)\|_2\, dt\,,
\end{array}
\end{equation}
where
$$\rho^T_k:=\max_{t^r_k\le T}\beta_{\C}\int_{t_k^{r-1}}^{t_k^r}\|E\dot 
\ww(t)\|_2\, dt \, \to 0$$
by the absolute continuity of the integral.
Iterating now inequality (\ref{b04X}) for $i+1\le r\le j$, we get 
(\ref{inc-ineq-2X}) 
with $\omega^T_k:=\rho^T_k\int_0^T \|E\dot \ww(t)\|_2\, dt$.
\end{proof}

We now prove a dual energy estimate, where $\alpha_\C$ and $M$ are the constants introduced in \eqref{boundsC} and \eqref{2der}.

\begin{lemma}\label{dualestimate} 
We have
\begin{equation}\label{newine}
\begin{array}{c}
\displaystyle
\alpha_\C \int_{[t_1]_k}^{[t_2]_k}\|\dot \ee^\smtr_k(t)\|_2^2\, dt
+ \HH_\e^*(\sigmaa_k(t_2)_D, \zetaa_k(t_2))- \HH_\e^*(\sigmaa_k(t_1)_D, \zetaa_k(t_1))
\le
\\
\\
\displaystyle
\le
M\int_{[t_1]_k}^{[t_2]_k}\|\dot \zz^\smtr_k(t)\|_2^2\, dt
+\int_{[t_1]_k}^{[t_2]_k}\langle\dot\sigmaa^\smtr_k(t),E\dot \ww(t)\rangle\, dt\,,
\end{array}
\end{equation}
for every $k$ and every
$t_1,t_2\in[0,+\infty)$ with $t_1<t_2$.
\end{lemma}

\begin{proof}
We have to prove that 
\begin{equation}\label{formula1000}
\begin{array}{c}
\displaystyle
\alpha_\C \int_{t^i_k}^{t^j_k}\|\dot \ee^\smtr_k(t)\|_2^2\, dt
+ \HH_\e^*((\sigma_k^j)_D, \zeta_k^j)- \HH_\e^*((\sigma_k^{i})_D, \zeta_k^{i}))
\le
\\
\displaystyle
\le
M\int_{t^i_k}^{t^j_k}\|\dot \zz^\smtr_k(t)\|_2^2\, dt
+\int_{t^i_k}^{t^j_k}\langle\dot\sigmaa^\smtr_k(t),E\dot \ww(t)\rangle\, dt\,,
\end{array}
\end{equation}
for every $k$ and every $0\le i\le j $. Let us fix an integer $r$ with $i < r \leq j$. We observe that $\div\,\sigma^{r-1}_k=0$ in $\Om$ and $[\sigma^{r-1}_k n]=0$ on $\Ga_1$.  This follows from \eqref{3.5} if $r>1$, and from the equilibrium condition (ev2)$\!_0$ if $r=1$. By Theorem~\ref{mindual} we have 
\begin{eqnarray*}
&\tfrac{1}{\tau_k^r} \QQ^*(\sigma^r_k-\sigma^{r-1}_k )+\HH^*_\e((\sigma^r_k)_D,\zeta^r_k)
-\HH^*_\e((\sigma^{r-1}_k)_D,\zeta^{r-1}_k) \leq \\
&\leq \tfrac{1}{\tau_k^r} \langle \zeta^r_k - \zeta^{r-1}_k , z^r_k - z^{r-1}_k
 \rangle
+\tfrac{1}{\tau_k^r} \langle \sigma^r_k - \sigma^{r-1}_k, Ew^r_k - Ew^{r-1}_k
  \rangle\,.
\end{eqnarray*}
Since $\alpha_\C\|e^r_k - e^{r-1}_k  \|^2_2 \leq \QQ^*(\sigma^r_k-\sigma^{r-1}_k)$ by \eqref{boundsC} and
$\langle \zeta^r_k - \zeta^{r-1}_k , z^r_k - z^{r-1}_k \rangle
 \leq M \|z^r_k - z^{r-1}_k  \|^2_2$
 by \eqref{2der}, we obtain
 \begin{eqnarray*}
&  \displaystyle
\alpha_\C \int_{t^{r-1}_k}^{t^{r}_k}\|\dot \ee^\smtr_k(t)\|_2^2\, dt
+ \HH_\e^*((\sigma_k^r)_D, \zeta_k^r)- \HH_\e^*((\sigma_k^{r-1})_D, \zeta_k^{r-1})
\le
\\
& \displaystyle
\le
M\int_{t^{r-1}_k}^{t^{r}_k}\|\dot \zz^\smtr_k(t)\|_2^2\, dt
+\int_{t^{r-1}_k}^{t^{r}_k}\langle\dot\sigmaa^\smtr_k(t),E\dot \ww(t)\rangle\, dt\,.
\end{eqnarray*}
Summing over $r$ we obtain \eqref{formula1000}.
\end{proof}

\noindent
{\it Proof of Theorem~\ref{young-main2}. Step 2. Proof of the bounds.\/}
Let us prove now that for every $T>0$ there exists a constant $C_T$, independent of $k$ and $\e$, such that
\begin{eqnarray}
& \label{b1005}
\displaystyle
 \sup_{t\in[0,T]}\|\ee_k(t)\|_2\leq C_T\,, \qquad 
\D_{\!H}(\pp_k,\zz_k;0,T)\leq C_T\,,
\\
& \label{B1005} 
\displaystyle
\int_0^T \|\dot \ee^\smtr_k(t)\|_2^2\, dt\leq \frac{C_T}{\e\land 1}\,, \quad 
\int_0^T \|\dot \pp^\smtr_k(t)\|_2^2\, dt\leq \frac{C_T}\e\,, \quad 
\int_0^T \|\dot \zz^\smtr_k(t)\|_2^2\, dt\le \frac{C_T}\e \,.
\end{eqnarray}
By (\ref{gammaM}) and (\ref{diss}) for every $t\in [0,+\infty)$  we have
\begin{equation}\label{bdiss10}
\begin{array}{c}
\D_{\!H}(\pp_k,\zz_k;0,t)+ \V(\zz_k(t)) \ge
\HH(\pp_k(t)-p_0,\zz_k(t)-z_0)+ \V(\zz_k(t)) \ge 
\smallskip
\\
\ge   C^K_V \|\pp_k(t)-p_0\|_1 +    C^K_V \|\zz_k(t)-z_0\|_1
+\V(z_0)\,.
\end{array}
\end{equation}
Let us fix $T>0$.
Using the discrete energy inequality (\ref{inc-ineqX}) with $t_1=0$ and $t_2=t\le T$ and inequalities
(\ref{boundsC}), (\ref{normC}), and (\ref{bdiss10}), we deduce that
$$
\alpha_{\C}\|\ee_k(t)\|^2_2 \le 
\beta_{\C}\|e_0\|^2_2  + 2\beta_\C\, \sup_{t\in[0,T]}\|\ee_k(t)\|_2
\int_0^T  \|E\dot\ww(t)\|_2 \, dt +\omega_k^T
$$
for every $k$ and every $t\in[0,T]$. The first estimate in (\ref{b1005}) can be obtained now by using the Cauchy inequality.

By (\ref{inc-ineqX}) and the first inequality in (\ref{b1005}) we have that 
\begin{equation}\label{10024}
 \D_{\!H}(\pp_k,\zz_k;0,t)+\V(\zz_k(t))\le C
\end{equation}
with $C$ independent of $k$ and $t$. By (\ref{bdiss10}) this implies the boundedness of
$\|\pp_k(t)\|_1$ and $\|\zz_k(t)\|_1$, and in turn, the boundedness of $\V(\zz_k(t))$. 
Now we can use  (\ref{10024}) to obtain the second estimate in (\ref{b1005}).
The last two inequalities in (\ref{B1005}) follow immediately from (\ref{inc-ineqX}) and 
\eqref{bdiss10}. The first inequality in (\ref{B1005}) can be obtained from \eqref{newine} with $t_1=0$ and $t_2=T$ thanks to \eqref{He*} and \eqref{initk}.

To continue the proof of Theorem~\ref{young-main2} we need the following lemma, based on Gronwall's inequality.

\begin{lemma}\label{Gronwall} 
For every $T>0$ the sequences $\ee_k$ and $\zz_k$ satisfy the Cauchy condition in
$L^\infty([0,T];L^2(\Om;\Mnn))$ and $L^\infty([0,T];L^2(\Om))$, respectively.
\end{lemma}

\begin{proof}
By \eqref{dHeNe} and \eqref{subPar} for every $i$ and $k$ we have
$$
(p^i_k,z^i_k)-(p^{i-1}_k, z^{i-1}_k) =\tau^i_k N^\e_K((\sigma^i_k)_D,-V'(z^i_k))\,.
$$
Since $Eu^i_k=e^i_k+p^i_k$ and $Eu^{i-1}_k=e^{i-1}_k+p^{i-1}_k$ a.e.\ in $\Om$,
we obtain
$$
(-e^i_k,z^i_k)-(-e^{i-1}_k,z^{i-1}_k)= \tau^i_k N^\e_K((\sigma^i_k)_D,-V'(z^i_k))
-(Eu^i_k,0)+(Eu^{i-1}_k,0)\,.
$$
Summing for $1\le i\le j$ we obtain
$$
(-e^j_k,z^j_k)-(-e_0,z_0)= 
\sum_{i=1}^j \tau^i_k N^\e_K((\sigma^i_k)_D,-V'(z^i_k))
-(Eu^j_k,0)+(Eu_0,0)\,.
$$
Since $N^\e_K$ is $1/\e$-Lipschitz and $V'$ is $M$-Lipschitz, we have
\begin{equation}\label{4567}
(-e^j_k,z^j_k)-(-e_0,z_0)= 
\sum_{i=1}^j \tau^i_k N^\e_K((\sigma^{i-1}_k)_D,-V'(z^{i-1}_k))+R^j_k
-(Eu^j_k,0)+(Eu_0,0)\,,
\end{equation}
where the rest $R^j_k\in L^2(\Om;\MD){\times}L^2(\Om)$ can be estimated by
\begin{eqnarray*}
&\displaystyle
\|R^j_k\|_2\le \frac{1}\e \sum_{i=1}^j \tau^i_k(\|\sigma^i_k-\sigma^{i-1}_k\|_2+
M\,\|z^i_k-z^{i-1}_k\|_2)\le
\\
&\displaystyle
\le \frac{\tau_k }\e \int_0^{t^j_k}(2\beta_\C\|\dot\ee_k^\smtr(t)\|_2+ M \|\dot\zz_k^\smtr(t)\|_2)\,dt\,,
\end{eqnarray*}
thanks to \eqref{normC}. Using \eqref{B1005} we can prove that for every $T>0$ there exists a constant $A^T$, depending on $\e$, but independent of $j$ and $k$, such that
\begin{equation}\label{Rk}
\|R^j_k\|_2\le A^T\tau_k
\end{equation}
for every $j$ and $k$ with $t^j_k\le T$.

If $t\in [0,T]$, for every $k$ there exists a unique $j$ such that $t^j_k\le t<t^{j+1}_k$.
{}From \eqref{4567} we obtain
\begin{equation}\label{4568}
\begin{array}{c}
\displaystyle
(-\ee_k(t),\zz_k(t))-(-e_0,z_0)= 
\int_0^t N^\e_K(\sigmaa_k(s)_D,-V'(\zz_k(s)))\,ds +{}
\\
\displaystyle\vphantom{\int_0^t}
{}+{\boldsymbol R}_k(t)
-(E\uu_k(t),0)+(Eu_0,0)\,,
\end{array}
\end{equation}
where
$$
{\boldsymbol R}_k(t)=
R^j_k-\int_{t^j_k}^t N^\e_K(\sigmaa_k(s)_D,-V'(\zz_k(s)))\,ds\,.
$$
In these formulas we are using Bochner integrals of functions with values in
$L^2(\Om;\MD)\times L^2(\Om)$.
It follows from \eqref{V'infty}, \eqref{partialHe*}, \eqref{b1005}, and \eqref{Rk} that for every $T>0$ there exists a constant
$B^T$, depending on $\e$, but independent of $t$ and $k$, such that
\begin{equation}\label{Rkt}
\|{\boldsymbol R}_k(t)\|_2\le B^T\tau_k
\end{equation}
for every $t\in [0,T]$ and every $k$.

We now consider two indices $h$ and $k$. Subtracting term by term the equations corresponding to \eqref{4568} we obtain
\begin{equation}\label{4569}
\begin{array}{c}
\displaystyle
(-\ee_k(t),\zz_k(t))-(-\ee_h(t),\zz_h(t))= {\boldsymbol I}_{hk}(t)+
{\boldsymbol R}_k(t)-{\boldsymbol R}_h(t)+{}
\\
\displaystyle\vphantom{\int_0^t}
{}+(E\ww_k(t)-E\uu_k(t),0)-(E\ww_h(t)-E\uu_h(t),0)
-(E\ww_k(t)-E\ww_h(t),0)\,,
\end{array}
\end{equation}
where
$$
{\boldsymbol I}_{hk}(t):=\int_0^t \{N^\e_K(\sigmaa_k(s)_D,-V'(\zz_k(s)))-
N^\e_K(\sigmaa_h(s)_D,-V'(\zz_h(s)))\}\,ds\,.
$$
Since $N^\e_K$ is $1/\e$-Lipschitz and $V'$ is $M$-Lipschitz, using \eqref{normC}
we obtain that there exists a constant $L$, depending on $\e$, but independent of $t$, $h$, $k$, and $T$, such that
\begin{equation}\label{4570}
\|{\boldsymbol I}_{hk}(t)\|_2\le L\int_0^t\varphi_{hk}(s)\,ds
\end{equation}
where
\begin{equation}\label{4571}
\varphi_{hk}(t):=\{\|\ee_k(t)-\ee_h(t)\|_2^2+ \|\zz_k(t)-\zz_h(t)\|_2^2\}^{1/2}\,.
\end{equation}
Since $\ww\in H^1([0,T];H^1(\Om;\Rn))$, there exists a constant $C^T$, independent of $t$ and $k$, such that
$\|E\ww(t)-E\ww(s)\|_2\le C^T|t-s|^{1/2}$ for every $s,t\in [0,T]$. It follows that
\begin{equation}\label{4572}
\|E\ww_k(t)-E\ww_h(t)\|_2\le C^T(\tau_h^{1/2}+\tau_k^{1/2})
\end{equation}
for every $t\in [0,T]$ and for every $h$ and $k$.

Since $\div(\sigmaa_k(t)-\sigmaa_h(t))=0$ in $\Om$ and $[(\sigmaa_k(t)-\sigmaa_h(t))n]=0$ on $\Gamma_1$ by \eqref{3.5}, while $\ww_k(t)-\uu_k(t)=\ww_h(t)-\uu_h(t)=0$ on $\Gamma_0$, 
integrating by parts we obtain
$$
\langle \sigmaa_k(t)-\sigmaa_h(t), E\ww_k(t)-E\uu_k(t)\rangle=
\langle \sigmaa_k(t)-\sigmaa_h(t),E\ww_h(t)-E\uu_h(t)\rangle=0\,.
$$
Therefore,
taking the scalar product of both sides of \eqref{4569} with 
$$
(-\sigmaa_k(t),\zz_k(t))-(-\sigmaa_h(t),\zz_h(t))
$$
and using \eqref{Rkt}, \eqref{4570}, and \eqref{4572}
we obtain
\begin{equation}\label{4573}
\begin{array}{c}
\displaystyle\vphantom{\int_0^t}
\langle \sigmaa_k(t)-\sigmaa_h(t), \ee_k(t)-\ee_h(t)\rangle +  \|\zz_k(t)-\zz_h(t)\|_2^2\le
\\
\displaystyle
\le  
\psi_{hk}(t)
 \{L\int_0^t\varphi_{hk}(s)\,ds + B_T(\tau_h+\tau_k)+
 C^T(\tau_h^{1/2}+\tau_k^{1/2})
\}\,,
\end{array}
\end{equation}
where
$$
\psi_{hk}(t):=\{\|\sigmaa_k(t)-\sigmaa_h(t)\|_2^2+ \|\zz_k(t)-\zz_h(t)\|_2^2\}^{1/2}\,.
$$
By \eqref{boundsC} and \eqref{normC} we obtain $\psi_{hk}(t)\le \beta\,\varphi_{hk}(t)$ and
$$
\alpha\,\varphi_{hk}(t)^2\le \langle \sigmaa_k(t)-\sigmaa_h(t), \ee_k(t)-\ee_h(t)\rangle + 
\|\zz_k(t)-\zz_h(t)\|_2^2\,,
$$
where $\alpha:=\min\{1,\alpha_\C\}$ and $\beta:=\max\{1,2\beta_\C\}$. 
Therefore \eqref{4573} gives
$$
\alpha\varphi_{hk}(t)\le \beta\, L\int_0^t\varphi_{hk}(s)\,ds + \beta\, B^T(\tau_h+\tau_k)+ 
\beta\, C^T(\tau_h^{1/2}+\tau_k^{1/2})\,.
$$
Using Gronwall's inequality we deduce that
$$
\alpha \varphi_{hk}(t)\le\beta  \{B^T(\tau_h+\tau_k)+ 
C^T(\tau_h^{1/2}+\tau_k^{1/2})\}\exp(\beta Lt/\alpha)
$$
for every $t\in [0,T]$ and every $h$ and $k$.
As $\tau_k\to 0$ by \eqref{fineX}, recalling \eqref{4571} we obtain that $\ee_k$ and $\zz_k$ satisfy the Cauchy condition in
$L^\infty([0,T];L^2(\Om;\Mnn))$ and $L^\infty([0,T];L^2(\Om))$, respectively.
\end{proof}

\noindent
{\it Proof of Theorem~\ref{young-main2}. Step 3. Convergence of the interpolations.\/}
By Lemma~\ref{Gronwall} there exist two functions
$\ee\in L^\infty_{loc}([0,+\infty);L^2(\Om;\Mnn))$ and 
$\zz\in L^\infty_{loc}([0,+\infty);L^2(\Om))$ such that for every $T>0$
\begin{eqnarray}
& \label{sek}
\displaystyle
\sup_{0\le t\le T}\|\ee_k(t)-\ee(t)\|_2\to 0\,,
\\
&
\displaystyle
\sup_{0\le t\le T}\|\zz_k(t)-\zz(t)\|_2\to 0\,.
 \label{szk}
\end{eqnarray}
Let $\sigmaa(t):=\C\ee(t)$ and $\zetaa(t):=-V'(\zz(t))$. By \eqref{V'infty} we have also
\begin{eqnarray}
& \label{ssigmak}
\displaystyle
\sup_{0\le t\le T}\|\sigmaa_k(t)-\sigmaa(t)\|_2\to 0\,,
\\
& \label{szetak}
\displaystyle
\sup_{0\le t\le T}\|\zetaa_k(t)-\zetaa(t)\|_2\to 0\,,
\end{eqnarray}
for every $T>0$.

By \eqref{initk} and \eqref{B1005} for every $T>0$ the sequences $\ee_k^\smtr$ and $\zz_k^\smtr$ are bounded in 
$H^1([0,T];L^2(\Om;\Mnn))$ and $ H^1([0,T];L^2(\Om))$,
respectively.
Hence, up to subsequences, we may assume
\begin{eqnarray}
& \label{wektr}
\ee_k^\smtr \wto \hat \ee \quad \text{weakly in } H^1([0,T];L^2(\Om;\Mnn))\,,
\\
&\zz_k^\smtr \wto \hat \zz \quad \text{weakly in } H^1([0,T];L^2(\Om))\,.
 \label{wzktr}
\end{eqnarray}
Moreover, using \eqref{B1005} and the identity
$$
\ee_k^{\smtr}(t)=\ee_k(t)+\int_{[t]_k}^t \dot \ee^\smtr_k(s)\, ds\,,
$$
which holds for every $t\in [0,T]$, we obtain
$$
\|\ee_k^{\smtr}(t)-\ee_k(t)\|_2\le C_T^{1/2}\e^{-1/2}\tau_k^{1/2}\,.
$$
Together with \eqref{sek} and \eqref{wektr} this implies $\hat \ee(t)=\ee(t)$ for every $t\in [0,T]$ and that
\begin{equation}\label{5270}
\sup_{0\le t\le T}\|\ee_k^{\smtr}(t)-\ee(t)\|_2\to 0\,.
\end{equation}
In the same way we deduce that $\hat \zz(t)=\zz(t)$  for every $t\in [0,T]$ and that
\begin{equation}\label{5271}
\sup_{0\le t\le T}\|\zz_k^{\smtr}(t)-\zz(t)\|_2\to 0\,.
\end{equation}
This implies that, without extracting any subsequence,
\begin{eqnarray}
& \label{wektr2}
\ee_k^\smtr \wto \ee \quad \text{weakly in } H^1([0,T];L^2(\Om;\Mnn))\,,
\\
&\zz_k^\smtr \wto \zz \quad \text{weakly in } H^1([0,T];L^2(\Om))\,,
 \label{wzktr2}
\end{eqnarray}
and that $\ee\in H^1([0,T];L^2(\Om;\Mnn))$ and $\zz\in H^1([0,T];L^2(\Om))$.

Since $\div\,\sigmaa_k(t)=0$ in $\Om$ and  $[\sigmaa_k(t)n]=0$ on $\Gamma_1$ by \eqref{3.5}, the strong convergence of $\sigmaa_k(t)$ to $\sigmaa(t)$ in $L^2(\Om;\Mnn)$ yields (ev2)$\!_\e$, taking \eqref{sigmanu}  into account.

By \eqref{dHeNe} and \eqref{subPar} we have
$$
(\dot\pp^\smtr_k(t), \dot\zz^\smtr_k(t))=N_K^\e(\sigmaa_k(t)_D+(\sigma_k^i)_D
-(\sigma_k^{i-1})_D,\zetaa_k(t)+\zeta_k^i-\zeta_k^{i-1})
\qquad\hbox{a.e.\ in\ }\Om
$$
for $t_k^{i-1}\le t<t^i_k$.
Using \eqref{B1005} and the Lipschitz continuity of $N_K^\e$ we obtain
\begin{equation}\label{flowk}
(\dot\pp^\smtr_k(t), \dot\zz^\smtr_k(t))=N_K^\e(\sigmaa_k(t)_D,\zetaa_k(t))
+\hat {\boldsymbol R}_k(t)\qquad\hbox{a.e.\ in\ }\Om\,,
\end{equation}
with $\|\hat {\boldsymbol R}_k(t)\|_2\le \hat C^T\tau_k^{1/2}$, where $ \hat C^T$ is a constant, depending on $\e$, but independent of $t$ and~$k$. By \eqref{ssigmak} and \eqref{szetak} this implies that
\begin{equation}\label{5273}
\sup_{0\le t\le T}
\|(\dot\pp^\smtr_k(t), \dot\zz^\smtr_k(t))-N_K^\e(\sigmaa(t)_D,\zetaa(t))\|_2\to 0
\end{equation}
for every  $T>0$. This shows, in particular, that $\dot\zz^\smtr_k$ converges strongly 
in $L^\infty([0,T];L^2(\Om))$ to a limit that, by
\eqref{5271}, must coincide with $\dot\zz$. Therefore
\begin{equation}\label{5274}
\sup_{0\le t\le T}\|\dot\zz^\smtr_k(t)-\dot\zz(t)\|_2\to 0
\end{equation}
for every $T>0$. Moreover,
since $\pp_k^\smtr(0)=p_0$, \eqref{5273} implies also that there exists 
$\pp\in H^1_{loc}([0,+\infty);L^2(\Om;\MD))$ such that
\begin{equation}
\label{7613}
\sup_{0\le t\le T}
\|\pp^\smtr_k(t)-\pp(t)\|_2\to 0\,,\qquad
\sup_{0\le t\le T}
\|\dot\pp^\smtr_k(t)-\dot\pp(t)\|_2\to 0
\end{equation}
for every $T>0$, and
$$
(\dot\pp(t), \dot\zz(t))=N_K^\e(\sigmaa(t)_D,\zetaa(t))\quad\text{a.e.\ in }\Omega
$$
for a.e.\ $t\in[0,+\infty)$. This proves (ev\^3)$\!_\e$.

As we did for $\ee_k^{\smtr}$ and $\ee_k$, using \eqref{B1005} we can prove that
$$
\|\pp_k^{\smtr}(t)-\pp_k(t)\|_2\le C_T^{1/2}\e^{-1/2}\tau_k^{1/2}\,.
$$
Therefore \eqref{7613} implies
\begin{equation}
\label{spk}
\sup_{0\le t\le T}
\|\pp_k(t)-\pp(t)\|_2\to 0
\end{equation}
for every $T>0$.

For every $t\in[0,T]$ and every $k$ we have
\begin{eqnarray*}
&
E\uu_k(t)=\ee_k(t)+\pp_k(t) \quad\hbox{a.e.\ in }\Om\,,
\\
&
\uu_k(t)=\ww_k(t)\quad\hn\hbox{-a.e.\ in }\Gamma_0\,.
\end{eqnarray*}
As $\ww_k(t)\to\ww(t)$ strongly in $H^1(\Om;\Rn)$, using the Korn-Poincar\'e inequality,
we deduce from \eqref{B1005}, \eqref{sek}, and \eqref{spk}  that there exists
$\uu\in H^1_{loc}([0,+\infty);H^1(\Om;\Rn))$ satisfying (ev1)$\!_\e$ such that
\begin{eqnarray}
&
\label{wu1}
\displaystyle
\sup_{0\le t\le T}\|\uu_k(t)- \uu(t)\|_2\to 0\,,\qquad\sup_{0\le t\le T}\|E\uu_k(t)- E\uu(t)\|_2\to 0\,,
\\
&
\label{wu2}
\displaystyle
\sup_{0\le t\le T}\|\uu_k^\smtr(t)- \uu(t)\|_2\to 0\,,\qquad
\sup_{0\le t\le T}\|E\uu_k^\smtr(t)- E\uu(t)\|_2\to 0\,,
\\
&
\label{wu3}
\displaystyle
\uu_k^\smtr \wto \uu \quad \text{weakly in } H^1([0,T];H^1(\Om;\Rn))\,,
\end{eqnarray}
for every $T>0$ .

\medskip

\noindent
{\it Step 4. Uniqueness of the solution.\/}
Suppose that $(\uu_\e^1, \ee_\e^1,\pp_\e^1, \zz_\e^1)$ and 
$(\uu_\e^2, \ee_\e^2,\pp_\e^2, \zz_\e^2)$ are two solutions of the $\e$-regularized evolution problem with boundary datum $\ww$ and initial condition
$(u_0,e_0,p_0,z_0)$. 
Let $\sigmaa_\e^1(t):=\C\ee_\e^1(t)$, $\sigmaa_\e^2(t):=\C\ee_\e^2(t)$,
$\zetaa_\e^1(t):=-V'(\zz_\e^1(t))$, and $\zetaa_\e^2(t):=-V'(\zz_\e^2(t))$.
By (ev\^3)$\!_\e$ for a.e.\ $t\in[0,+\infty)$ we have
\begin{eqnarray}
&\label{7681}
(\dot\pp_\e^1(t),\dot\zz_\e^1(t))=N^\e_K(\sigmaa_\e^1(t)_D,\zetaa_\e^1(t))
\quad \text{a.e.\ in }\Om\,,
\\
&\label{7682}
(\dot\pp_\e^2(t),\dot\zz_\e^2(t))=N^\e_K(\sigmaa_\e^2(t)_D,\zetaa_\e^2(t))
\quad \text{a.e.\ in }\Om\,.
\end{eqnarray}
Using (ev1)$\!_\e$ we get
\begin{eqnarray*}
&
(-\dot\ee_\e^1(t),\dot\zz_\e^1(t))=N^\e_K(\sigmaa_\e^1(t)_D,\zetaa_\e^1(t))-
(E\dot\uu_\e^1(t),0)
\quad \text{a.e.\ in }\Om\,,
\\
&
(-\dot\ee_\e^2(t),\dot\zz_\e^2(t))=N^\e_K(\sigmaa_\e^2(t)_D,\zetaa_\e^2(t))-
(E\dot\uu_\e^2(t),0)
\quad \text{a.e.\ in }\Om\,.
\end{eqnarray*}
Subtracting term by term we obtain
\begin{equation}\label{7685}
\begin{array}{c}
(-\dot\ee_\e^1(t)+\dot\ee_\e^2(t),\dot\zz_\e^1(t)-\dot\zz_\e^2(t))=
\\
=
N^\e_K(\sigmaa_\e^1(t)_D,\zetaa_\e^1(t))-N^\e_K(\sigmaa_\e^2(t)_D,\zetaa_\e^2(t))
-(E\dot\uu_\e^1(t)-E\dot\uu_\e^2(t),0)
\end{array}
\end{equation}
a.e.\ in $\Om$. Since $N^\e_K$ is $1/\e$-Lipschitz and $V'$ is $M$-Lipschitz, there exists a contant $L_\e$, independent of $t$, such that
\begin{equation}\label{7686}
\|N^\e_K(\sigmaa_\e^1(t)_D,\zetaa_\e^1(t))-N^\e_K(\sigmaa_\e^2(t)_D,\zetaa_\e^2(t))\|_2\le
L_\e\{ \|\ee_\e^1(t)-\ee_\e^2(t)\|_2+ \|\zz_\e^1(t)-\zz_\e^2(t)\|_2\}\,. 
\end{equation}
Integrating by parts and using the equilibrium condition (ev2)$\!_\e$ and the boundary condition $\dot \uu_\e^1(t)=\dot \uu_\e^2(t)=\dot \ww(t)$ $\hn$-a.e.\ on $\Gamma_0$, we obtain
$$
\langle E\dot\uu_\e^1(t)-E\dot\uu_\e^2(t), \sigmaa_\e^1(t)-\sigmaa_\e^2(t)\rangle=0\,. 
$$
Therefore, taking the scalar product of both sides of \eqref{7685} with
$$
(-\sigmaa_\e^1(t)+\sigmaa_\e^2(t),\zz_\e^1(t)-\zz_\e^2(t))
$$
we obtain,  using \eqref{7686},
$$
\begin{array}{c}
\langle \sigmaa_\e^1(t)-\sigmaa_\e^2(t), \dot\ee_\e^1(t)-\dot\ee_\e^2(t)\rangle+
\langle\zz_\e^1(t)-\zz_\e^2(t),  \dot\zz_\e^1(t)-\dot\zz_\e^2(t)\rangle
\le
\\
\le
L_\e\{ \|\ee_\e^1(t)-\ee_\e^2(t)\|_2+ \|\zz_\e^1(t)-\zz_\e^2(t)\|_2\}^2\,,
\end{array}
$$
which gives, thanks to \eqref{boundsC},
$$
\tfrac{d}{dt}\{\QQ(\ee_\e^1(t)-\ee_\e^2(t))+ \tfrac12\|\zz_\e^1(t)-\zz_\e^2(t)\|_2\}\le
\hat L_\e\{\QQ(\ee_\e^1(t)-\ee_\e^2(t))+ \tfrac12\|\zz_\e^1(t)-\zz_\e^2(t)\|_2^2\}\,,
$$
for a suitable constant $\hat L_\e$ independent of $t$.

Since $\QQ(\ee_\e^1(0)-\ee_\e^2(0))=0$  and $\|\zz_\e^1(0)-\zz_\e^2(0)\|_2^2=0$, from Gronwall's inequality we obtain  $\QQ(\ee_\e^1(t)-\ee_\e^2(t))=0$  and 
$\|\zz_\e^1(t)-\zz_\e^2(t)\|_2^2=0$ for every $t\in[0,+\infty)$, hence 
$\ee_\e^1(t)=\ee_\e^2(t)$,  $\zz_\e^1(t)=\zz_\e^2(t)$, 
$\sigmaa_\e^1(t)=\sigmaa_\e^2(t)$, and $\zetaa_\e^1(t)=\zetaa_\e^2(t)$  for every $t\in[0,+\infty)$. {}From \eqref{7681} and
\eqref{7682} we deduce that $\dot\pp_\e^1(t)=\dot\pp_\e^2(t)$ for a.e.\ $t\in [0,+\infty)$. As $\pp_\e^1(0)=\pp_\e^2(0)$, this implies $\pp_\e^1(t)=\pp_\e^2(t)$ 
for every $t\in[0,+\infty)$. The equality $\uu_\e^1(t)=\uu_\e^2(t)$ follows now from the kinematic admissibility condition (ev1)$\!_\e$.
\end{proof}

\begin{remark}\label{rem:convergence}
As the solution of the $\e$-regularized evolution problem is unique,
the proof of Theorem~\ref{young-main2} shows that
it can be approximated by the sequences $(\uu_k,\ee_k,\pp_k,\zz_k)$ and $(\uu_k^\smtr,\ee_k^\smtr,\pp_k^\smtr,\zz_k^\smtr)$ obtained by interpolating the solutions to the incremental problems. The precise convergence properties are given by 
\eqref{sek}--\eqref{szetak}, \eqref{5270}--\eqref{wzktr2}, and \eqref{5274}--\eqref{wu3}\,.
\end{remark}

\begin{remark}\label{rem1d}
When $d=1$
one can express the conditions of Definition~\ref{def:reym} in a simpler way.

We begin by observing that condition (ev2)$\!_\e$ is equivalent to the fact that for every $t\in(0,+\infty)$ the functions
$\ee_\e(t)$ and  $\sigmaa_\e(t)$ are constant on $[-\tfrac12,\tfrac12]$.

Let us write $N_{K}^\e\colon\R{\times}\R\to\R{\times}\R$ componentwise, 
setting $N_K^\e=
(N_{K}^{\e,1}, N_{K}^{\e,2})$. 
We select an arbitrary representative in the equivalence class of
$z_0$. For every $y\in[-\tfrac12,\tfrac12]$ we can consider the solution 
$\psi(\cdot,y)$ of the ordinary differential equation
$$
\partial_t\psi(t,y)= N_{K}^{\e,2}(\sigmaa_\e(t),-V'(\psi(t,y)))\,,\qquad
\psi(0,y)=z_0(y)\,.
$$
As $\sigmaa_\e(t)$ and $V'(\psi(t,y))$ are bounded uniformly with respect to $t$ and $y$,
the same property holds for $\partial_t\psi(t,y)$.
By the dominated convergence theorem we have
$$
\frac{\psi(t+h,\cdot)-\psi(t,\cdot)}{h}\to \partial_t\psi(t,\cdot)
$$
strongly in $L^2([-\tfrac12,\tfrac12])$, so that $t\mapsto \psi(t,\cdot)$ is a $C^1$
solution of the Cauchy problem in $L^2([-\tfrac12,\tfrac12])$
$$
\dot\psii(t)=N_{K}^{\e,2}(\sigmaa_\e(t),V'(\psii(t)))\,,\qquad
\psii(0)=z_0\,.
$$
By condition (ev\^3)$\!_\e$ the function $\zz_\e$ is a solution of the same Cauchy problem. Since the function
$\psi\mapsto N_{K}^{\e,2}(\sigmaa_\e(t),\psi)$ is a Lipschitz continuous function from 
$L^2([-\tfrac12,\tfrac12])$ into $L^2([-\tfrac12,\tfrac12])$, by uniqueness we conclude that for every $t\ge 0$ we have $\zz_\e(t,y)=\psi(t,y)$ for a.e.\ 
$y\in [-\tfrac12,\tfrac12]$.

This shows that for every $t\ge 0$ we can choose  representatives of $\zz_\e(t)$
and $\dot\zz_\e(t)$ in their equivalence classes such that
$\dot\zz_\e(t,y)=\partial_t\zz_\e(t,y)$ for every $y\in [-\tfrac12,\tfrac12]$.  
Using (ev1)$\!_\e$, for these representatives
 condition (ev\^3)$\!_\e$ is equivalent to
\begin{eqnarray}
&\displaystyle
\dot\ee_\e(t)=-\int_{-\frac12}^\frac12 N_{K}^{\e,1}(\sigmaa_\e(t),\zetaa_\e(t,y))\,dy + 
\dot \ww(t,\tfrac12)- \dot \ww(t,-\tfrac12)\,,
\label{oneD1}
\\
&\displaystyle
\dot\zz_\e(t,y)=N_{K}^{\e,2}(\sigmaa_\e(t),\zetaa_\e(t,y))\qquad\hbox{for every }
y\in (-\tfrac12,\tfrac12) \,.\label{oneD2}
\end{eqnarray}
If, in addition, $z_0\in C([-\tfrac12,\tfrac12])$, then 
$\zz_\e\in C([0,+\infty); C([-\tfrac12,\tfrac12]))$, since 
$\sigmaa_\e$ is continuous and  $\zz_\e(\cdot,y)$ is a solution of an 
ordinary differential equation with initial condition $z_0(y)$ depending 
continuously on~$y$.
\end{remark}

\end{section}

\begin{section}{Approximable quasistatic evolution}\label{appr}

In this section we state and prove the main results of the paper on approximable quasistatic evolutions.

\subsection{Definition and properties}
We recall that a function $\vv$ from $[0,+\infty)$ into the dual $Y$ of a Banach space is weakly$^*$ left-continuous if 
$$
\vv(s)\wto \vv(t)\quad\hbox{weakly}^*\hbox{ in }Y
$$
as $s\to t$, with $s\le t$. 
We recall that a system $\muu\in SGY([0,+\infty),\ol\Om;\MD{\times}\R)$ is said to be
{\it weakly$^*$ left-continuous\/} if for every finite sequence $t_1,\dots,t_m$ in $[0,+\infty)$ with 
$t_1 <\dots < t_m$ the following  continuity property holds:
\begin{equation} \label{leftlim}
 \muu_{s_1\dots s_m}\wto \muu_{t_1\dots t_m} 
 \quad\hbox{weakly}^*\hbox{ in }GY(\ol\Om; (\MD{\times}\R)^m)
\end{equation}
as $s_i\to t_i$, with $s_i\in [0,+\infty)$ and $s_i\le t_i$.

The notion of approximable quasistatic evolution is made precise by the following definition.

\begin{definition}\label{def:qsym}
Let  $\ww\in H^1_{loc}([0,+\infty);H^1(\Om;\Rn))$, $u_0\in H^1(\Om)$, 
$e_0\in L^2(\Om; \Mnn)$,  $p_0\in L^2(\Om; \MD)$, and $z_0\in L^2(\Om)$. 
Suppose that the kinematic admissibility condition (ev1)$\!_0$ and the equilibrium condition (ev2)$\!_0$ of Theorem \ref{young-main2} are satisfied.
An {\it approximable quasistatic evolution\/} with boundary datum $\ww$
and initial condition $(u_0, e_0, p_0, z_0)$ is a triple $(\uu,\ee,\muu)$, with 
\begin{eqnarray*}
&\uu\in L^\infty_{loc}([0,+\infty);BD(\Om))\,, \qquad
 \ee\in L^\infty_{loc}([0,+\infty);L^2(\Om;\Mnn))\,,
\\
&\muu\in SGY([0,+\infty), \ol\Omega; \MD{\times}\R)\,,
\end{eqnarray*}
and $\uu$, $\ee$, $\muu$ weakly$^*$ left-continuous,
for which  there exist a positive sequence $\e_k\to 0$ and
 a set $\Theta\subset[0,+\infty)$, containing $0$ and with $\LL^1([0,+\infty)\setminus \Theta)=0$,  such that 
the solutions
$(\uu_{\e_k},\ee_{\e_k},\pp_{\e_k}, \zz_{\e_k})$ of the $\e_k$-regularized evolution problems with boundary datum $\ww$ and initial condition $(u_0, e_0, p_0, z_0)$  satisfy\begin{eqnarray}
&\uu_{\e_k}(t)\wto \uu(t)\qquad\text{weakly$^*$ in }BD(\Omega)\,,
\label{u55}\\
&\ee_{\e_k}(t)\wto \ee(t)\qquad\text{weakly in }L^2(\Omega;\Mnn) \label{e55}
\end{eqnarray}
for every $t\in\Theta$, and 
\begin{equation}\label{p55}
\big(\deltaa_{(\pp_{\e_k},\zz_{\e_k})}\big)_{t_1\dots t_m}\wto \muu_{t_1\dots t_m}^{}\qquad\text{weakly$^*$ in }GY(\ol\Om;(\MD{\times}\R)^m)
\end{equation}
for every finite sequence $t_1,\dots,t_m$ in $\Theta$ with $t_1<\dots<t_m$.

If, in addition, $\muu=\deltaa_{(\pp,\zz)}$ for some functions
$\pp\in L^\infty_{loc}([0,+\infty); M_b(\ol\Om;\MD))$ and 
$\zz\in L^\infty_{loc}([0,+\infty); L^1(\Om))$, we say that $(\uu,\ee,\pp,\zz)$ is a {\em  spatially regular} approximable quasistatic evolution.
\end{definition}

The first estimates for approximable quasistatic evolutions are proved in the next lemma.

\begin{lemma}\label{lm:prime stime}
Let  $\ww$, $u_0$, $e_0$, $p_0$, $z_0$, $(\uu,\ee,\muu)$, $\e_k$, and 
$(\uu_{\e_k},\ee_{\e_k},\pp_{\e_k}, \zz_{\e_k})$ be as in Definition~\ref{def:qsym}. Then for every $T>0$ there exists a constant $C_T<+\infty$ such that 
\begin{eqnarray}
&\displaystyle
\sup_{t\in [0,T]}\|\ee_{\e_k}(t)\|_2\leq C_T\,, \quad
\D_{\!H}(\pp_{\e_k}, \zz_{\e_k};0,T)\leq C_T\,, \label{primest2}
\\
&\displaystyle\e_k\int_0^T\|\dot \pp_{\e_k}(t)\|_2^2\, dt\leq C_T\,, \quad\text{and}\quad 
\e_k\int_0^T\|\dot \zz_{\e_k}(t)\|_2^2\, dt\leq C_T\,. \label{primest1}
\end{eqnarray}
\end{lemma}
\begin{proof}
By the energy equality (ev4)$\!_\e$ and by \eqref{diss-reg} for every $T>0$ we have 
\begin{eqnarray}
&\displaystyle\QQ(\ee_{\e_k}(T))+
\D_{\!H}(p_{\e_k}, z_{\e_k}; 0,T) +
\V(\zz_{\e_k}(T))\nonumber\\
&\displaystyle+\e_k\int_{0}^{T}\|\dot \pp_{\e_k}(t)\|_2^2\, dt + \e_k
\int_{0}^{T} \|\dot \zz_{\e_k}(t)\|_2^2\, dt
=\label{dis-en}
 \\
&\displaystyle= \QQ(e_0)+\V(z_0)+
\int_{0}^{T} \langle\sigmaa_{\e_k}(t), E\dot \ww(t)\rangle\, dt\,,\nonumber
\end{eqnarray}
where $\sigmaa_{\e_k}(t):=\C\ee_{\e_k}(t)$. 

The proof can now be concluded as in Step 2 of the proof of Theorem~\ref{young-main2}.
\end{proof}

To study the properties of approximable quasistatic evolutions we need to introduce some definitions.
Given $\muu\in SGY([0,+\infty),\ol\Om; \MD{\times}\R)$, its dissipation $\D_{\!H}(\muu;a,b)$ on the time interval $[a,b]\subset [0,+\infty)$ is defined as 
\begin{equation}\label{nudiss}
\sup \sum_{i=1}^k 
\langle H(\xi_i-\xi_{i-1},\theta_i-\theta_{i-1}), 
\muu_{t_0t_1\dots t_k}(x,\xi_0,\theta_0,\dots,\xi_k,\theta_k,\eta) \rangle\,,
\end{equation}
where the supremum is taken over all finite families $t_0,t_1,\dots,t_k$ such that 
$a=t_0<t_1<\dots<t_k=b$.
As in the case of the variation, we have
\begin{equation}\label{nudiss2}
\D_{\!H}(\muu;a,b) = 
\sup \sum_{i=1}^k \langle H(\xi_i-\xi_{i-1},\theta_i-\theta_{i-1}),
\muu_{t_{i-1}t_i}(x,\xi_{i-1},\theta_{i-1},\xi_i,\theta_i,\eta)\rangle\,,
\end{equation}
where the supremum is taken over all finite families $t_0,t_1,\dots,t_k$ such that 
$a=t_0<t_1<\dots<t_k=b$.

If $\muu=\deltaa_{(\pp,\zz)}$ for some $\pp\colon[0,+\infty)\to M_b(\ol\Om;\MD)$
and $\zz\colon[0,+\infty)\to M_b(\ol\Om)$, we have
\begin{equation}\label{eq=diss}
\D_{\!H}(\muu;a,b)= \D_{\!H}(\pp,\zz;a,b)\,,
\end{equation}
where $\D_{\!H}(\pp,\zz;a,b)$ is defined in \eqref{diss}.

We also introduce the homogeneous functions
$\{V\}\colon \R{\times}\R\to\R$ and $\{V'\}\colon \R{\times}\R\to\R$ defined by
$$
\{V\}(\theta,\eta):=
\begin{cases}
\eta\, V(\theta/\eta) & \text{if }\eta>0\,,
\\
V^\infty(\theta) &  \text{if }\eta\le 0 \,,
\end{cases}
$$
and
$$
\{V'\}(\theta,\eta):=
\begin{cases}
\eta\, V'(\theta/\eta) & \text{if }\eta>0\,,
\\
0 &  \text{if }\eta\le 0\,.
\end{cases}
$$

According to \cite[Definition~2.15]{DM-DeS-Mor-Mor-1},
for every $\mu\in GY(\ol\Om; \MD{\times}\R)$  we can define
the measure $\pi_{\ol\Om}(\{V'\}\mu)\in M_b(\ol\Om)$ as
\begin{equation}\label{7220}
\langle\varphi,\pi_{\ol\Om}(\{V'\}\mu) \rangle = 
\langle \varphi(x) \{V'\}(\theta,\eta),\mu(x,\xi,\theta,\eta)\rangle
\end{equation}
for every $\varphi\in C(\ol\Om)$.

\begin{lemma}\label{lm:V'nu}
For every $\mu\in GY(\ol\Om; \MD{\times}\R)$ we have
$\pi_{\ol\Om}(\{V'\}\mu)\in L^{\infty}(\Om)$ and, moreover,
$\|\pi_{\ol\Om}(\{V'\}\mu)\|_{\infty}\le\|V'\|_\infty$.
\end{lemma} 

\begin{proof}
By (\ref{mupos}) we have
$$
| \langle \varphi(x) \{V'\}(\theta,\eta),\mu(x,\xi,\theta,\eta)\rangle|\le
\|V'\|_\infty\langle |\varphi(x)|\eta,\mu(x,\xi,\theta,\eta)\rangle
=\|V'\|_\infty\|\varphi\|_1
$$
for every $\varphi\in C(\ol\Om)$. The conclusion follows.
\end{proof}

The main properties of approximable quasistatic evolutions are proved in the following theorem, where $\K(\Om)$ is the convex set defined in~\eqref{KOm}.

\begin{theorem}\label{Thm55}
Let  $\ww$, $u_0$, 
$e_0$,  $p_0$, and $z_0$ be as in Definition~\ref{def:qsym}.
Let $(\uu,\ee,\muu)$ be an approximable quasistatic evolution with boundary datum
$\ww$ and initial condition $(u_0, e_0, \mu_0)$, and let 
\begin{equation}\label{7221}
(\pp(t), \zz(t))=\bary(\muu_t)\,, \quad\sigmaa(t):=\C \ee(t)\,, \quad \zetaa(t):=-\pi_{\ol\Om}(\{V'\}\muu_t)\,.
\end{equation}
Then the following conditions are satisfied:
\begin{itemize}
\smallskip
\item[(ev1)] {\it weak kinematic admissibility:\/} for every $t\in(0,+\infty)$ 
$$
\begin{array}{c}
E\uu(t)=\ee(t)+\pp(t) \quad \text{ in }\Om\,,
\smallskip
\\
\pp(t)=(\ww(t)-\uu(t))\odot n \hn
\quad \text{on } \Gamma_0\,;
\end{array}
$$
\smallskip 
\item[(ev2)] {\it equilibrium condition:\/} for every $t\in(0,+\infty)$ 
$$
\begin{array}{c}
\div\,\sigmaa(t)=0 \text{ in } \Om\,, \quad
[\sigmaa(t) n]=0 \text{ on } \Ga_1\,;
\end{array}
$$
\smallskip
\item[(ev3)] {\it stress constraint:\/} for every $t\in(0,+\infty)$ 
$$ 
 (\sigmaa_D(t),\zetaa(t))\in \K(\Om)\,;
$$
\smallskip
\item[(ev4)] {\it energy inequality:\/}  
for every $T\in(0,+\infty)$  
\begin{eqnarray}
& \displaystyle 
\QQ(\ee(T)) + \D_{\!H}(\muu;0,T) + 
\langle \{V\}(\theta,\eta), \muu_{T}(x, \xi, \theta, \eta)\rangle 
\le \nonumber \smallskip
\\
& \displaystyle
\le \QQ(e_0)+\V(z_0)+\int_{0}^{T} \langle\sigmaa(t), E\dot \ww(t)\rangle\, dt\,. \nonumber
\end{eqnarray}
\end{itemize}
\end{theorem}

\begin{remark}\label{RemThm55}
If $(\uu,\ee, \pp,\zz)$ is a spatially regular approximable quasistatic evolution, then 
$\zetaa(t)=-V'(\zz(t))$ by \eqref{mup}, \eqref{7220}, and \eqref{7221}. Therefore the energy inequality (ev4) reduces to 
\begin{equation}\label{ineqreduced}
\begin{array}{c}
\displaystyle\QQ(\ee(T))+\D_{\!H}(\pp,\zz; 0,T)+\V(\zz(T))\leq
\smallskip\\
\displaystyle\leq \QQ(e_0)+\V(z_0)+\int_{0}^{T} \langle\sigmaa(t), E\dot \ww(t)\rangle\, dt\,
\end{array}
\end{equation}
thanks to \eqref{eq=diss}.
\end{remark}

\begin{proof}[Proof of Theorem~\ref{Thm55}]
Let $\Theta$, $\e_k$,  and $(\uu_{\e_k},\ee_{\e_k},\pp_{\e_k},\zz_{\e_k})$ be as in Definition~\ref{def:qsym}.
By \eqref{weakbar} and \eqref{u55}--\eqref{p55}  we have
\begin{eqnarray*}
&\uu_{\e_k}\!(t)\wto \uu(t)\qquad\text{weakly$^*$ in }BD(\Omega)\,,
\\
&\ee_{\e_k}\!(t)\wto \ee(t)\qquad\text{weakly in }L^2(\Omega;\Mnn)\,,
\\
&\pp_{\e_k}\!(t)\wto \pp(t)\qquad\text{weakly$^*$ in }M_b(\ol\Omega;\MD)
\end{eqnarray*}
for every $t\in\Theta$.
By applying \cite[Lemma~2.1]{DM-DeS-Mor} we  prove that (ev1)$\!_\e$ imples (ev1) for every $t\in\Theta$. These equalities are then extended to every $t\in(0,+\infty)$ by left-continuity.

For every $t\in\Theta$ the equilibrium condition (ev2) follows from the equilibrium condition (ev2)$\!_\e$ satisfied by $\sigmaa_{\!\e_k}\!(t)$ and from \eqref{sigmanu}. The same result can be obtained for every $t\in(0,+\infty)$ by left-continuity.

By the modified stress constraint (ev3)$\!_\e$,
for every $k$ and a.e.\ $t\in(0,+\infty)$ we have that
\begin{equation}\label{AkEulereq-10}
(\sigmaa_{\!\e_k}\!(t)_D-\e_k\dot \pp_{\e_k}\!(t),
\zetaa_{\e_k}(t)-\e_k\dot \zz_{\e_k}\!(t))\in \K(\Om)\,.
\end{equation}
For every $t\in\Theta$ we have that $\delta_{(\pp_{\e_k}\!(t),\zz_{\e_k}\!(t))}\wto\muu_t$ weakly$^*$ in $GY(\ol\Om;\MD{\times}\R)$, which, by \eqref{weakbar} and  Lemma~\ref{lm:V'nu}, implies that $\zetaa_{\e_k}\!(t)\wto \zetaa(t)$ weakly$^*$ in $L^\infty(\Om)$.
By (\ref{primest1}) the sequences $\e_k\|\dot \pp_{\e_k}\!(\cdot)\|_2$ and
$\e_k\|\dot \zz_{\e_k}\!(\cdot)\|_2$ converge to $0$ strongly in $L^2_{loc}([0,+\infty))$, 
hence there exists a sequence $k_j\to\infty$ such that $\e_{k_j}\dot \pp_{\e_{k_j}}\!(t)\to 0$ in $L^2(\Om;\MD)$ and
$\e_{k_j}\dot \zz_{\e_{k_j}}\!(t)\to 0$ in $L^2(\Om)$ 
for a.e.\ $t\in[0,+\infty)$.
This implies
\begin{equation}\label{AkEulereq-10f}
(\sigmaa(t), \zetaa(t))\in \K(\Om)
\end{equation}
for a.e.\ $t\in(0,+\infty)$.

The left continuity of $\muu$, together with Lemma~\ref{lm:V'nu}, implies that the function $\zetaa$ is weakly$^*$ left-continuous in $L^\infty(\Om)$. As $\sigmaa$ is left-continuous too, by approximation (\ref{AkEulereq-10f}) holds for every $t\in(0,+\infty)$.
This concludes the proof of (ev3).

It remains to prove the energy inequality (ev4).  
By \eqref{e55} and by the lower semicontinuity of $\QQ$, for every $T\in\Theta$  we have
\begin{equation}\label{Qlsc-10}
\QQ(\ee(T)) \le\,\liminf_{k\to\infty}\,\QQ(\ee_{\e_k}(T))\,, 
\end{equation}
and, by left-continuity, for every $T\in [0,+\infty)$ we have
\begin{equation}\label{Qlsc2-10}
\QQ(\ee(T)) \le \liminf_{\genfrac{}{}{0pt}2{\scriptstyle S\to T}
{\scriptstyle S\le T}} \QQ(\ee(S))\,. 
\end{equation}

As the dissipation is lower semicontinuous (see \cite[Theorem~8.11]{DM-DeS-Mor-Mor-1}), using (\ref{p55}) and the left-continuity, and taking into account \eqref{eq=diss}, one can check that for every $T\in\Theta$
\begin{equation}\label{dislsc-10}
\D_{\!H}(\muu; 0,T) \le\,\liminf_{k\to\infty}\,\D_{\!H}(\pp_{\e_k},\zz_{\e_k}; 0,T)\,, 
\end{equation}
and for every $T\in (0,+\infty)$
\begin{equation}\label{dislsc2-10}
\D_{\!H}(\muu; 0,T) \le \liminf_{\genfrac{}{}{0pt}2{\scriptstyle S\to T}
{\scriptstyle S\le T}} \D_{\!H}(\muu;0,S)\,. 
\end{equation}
Moreover, for every $T\in\Theta$ we have 
\begin{equation}\label{Vnu-10}
\V(\zz_{\e_k}(T))
\to \langle \{V\}(\theta,\eta),\muu_T(x,\xi,\eta,\theta)\rangle
\end{equation}
as $k\to\infty$, and for every $T\in(0,+\infty)$
\begin{equation}\label{Vnus-10}
\langle \{V\}(\theta,\eta),\muu_S(x,\xi,\eta,\theta)\rangle
 \to  \langle \{V\}(\theta,\eta),\muu_T(x,\xi,\eta,\theta)\rangle
\end{equation}
as $S\to T$, with $S\le T$.

Finally, by \eqref{e55} we obtain
$$
\langle \sigmaa_{\!\e_k}\!(t), E\dot \ww(t)\rangle \to \langle \sigmaa(t), E\dot \ww(t)\rangle 
$$
for a.e.\ $t\in[0,+\infty)$. Thanks to \eqref{primest2} we can apply the dominated convergence theorem and we obtain
\begin{equation}\label{INT-10}
\int_0^T\langle \sigmaa_{\!\e_k}\!(t), E\dot \ww(t)\rangle\, dt \to
\int_0^T \langle \sigmaa(t), E\dot \ww(t)\rangle \, dt\,.
\end{equation}

Combining (\ref{Qlsc-10}), (\ref{dislsc-10}), (\ref{Vnu-10}), and (\ref{INT-10}), we can pass to the limit in the energy equality (ev4)$\!_\e$ satisfied by $(\uu_{\e_k},\ee_{\e_k},\pp_{\e_k},\zz_{\e_k})$ and we obtain the energy inequality (ev4) for every $T\in\Theta$. The result for every $T\in[0,+\infty)$ can be obtained by left-continuity thanks to (\ref{Qlsc2-10}), (\ref{dislsc2-10}), and (\ref{Vnus-10}).
\end{proof}

\subsection{The existence result}
We now prove the main existence result of the paper.

\begin{theorem}\label{young-main3}
Let  $\ww\in H^1_{loc}([0,+\infty);H^1(\Om;\Rn))$, $u_0\in H^1(\Om)$, 
$e_0\in L^2(\Om; \Mnn)$,  $p_0\in L^2(\Om; \MD)$, and $z_0\in L^2(\Om)$.  
Suppose that the kinematic admissibility condition $\rm(ev1)_0$ and the equilibrium condition $\rm(ev2)_0$ of Theorem \ref{young-main2} are satisfied.
Then, there exists an approximable quasistatic evolution  with boundary datum $\ww$
and initial condition $(u_0, e_0, p_0, z_0)$.
\end{theorem}

\begin{proof}
Let us fix a positive sequence $\e_k\to 0$.
For every $k$ let $(\uu_{\e_k},\ee_{\e_k},\pp_{\e_k}, \zz_{\e_k})$ be the
solution of the $\e_k$-regularized evolution problem with boundary datum $\ww$ and initial condition $(u_0, e_0, p_0, z_0)$, and let $\sigmaa_{\e_k}(t):=\C \ee_{\e_k}(t)$.
{}From the energy equality (ev4)$\!_\e$ for every $T>0$ we obtain
$$
\QQ(\ee_{\e_k}(T)) + \D_{\!H}(\pp_{\e_k}, \zz_{\e_k};0,T)
+\V( \zz_{\e_k}) \le 
 \QQ(e_0) + \V(z_0)
 + \int_{0}^{T} \langle\sigmaa_{\!\e_k}\!(t), E\dot \ww(t)\rangle\, dt\,.
$$

By \eqref{pvar} and (\ref{gammaM}) we have
\begin{equation}\label{dhvbdd}
\begin{array}{c}
\D_{\!H}(\pp_{\e_k}, \zz_{\e_k};0,T)
+\V( \zz_{\e_k})-\V(z_0)\ge \smallskip
\\
\ge
  C^K_V {\rm Var}(\pp_{\e_k}, \zz_{\e_k};0,T)=
  C^K_V {\rm Var}(\deltaa_{(\pp_{\e_k}, \zz_{\e_k})};0,T)\,.
\end{array}
\end{equation}
Therefore, 
$$
\QQ(\ee_{\e_k}(T)) \le  \QQ(e_0) +
\int_0^T \langle\sigmaa_{\!\e_k}\!(s), E\dot \ww(s)\rangle\, ds\,.
$$
Arguing as in the proof of (\ref{b1005}) and \eqref{B1005}, we obtain that there exists a constant $C_T$, depending on $T$ but independent of $t$ and $k$, such that
\begin{eqnarray}
&\label{eq10.5}
\|\ee_{\e_k}\!(t)\|_2\le C_T\,, \qquad \|\sigmaa_{\!\e_k}\!(t)\|_2\le C_T\,,
\\
&\label{eq10.6}
{\rm Var}(\deltaa_{(\pp_{\e_k}, \zz_{\e_k})}; 0,T)={\rm Var}(\pp_{\e_k}, \zz_{\e_k}; 0,T) \le C_T\,,
\end{eqnarray}
for every $t\in [0,T]$ and for every $k$.
Inequality \eqref{eq10.6}, together with the initial conditions,
implies that  $\|\pp_{\e_k}(t)\|_1\le \|p_0\|_1+\|z_0\|_1+C_T$ and
$\|\zz_{\e_k}(t)\|_1\le \|p_0\|_1+\|z_0\|_1+C_T$ for every $t\in[0,T]$.
 
Using \eqref{seminorm} and the kinematic admissibility condition (ev1)$\!_\e$ of Definition~\ref{def:reym},
from the estimate of $\|\pp_{\e_k}(t)\|_1$ and from (\ref{eq10.5}) we deduce that for every $T>0$ the sequence $\uu_{\e_k}\!(t)$ is bounded in
$BD(\Om)$ uniformly with respect to $t\in[0,T]$ and $k$.  

By applying Helly's theorem for systems of generalized Young measures (see \cite[Theorem~8.10]{DM-DeS-Mor}, together with a standard diagonal argument, we construct a subsequence, still denoted $\e_k$, a set $\Theta\subset[0,+\infty)$, 
containing $0$ and with $[0,T]\setminus\Theta$ at most countable, 
and a left continuous $\muu\in SGY([0,+\infty),\ol\Om ;\MD{\times}\R)$, with
\begin{eqnarray}
& {\rm Var}(\muu;0,T)\le C_T\,,\label{bv0}
\\
& \|\muu_t\|_*\le \|p_0\|_1+\|z_0\|_1+C_T \quad \hbox{for every }t\in[0,T] \,,\label{bn0}
\end{eqnarray}
such that
\begin{equation}\label{conv-10}
(\deltaa_{(\pp_{\e_k}, \zz_{\e_k})})_{t_1\dots t_m} \wto\muu_{t_1\dots t_m} \quad\hbox{weakly}^*\hbox{ in }GY(\ol\Om ;(\MD{\times}\R)^m) 
\end{equation}
for every finite sequence $t_1,\dots,t_m$ in $\Theta$ with $ t_1 <\dots < t_m$. For every $t\in [0,+\infty)$ we set $(\pp(t), \zz(t)):=\bary(\muu_t)$.

Let us fix $t\in\Theta$. Since $\uu_{\e_k}(t)$ and $\ee_{\e_k}(t)$ are bounded in $BD(\Om)$ and $L^2(\Om;\Mnn)$,  there exist an increasing sequence
$k_j$  (possibly depending on $t$) and two functions $\uu(t)\in BD(\Om)$ and 
$\ee(t)\in L^2(\Om;\Mnn)$ such that 
$\uu_{\e_{k_j}}\!(t)\wto \uu(t)$ weakly$^*$ in $BD(\Om)$ and
$\ee_{\e_{k_j}}\!(t)\wto \ee(t)$ weakly in $L^2(\Om;\Mnn)$. 
As in the proof of Theorem~\ref{Thm55} we can show that $\uu(t)$, $\ee(t)$, $\pp(t)$, and $\ww(t)$ satisfy the weak kinematic admissibility condition~(ev1).

By the equilibrium condition (ev2)$\!_\e$ we have
$$
\div\,\sigmaa_{\!\e_k}\!(t)=0 \text{ in } \Om \quad\text{and}\quad
[\sigmaa_{\e_k}(t) n]=0 \text{ on } \Ga_1\,.
$$
Using the weak definition of divergence and normal trace \eqref{sigmanu}, we can pass to the limit and we obtain
\begin{equation}\label{eq-sig}
\div\,\sigmaa(t)=0 \text{ in } \Om \quad\text{and}\quad
[\sigmaa(t) n]=0 \text{ on } \Ga_1\,.
\end{equation}
For every $\varphi\in H^1(\Om;\Rn)$ with $\varphi=0$ $\hn$-a.e.\ on $\Gamma_0$ we have
$$
\QQ(\ee(t)+ E \varphi)- \QQ(\ee(t))=
\langle \sigmaa(t), E\varphi\rangle +\QQ(E\varphi)
= \QQ(E\varphi)\,,
$$
where the last equality follows from \eqref{sigmanu} and  (\ref{eq-sig}).
Therefore
\begin{equation}\label{secpb-10}
\QQ(\ee(t)) \le
 \QQ(\ee(t)+ E \varphi)\,.
\end{equation}
If $\hat k_j$ is another sequence such that $\uu_{\e_{\hat k_j}}\!(t)\wto \hat \uu(t)$ weakly$^*$ in $BD(\Om)$ and $\ee_{\e_{\hat k_j}}\!(t)\wto \hat \ee(t)$ weakly in $L^2(\Om;\Mnn)$, then $\hat\uu(t)$, $\hat\ee(t)$, $\pp(t)$, and $\ww(t)$ satisfy the weak kinematic admissibility condition~(ev1) and
\begin{equation}\label{secpb-10h}
\QQ(\hat \ee(t)) \le
 \QQ(\hat \ee(t)- E \varphi)
\end{equation}
for every $\varphi\in H^1(\Om;\Rn)$ with $\varphi=0$ $\hn$-a.e.\ on $\Gamma_0$. 

{}From the weak kinematic admissibility condition~(ev1) it follows that $\hat \uu(t)-
\uu(t)=0$ $\hn$-a.e.\ on $\Gamma_0$, and $\hat \ee(t)-\ee(t)=E\hat \uu(t)-E\uu(t)$ a.e.\ on $\Om$, hence  $\hat \uu(t)-\uu(t)\in  H^1(\Om;\Rn)$. If we take $\varphi=\frac12(\hat \uu(t)-\uu(t))$ in (\ref{secpb-10}) and (\ref{secpb-10h}), by adding the inequalities we obtain 
$$
\QQ(\ee(t))+\QQ(\hat \ee(t))\le 2 \QQ(\textstyle\frac12 (\ee(t)+\hat \ee(t)))\,,
$$
which implies that $\hat \ee(t)=\ee(t)$ a.e.\ on $\Om$ by the strict convexity of $\QQ$. Therefore $\hat \uu(t)=\uu(t)$  a.e.\ on $\Om$, since $E\hat\uu(t)=E\uu (t)$ a.e.\ on $\Om$ and $\hat\uu (t)=\uu (t)$ $\hn$-a.e.\ on $\Ga_0$.
As these limits do not depend on the subsequences, the convergence results hold for the whole sequences. This proves \eqref{u55} and~\eqref{e55}. Moreover, (\ref{eq10.5}) and the uniform estimate on $\uu_{\e_k}\!(t)$ in $BD(\Om)$ give
\begin{equation}\label{BOUND-10}
\sup_{t\in\Theta\cap[0,T]}\|\ee(t)\|_2<+\infty\qquad\text{and}\qquad \sup_{t\in\Theta\cap[0,T]}\|E\uu(t)\|_1<+\infty
\end{equation}
for every $T>0$

We now show that for every $t\in\Theta$ we have
\begin{equation}\label{leftlimXX-10}
\begin{array}{c}
\ee(s)\wto \ee(t)\quad\hbox{weakly}\hbox{ in }L^2(\Om ;\Mnn)\,,\\
\uu(s)\wto \uu(t)\quad\hbox{weakly}^*\hbox{ in }BD(\Om)\, ,
\end{array}
\end{equation}
as $s\to t$, with $s\in\Theta$  and $s\le t$. By (\ref{BOUND-10}), for every $t\in\Theta$ there exist $\hat e\in L^2(\Om ;\Mnn)$, $\hat u\in BD(\Om)$, and a sequence $s_j\in\Theta$, with $s_j\to t$ and $s_j\le t$, such that 
$\ee(s_j)\wto \hat e$ weakly in $L^2(\Om ;\Mnn)$ and 
$\uu(s_j)\wto \hat u$ weakly$^*$ in $BD(\Om)$. 
As each $\ee(s_j)$ satisfies the minimality property (\ref{secpb-10}), we deduce 
\begin{equation}\label{secpb55-10}
\QQ(\hat e) \le \QQ(\hat e-E \varphi)
\end{equation}
for every $\varphi\in H^1(\Om;\Rn)$ with $\varphi=0$ $\hn$-a.e.\ on $\Gamma_0$.  Moreover by left-continuity we have $\pp(s_j)\wto \pp(t)$ weakly$^*$ in 
$M_b(\ol\Om; \MD)$, which implies that
$\hat u$, $\hat e$, $\pp(t)$, and $\ww(t)$ satisfy the weak kinematic admissibility condition (ev1) of Theorem~\ref{Thm55}. 
As before we can take $\varphi:=\frac{1}{2}(\hat u-\uu(t))$ 
as a test function in (\ref{secpb-10}) 
and (\ref{secpb-10h}) to deduce from the strict convexity of 
$\QQ$ that $\hat e=\ee(t)$ and, in turn, 
$\hat u=\uu(t)$.
Since the limit is independent of the sequence $s_j$,  we have proved (\ref{leftlimXX-10}).

The same argument shows that for every $t\in[0,T]$ we can define in a unique way 
$\ee(t)\in L^2(\Om;\Mnn)$ and $\uu(t)\in BD(\Om)$ such that (\ref{leftlimXX-10}) holds as $s\to t$, with $s\in [0,T]$ and $s\le t$.
\end{proof}

\begin{remark}\label{rm:hom}
Assume that $\ww$, $\Ga_0$, $u_0$, $e_0$, $p_0$,  $z_0$ satisfy 
\eqref{bd-hom}--\eqref{initial-hom3}. Then by Proposition~\ref{prop:sp-hom} the solutions of the $\e$-regularized evolution problems are spatially regular. Therefore
every approximable quasistatic evolution with boundary datum $\ww$ and initial condition
$(u_0, e_0, p_0, z_0)$ is spatially regular and has the special form
\begin{equation}\label{special}
\uu(t,x)=\xi(t)x\,,\quad \ee(t,x)=\xi^e(t)\,,\quad\pp(t,x)=\xi^p(t)\,,\quad \zz(t,x)=\theta(t)\,,
\end{equation}
with $\xi^e\colon[0,+\infty)\to \Mnn$, $\xi^p\colon[0,+\infty)\to \MD$, $\theta\colon[0,+\infty)\to \R$ left-continuous and 
\begin{equation}\label{xip-bis}
\xi^e(t)+\xi^p(t)=\xi^s(t)
\end{equation}
for every $t\in [0,+\infty)$, where $\xi^s(t)$ is the symmetric part of $\xi(t)$.
\end{remark}

\end{section}

\begin{section}{The case of simple shear}\label{simple-shear}

We analyze in this section the case of simple shear in the Dirichlet-Periodic case (DP) for an isotropic material in dimension $d\ge 2$ with shear modulus $\mu$ (see \eqref{iso}). We consider $\ol x$-periodic solutions  with boundary data of the form
\begin{equation}\label{bc}
\ww(t, x_1, \ol x) := \sqrt2 \, \ww^R( t,  x_1) e_2
\end{equation}
and initial conditions of the form
\begin{equation}\label{ic}
\begin{array}{c}
u_0(0,x_1, \ol x):=\sqrt2 \,u_0^R(x_1)e_2\,, \quad 
e_0(0,x_1, \ol x):=M(e_0^R)\,,
\\
p_0(0,x_1, \ol x):=M(p_0^R(x_1))\,, \quad
z_0(0,x_1, \ol x):= z_0^R(x_1)\,, 
\end{array}
\end{equation}
where $\ww^R\in H^1_{loc}([0,+\infty); H^1([-\frac12 , \frac12 ]))$, 
$u_0^R\in H^1([-\frac12,\frac12])$,
$e_0^R\in\R$, $p_0^R$, $z_0^R\in L^2([-\frac12,\frac12])$, $e_2$ is the
second element of the canonical basis of $\Rn$, and $M\colon \R \to \MD$ is the 
linear isometry defined by
\begin{equation}\label{M-oper}
M(\alpha):= \left(\begin{array}{cccc}0 & \frac{\alpha}{\sqrt2} & \cdots & 0 
\\\frac{\alpha}{\sqrt2} & 0 & \cdots & 0
\\ \vdots & \vdots & \ddots & \vdots 
\\0 & 0 & \cdots & 0\end{array}\right)
\,,
\end{equation}
which connects the amount of a pure shear with its strain.
 
We assume that the elastic domain $K$ depends on the stress $\sigma$
 only through its euclidean 
norm $|\sigma|$. More precisely, we assume that there exisits a closed convex set 
$K^{\!R}\subset \R{\times}\R$ such that
\begin{equation}\label{red-K}
K=\{(\sigma , \zeta) \in \MD {\times} \R : (|\sigma|,\zeta)\in K^{\!R} \}\,;
\end{equation}
it is not restrictive to assume that
\begin{equation}\label{plusminus} 
 (\alpha, \zeta)\in K^{\!R} \quad \iff \quad (-\alpha, \zeta)\in K^{\!R}\,.
 \end{equation}

\begin{theorem}\label{thm:ss}
Assume \eqref{bc}, \eqref{ic}, \eqref{red-K}, and \eqref{plusminus}.
Then $(\uu_\e , \ee_\e, \pp_\e, \zz_\e)$ is the  solution of the $\e$-regularized evolution problem with boundary datum $\ww$ and initial condition $(u_0, e_0, p_0,z_0)$ if and only if
\begin{eqnarray}
& \uu_\e(t,x_1, \ol x)=\sqrt2\,\uu^R_\e(t,x_1) e_2\,,\quad 
\ee_\e(t,x_1, \ol x)= M(\ee^R_{\e}(t)) \,,\label{uuee}\\
&
\pp_\e(t,x_1, \ol x)= M(\pp^R_{\e}(t, x_1)) \,, \quad
\zz_\e(t,x_1, \ol x)=\zz^R_\e (t, x_1)\,, \label{ppzz}
\end{eqnarray}
where $(\uu_\e^R , \ee_\e^R, \pp_\e^R, \zz_\e^R)$ is the  solution of the
$\e$-regularized evolution problem in dimension $d=1$,
corresponding to $\C\xi=2\mu\xi$, with boundary datum $\ww^R$ 
and initial condition $(u_0^R, e_0^R, p_0^R,z_0^R)$. 
\end{theorem}
\begin{proof}
Let us first prove that for every $(\alpha , \zeta)\in \R{\times}\R$ we have
\begin{equation}\label{equiv-P}
P_K(M(\alpha), \zeta) = (M(\hat \alpha), \hat \zeta) \quad \iff \quad
P_{K^{\!R}}( \alpha, \zeta)=(\hat \alpha, \hat \zeta)\,.
\end{equation}
Indeed on the one hand, if $(\hat \alpha, \hat \zeta):=P_{K^{\!R}}( \alpha, \zeta)$, then the distance between $(M( \alpha), \zeta)$ and $(M(\tilde \alpha),\tilde \zeta)$ is larger than the distance between $(M(\alpha), \zeta)$ and $(M(\hat\alpha),\hat \zeta)$. This is an easy consequence of the fact that $M$ is a linear isometry.
On the other hand, if $(\check \sigma, \check \zeta):= P_K(M(\alpha) , \zeta)$,
it follows from \eqref{red-K} that $\check \sigma$ is the projection of $M(\alpha)$ onto the ball $ \{  \sigma\in \MD : |\sigma| \leq |\check \sigma| \}$, which implies that $\check\sigma=M(\check \alpha)$ for some $\check \alpha\in \R$. These two facts together allow to establish \eqref{equiv-P}.

Using the linearity of $M$ we can prove that \eqref{equiv-P} implies 
\begin{equation}\label{equiv-N}
(M(\hat \alpha), \hat \zeta) = N_K^\e (M(\alpha), \zeta) \quad \iff \quad
(\hat\alpha, \hat\zeta) = N_{K^{\!R}}^\e (\alpha, \zeta)
\end{equation}
for every $\alpha$, $ \zeta$, $\hat\alpha$, $\hat\zeta$ in $\R$.

By \eqref{ic} the kinematic admissibility condition (ev1)$\!_0$ 
(Theorem~\ref{young-main2}) for $(u_0, e_0, p_0,z_0)$ is
equivalent to (ev1)$\!_0$ in dimension $d=1$ for $(u_0^R, e_0^R, p_0^R,z_0^R)$, while the equilibrium condition (ev2)$\!_0$ is always satisfied, since $e_0$ and $e_0^R$ are constant.
Assuming (ev1)$\!_0$, let $(\uu_\e^R,\ee_\e^R,\pp_\e^R,\zz_\e^R)$ be the solution of the one dimensional
$\e$-regularized evolution problem with boundary datum $\ww^R$ and initial condition
$(u_0^R, e_0^R, p_0^R,z_0^R)$, and let $(\uu_\e , \ee_\e, \pp_\e, \zz_\e)$ be defined by \eqref{uuee} and \eqref{ppzz}. It is easy to check that $(\uu_\e , \ee_\e, \pp_\e, \zz_\e)$ satisfies the initial conditions (ev0)$\!_\e$ of Definition \ref{def:reym}. Moreover, the kinematic admissibility (ev1)$\!_\e$ follows from the fact that 
$E(\uu_\e(t))=M(D\uu^R_\e (t))$, and from the kinematic admissibility in dimension one
($D$ denotes the distributional derivative with respect to the space variable). 
Equilibrium condition  (ev2)$\!_\e$ follows from the fact that 
$\sigmaa^R_\e (t)$ is constant on $(-\frac12,\frac12)$ by Remark~\ref{rem1d}. Finally the regularized flow rule 
(ev\^3)$\!_\e$ follows from the condition in dimension one and from \eqref{equiv-N}.
\end{proof}

\begin{remark}\label{rem:symmetry}
The one-dimensional problem has also  different interpretations, for example in the study of uniaxial loading of cylindrical bodies. In this case, the symmetry condition \eqref{plusminus} is no longer a natural assumption since one may envisage a different behaviour in tension and compression. 
Theorem~\ref{young-main2} holds true also in this general case, since the only hypotheses on $K$ are \eqref{rk} and \eqref{0zeta}.
\end{remark}

In order to prove  the equivalence between the approximable quasistatic evolutions for simple shears and the solutions to the corresponding reduced one-dimesional problems,
for every $m$ we consider the linear operator 
$\psi_m\colon (\R{\times}\R)^m\to (\MD{\times}\R)^m$ defined by
$$
\psi_m((\beta_1,\theta_1),\dots,(\beta_m,\theta_m))=
((M(\beta_1),\theta_1),\dots,(M(\beta_m),\theta_m))\,.
$$ 
In the following theorem we shall prove that the system of Young measures $\muu$
in the  approximable quasistatic evolution is related to the system Young measure 
$\muu^R$ of the reduced problem by the formula
\begin{equation}\label{ppzzy}
\langle f, 
\muu^{}_{t_1\dots t_m}\rangle
=\!\!\!\!\!\!\!\!\!\! \!\!\!\!\!\int\limits_{\quad\quad(-\frac12,\frac12)^{d-1}}
\!\!\!\!\!\!\!\!\!\!\!\!\!\!\!\!
\langle f((x_1,\ol x),\psi_m(\xi),\eta), 
\muu^R_{t_1\dots t_m}(x_1,\xi,\eta)\rangle
\, d\ol x
\end{equation}
for every $f\in C^{\hom}(\ol Q{\times}(\R{\times}\MD)^m{\times}\R)$ and
every finite sequence $t_1,\dots,t_m$ in $[0,+\infty)$ with $t_1<\dots<t_m$.

\begin{theorem}\label{th:equiv2}
Assume \eqref{bc}, \eqref{ic}, \eqref{red-K}, and \eqref{plusminus}.
Then $(\uu , \ee, \muu)$ is  an approximable quasistatic evolution  with boundary datum
$\ww$ and initial conditions $(u_0, e_0, p_0,z_0)$ if and only if \eqref{uuee} and 
\eqref{ppzzy} hold and $(\uu^R , \ee^R, \muu^R)$ is an approximable quasistatic evolution
for the problem in dimension $d=1$, corresponding to $\C\xi=2\mu\xi$, 
with boundary datum
$\ww^R$ and initial condition $(u_0^R, e_0^R, p_0^R,z_0^R)$. 
\end{theorem}

\begin{proof} By Theorem~\ref{thm:ss}, the solution $(\uu_\e , \ee_\e, \pp_\e, \zz_\e)$ of the $\e$-regularized evolution problem with boundary datum $\ww$ and initial condition $(u_0, e_0, p_0,z_0)$ satisfies~ \eqref{ppzz}. Therefore, if $(\uu_\e^R,\ee_\e^R,\pp_\e^R,\zz_\e^R)$ is the solution of the reduced $\e$-regularized evolution problem
with boundary datum $\ww^R$ and initial condition $(u_0^R, e_0^R, p_0^R,z_0^R)$, then \eqref{ppzzy} holds with $\muu$ replaced by 
$\deltaa_{(\pp_\e,\zz_\e)}$ and $\muu^R$ replaced by $\deltaa_{(\pp^R_\e,\zz^R_\e)}$.
The conclusion follows.
\end{proof}

\end{section}

\begin{section}{A spatially homogeneous example}\label{examples}

In this section we
assume that $d\ge 2$ and  $\C$ is isotropic, which implies that
$$
\C\xi=2\mu\xi_D+\kappa(\tr\,\xi)I
$$ 
for some constants $\mu>0$ and $\kappa>0$.
We also assume that 
\begin{eqnarray}
&K:=\{(\xi,\theta)\in\MD{\times}\R: |\xi|^2+\theta^2\le 1\}\,,\label{K-sphere}
\\
&V(\theta):=\frac12-\frac12\sqrt{1+\theta^2}\,, \qquad \Ga_0:=\partial\Om\,, 
\qquad\Ga_1:=\emptyset\,.\nonumber
\end{eqnarray}
Let us fix a constant $\theta_0>0$ and a $d{\times}d$ matrix $\xi_0$ with $\tr\, \xi_0 = 0$. 
We assume that the symmetric part $\xi_0^s$ of $\xi_0$ is different from $0$. We will examine the approximable quasistatic evolution corresponding to the boundary datum
$$
\ww(t,x):=t\xi_0x\,,
$$
and to the initial conditions 
$$
u_0(x)=0\,,\quad e_0(x)=0\,, \quad p_0(x)=0\,, \quad
z_0(x)=\theta_0\,.
$$

\subsection{The special form of the solution} 
The following theorem shows that, in this case,  
any approximable quasistatic evolution is spatially regular and satisfies 
$\uu(t)=\ww(t)$ for every $t\in [0,+\infty)$,
while $\ee(t)$, $\pp(t)$, and $\zz(t)$ do not depend on $x$. Moreover, $\ee(t)$  
and $\pp(t)$ are proportional to the symmetric part $\xi_0^s$ of $\xi_0$. 
\begin{theorem}\label{thm:ex1}
Let $\C$, $K$, $V$, $\Gamma_0$, $\Gamma_1$, $\theta_0$, $\xi_0$, $\xi_0^s$, $\ww$, $u_0$, $e_0$, $p_0$, and $z_0$ satisfy the assumptions considered at the beginning of this section. Then any approximable quasistatic evolution
 with boundary datum $\ww$ and initial condition  $(u_0, e_0, p_0, z_0)$ is spatially regular and has the special form 
 \begin{equation}\label{expl}
\uu(t,x) :=  t\xi_0x\,,\quad \ee(t,x) := (t-\psi(t))\xi_0^s\,,\quad \pp(t,x) := \psi(t)\xi_0^s\,, \quad
\zz(t,x)=\theta(t)\,,
\end{equation}
where $\psi\colon[0,+\infty)\to[0,+\infty)$ and 
$\theta\colon[0,+\infty)\to[0,+\infty)$ are functions
depending on $\mu$ and $\theta_0$. We have $\psi(t)=0$ and
$\theta(t)=\theta_0$ for $0\le t\le t_0$, where 
\begin{equation}\nonumber
t_0:=\frac1{2\mu|\xi_0^s|}\sqrt{1-V'(\theta_0)^2}\,,
\end{equation}
while
\begin{equation}\label{psi}
\psi(t)=t-\frac1{2\mu|\xi_0^s|}\sqrt{1-V'(\theta(t))^2}
\end{equation}
for $t>t_0$.
\end{theorem}

\begin{proof}
By Proposition~\ref{prop:sp-hom} the  solution $(\uu_\e,\ee_\e,\pp_\e,\zz_\e)$
of the $\e$-regularized evolution problem with boundary datum $\ww$ and initial condition $(u_0, e_0, p_0, z_0)$ has the form 
$$
\uu_\e(t,x)=\ww(t,x)\,,\quad\ee_\e(t,x) := (t-\psi_\e(t))\xi_0^s\,,\quad \pp_\e(t,x) := \psi_\e(t)\xi_0^s\,, \quad
\zz_\e(t,x):=\theta_\e(t)\,,
$$
where $\psi_\e$, $\theta_\e\in H^1_{loc}([0,+\infty))$ is the unique solution of  the Cauchy problem 
\begin{equation}\label{eq-flowrule1000}
(\dot\psi_\e(t)\xi_0^s,\dot\theta_\e(t))\in  N_K^\e(2\mu(t-\psi_\e(t))\xi_0^s,-V'(\theta_\e(t)))
\end{equation}
\begin{equation}\label{initialconditions}
\psi_\e(0)=0\,,\qquad \theta_\e(0)=\theta_0\,.
\end{equation}
According to \eqref{partialHe*} and \eqref{K-sphere} we have
\begin{equation}\label{eq-flowrule1001}
N_K^\e(\sigma,\zeta)=\frac{1}{\e}\frac{(\sigma,\zeta)}{\sqrt{|\sigma|^2+\zeta^2}}(\sqrt{|\sigma|^2+\zeta^2}-1)^+\,.
\end{equation}
As $N_K^\e$ is Lipschitz continuous we have $\psi_\e$, $\theta_\e \in C^1([0,+\infty))$. 
It is easy to see that in the interval $[0,t_0]$ the pair 
$(\psi_\e(t),\theta_\e(t)):=(0,\theta_0)$ is a solution to \eqref{eq-flowrule1000} and \eqref{initialconditions}. 

By Remark~\ref{rmk:eq-def} we have 
\begin{equation}\label{eq-flowrule}
(2\mu(t-\psi_\e(t))\xi_0^s-\e\dot\psi_\e(t)\xi_0^s,-V'(\theta_\e(t))-\e\dot\theta_\e(t) )\in \partial H(\dot\psi_\e(t) \xi_0^s,\dot\theta_\e(t))
\end{equation}
for every $t\in[0,+\infty)$.
Let us prove that 
\begin{equation}\label{G2000}
t-\psi_\e(t)>0
\end{equation}
for every $t\geq t_0$. We argue by contradiction. If \eqref{G2000} is not true, let 
$$
t_1:=\inf\{t>t_0: t-\psi(t)\leq 0\}\,.
$$
Since $\psi_\e(t_0)=0$, by continuity we have $t_1>t_0$, $t_1-\psi_\e(t_1)=0$, and $t- \psi_\e(t)>0$ for $t\in[t_0,t_1)$. It follows that $1-\dot \psi_\e(t_1)\leq0$.
On the other hand,  \eqref{eq-flowrule1000} and \eqref{eq-flowrule1001}  imply that $\dot\psi_\e(t_1)=0$, which contradicts the previous inequality. This concludes the proof of \eqref{G2000}.

We now show that for every $t\ge t_0$ we have
\begin{equation}\label{on-boundary}
(2\mu(t-\psi_\e(t))\xi_0^s-\e\dot\psi_\e(t)\xi_0^s,-V'(\theta_\e(t))-\e\dot\theta_\e(t) )\in \partial K\,,
\end{equation}
which is equivalent to
\begin{equation}\label{on-boundary2}
|\xi_0^s|^2[2\mu(t-\psi_\e(t))-\e\dot\psi_\e(t)]^2+[V'(\theta_\e(t))+\e\dot\theta_\e(t)]^2=1\,.
\end{equation}
Arguing by contradiction, let us consider a maximal open interval $(\tau_1, \tau_2)$ in which \eqref{on-boundary} does not hold. Since $\partial H(\xi,\theta) \subset \partial K$
for every $(\xi,\theta)\neq(0,0)$, we deduce from   \eqref{eq-flowrule} that $\dot\psi_\e(t)=\dot\theta_\e(t)=0$ for every $t\in(\tau_1, \tau_2)$.
Therefore there exist two constants $c_1$ and $c_2$ such that $\psi_\e(t)=c_1$ and
$\theta_\e(t)=c_2$  for every $t\in(\tau_1, \tau_2)$. It follows that
$$
(2\mu(t-\psi_\e(t))\xi_0^s-
\e\dot\psi_\e(t)\xi_0^s,-V'(\theta_\e(t))-\e\dot\theta_\e(t) )=
(2\mu(t-c_1)\xi_0^s,-V'(c_2) )\,,
$$
which by \eqref{eq-flowrule} gives
\begin{equation}\label{contr2000}
\sqrt{4\mu^2(t-c_1)^2|\xi_0^s|^2 + V'(c_2)^2 }\leq1
\end{equation}
By the maximality of $(\tau_1, \tau_2)$ we have $(2\mu(\tau_1-c_1)\xi_0^s,-V'(c_2) )\in \partial K$ and hence
$$
\sqrt{4\mu^2(\tau_1-c_1)^2|\xi_0^s|^2 + V'(c_2)^2 }=1\,.
$$
On the other hand, on the interval $(\tau_1, \tau_2)$ we have $t-c_1>0$ by \eqref{G2000}. This implies that 
$$
\sqrt{4\mu^2(t-c_1)^2|\xi_0^s|^2 + V'(c_2)^2 }>1
$$
for every $t\in (\tau_1, \tau_2)$ which contradicts \eqref{contr2000} and concludes the proof of \eqref{on-boundary}.

Using Definition~\ref{def:qsym}  and the special form of the function $(\uu_\e, \ee_\e,\pp_\e, \zz_\e)$, it is easy to prove that there exist  two left-continuous functions $\psi\colon[0,+\infty)\to\R$ and $\theta\colon[0,+\infty)\to\R$
such that
$$
\psi_{\e_k}(t)\to\psi(t)\,, \qquad \theta_{\e_k}(t)\to\theta(t)
$$ 
for a.e.\ $t\in[0,+\infty)$.
This property, together with the weak$^*$ left-continuity of $\ee$, $\pp$, and $\zz$, implies that
\eqref{expl} holds for every $t\geq 0$.
Using \eqref{primest1} we obtain also
$\e_k\dot\psi_{\e_k}(t)\to 0$ and $\e_k\dot\theta_{\e_k}(t)\to 0$
for a.e.\ $t\geq 0$.

Passing to the limit in \eqref{on-boundary2} and using again left-continuity,
we get
$$
4\mu^2|\xi_0^s|^2 (t-\psi(t))^2+ V'(\theta(t))^2=1
$$
for every $t\geq 0$. Since $t-\psi(t)\geq 0$ by \eqref{G2000}, we obtain~\eqref{psi}.
\end{proof}

\subsection{A second order equation for the regularized evolution}
Thanks to \eqref{expl} and \eqref{psi} the behaviour of an approximable quasistatic evolution is completely determined by a function $\theta$, which is the limit of a sequence 
of funtions $\theta_{\e}$ related to the $\e$-regularized evolution problem.
The following lemma provides a careful analysis of a second order equation satisfied by 
the functions $\theta_\e$.
The special properties of the different terms appearing in this equation will be important in the proof of the asymptotic behaviour of $\theta_\e$ as $\e\to 0$, and consequently in the description of the limit function $\theta$. 
\begin{lemma}\label{lm:td}
Let $\e>0$ and let $\theta_\e$ be the function introduced in the proof of Theorem~\ref{thm:ex1}. Then $\theta_\e$ is the solution of the Cauchy problem
\begin{equation}\label{neweq}
\e\ddot \theta_\e=F_\e(\theta_\e,\dot\theta_\e)
\quad \text{on } [t_0,+\infty)\,,
\qquad
\theta_\e(t_0)=0\,, \quad \dot\theta_\e(t_0)=0\,,
\end{equation}
where 
\begin{eqnarray}
F_\e(\theta,\dot\theta) \!\!\!\! & :=  & \!\!\!\! 
A(\theta,\e\dot\theta)+B(\theta,\e\dot\theta)\dot\theta\,,
\vphantom{\frac{2\mu|\xi_0^s|}{V'(\theta)}}
\label{newf}
\\
A(\theta,v) \!\!\!\! & :=  & \!\!\!\! A_0(\theta)-A_1(\theta,v)v+A_2(\theta,v)v^2\,,
\vphantom{\frac{2\mu|\xi_0^s|}{V'(\theta)}}
\label{asmall}
\\
A_0(\theta) \!\!\!\! & :=  & \!\!\!\! -2\mu|\xi_0^s|V'(\theta)\sqrt{1-V'(\theta)^2}\,,
\vphantom{\frac{2\mu|\xi_0^s|}{V'(\theta)}}
\label{BG}
\\
A_1(\theta,v) \!\!\!\! & :=  & \!\!\!\! 4\mu|\xi_0^s| \sqrt{1-(V'(\theta)+v)^2}+{}
\nonumber
\\
& &
{}+ 2\mu|\xi_0^s|V'(\theta) \frac{\sqrt{1-(V'(\theta)+v)^2} - \sqrt{1-V'(\theta)^2}}{v} \,,
\label{A1}
\\
A_2(\theta,v) \!\!\!\! & :=  & \!\!\!\! -\frac{2\mu|\xi_0^s|}{V'(\theta)}\sqrt{1-(V'(\theta)+v)^2}\,,
\label{A2}
\\
B(\theta,v) \!\!\!\! & :=  & \!\!\!\! -B_0(\theta)+B_1(\theta)v +B_2(\theta)v^2 -B_3(\theta)v^3\,,
\vphantom{\frac{2\mu|\xi_0^s|}{V'(\theta)}}
\label{bsmall}
\\
B_0(\theta) \!\!\!\! & :=  & \!\!\!\! 2\mu(1-V'(\theta)^2)+V'(\theta)^2V''(\theta)\,,
\vphantom{\frac{2\mu|\xi_0^s|}{V'(\theta)}}
\label{AG}
\\
B_1(\theta) \!\!\!\! & :=  & \!\!\!\! -\frac{(2\mu- V''(\theta))(1-3V'(\theta)^2)}{V'(\theta)}\,,
\label{B1}
\\
B_2(\theta) \!\!\!\! & :=  & \!\!\!\! 3(2\mu- V''(\theta))\,,
\vphantom{\frac{2\mu|\xi_0^s|}{V'(\theta)}}
\label{B2}
\\
B_3(\theta) \!\!\!\! & :=  & \!\!\!\! -\frac{2\mu- V''(\theta)}{V'(\theta)}\,.
\label{B3}
\end{eqnarray}
Moreover, $\dot\theta_\e(t)>0$ for every $t>t_0$ and
$\theta_\e(t)\to +\infty$ as $t\to+\infty$.
\end{lemma}

\begin{remark}\label{5091}
We observe that for every $\theta\geq\theta_0$ and every $v\in[0,-V'(\theta)]$ all functions $A_0,A_1,A_2,B_1,B_2,B_3$ are positive, except $B_0$.
In particular, in this range the following inequalities hold:
\begin{eqnarray}
& \displaystyle
 -\sqrt3 \mu|\xi_0^s| V'(\theta_0)\le A_0(\theta)\le \mu|\xi_0^s|\,,
\label{A0est}
\\
& \displaystyle \vphantom{\frac{2\mu|\xi_0^s|}{V'(\theta_0)}}
2\mu|\xi_0^s|\le A_1(\theta,v) \le 4\mu|\xi_0^s|\,,
\label{A1est}
\\
& \displaystyle
2\sqrt3 \mu|\xi_0^s|\le A_2(\theta,v) \le -\frac{2\mu|\xi_0^s|}{V'(\theta_0)}\,,
\label{A2est}
\\
& \displaystyle  
\frac32 \mu-\frac18 \le B_0(\theta)\le 2\mu\,,
\label{BBest}
\\
& \displaystyle  
\mu \le B_1(\theta)\le -\frac{1}{V'(\theta_0)}(2\mu+\frac12)\,,
\label{B1est}
\\
& \displaystyle
6\mu\le B_2(\theta) \le 6\mu+\frac32\,,
\label{B2est}
\\
& \displaystyle  
4\mu\le B_3(\theta) \le -\frac{1}{V'(\theta_0)}(2\mu+\frac12)\,.
\label{B3est}
\end{eqnarray}
\end{remark}

\begin{proof}[Proof of Lemma~\ref{lm:td}]
Let $\psi_\e$ be defined as in the proof of Theorem~\ref{thm:ex1}.
We observe that \eqref{eq-flowrule1000}, \eqref{eq-flowrule1001}, and \eqref{G2000} imply that $\dot\psi_\e(t) \geq 0$
for every $t\geq t_0$. 
Let us prove that also 
\begin{equation}\label{contr2002}
\dot\theta_\e(t) \geq 0 \qquad\hbox{for every }t\ge t_0\,.
\end{equation}
If not, we can define $t_2=\inf E$, where $E:=\{t>t_0: \dot\theta_\e(t)< 0\}$.
By \eqref{eq-flowrule1000} and \eqref{eq-flowrule1001} we have $-V^\prime(\theta_\e(t))<0$ for every $t\in E$. By continuity, 
\begin{equation}\label{contr2001}
-V^\prime(\theta_\e(t_2))\leq0\,.
\end{equation}
On the other hand $\dot\theta_\e(t) \geq 0$ for every $t \in [t_0, t_2]$. Hence $0< \theta_0=\theta_\e(t_0)\leq  \theta_\e (t_2)$. As $-V^\prime$ is increasing, we dedude that $0<-V^\prime(\theta_0)\leq V^\prime(\theta_\e (t_2))$, which contradicts \eqref{contr2001} and concludes the proof of \eqref{contr2002}.

By \eqref{on-boundary}, the flow rule \eqref{eq-flowrule} is equivalent to saying that for every $t\ge t_0$ there exists $\lambda(t)\ge 0$ such that 
\begin{eqnarray}
\lambda(t)(2\mu(t-\psi_\e(t))-\e\dot\psi_\e(t))=
\dot\psi_\e(t)\,,\label{system-1}
\\
\lambda(t)(-V'(\theta_\e(t))-\e\dot\theta_\e(t))=
\dot\theta_\e(t)\,.\label{system-2}
\end{eqnarray}
To simplify the notation it is convenient to introduce the function  $W_\e\colon\R{\times}\R\to\R$ defined by 
$$
W_\e(\theta,y):=V'(\theta)+\e y. 
$$
We shall frequently use the inequalities
$$
-\frac12 < W_\e(\theta,y)<0
$$
for every $\theta>0$ and $y\in[0,-V'(\theta)/\e)$.
Let us prove that 
\begin{equation} \label{adesso2000}
2\mu(t-\psi_\e(t))-\e\dot\psi_\e(t)>0 \,, \qquad W_\e(\theta_\e(t),\dot\theta_\e(t))<0
\end{equation}
for every $t\ge t_0$.
If $\dot\psi_\e(t)>0$, the first inequality in \eqref{adesso2000} follows from \eqref{system-1}.
If $\dot\psi_\e(t)=0$, the same inequality follows from \eqref{G2000}. If $\dot\theta_\e(t)>0$, the second inequality in \eqref{adesso2000} follows from \eqref{system-2}.
Finally, as $\theta_\e$ is nondecreasing by \eqref{contr2002}, we have $0<-V^\prime(\theta_0)\leq -V^\prime(\theta_\e (t))$ for every $t\geq t_0$. This implies the second inequality in \eqref{adesso2000} when $\dot\theta_\e(t)=0$.

From \eqref{system-2} and \eqref{adesso2000} we obtain
$$
\lambda(t)=
-\frac{\dot\theta_\e(t)}{W_\e(\theta_\e(t),\dot\theta_\e(t))}\,.
$$
Substituting $\lambda(t)$ in \eqref{system-1}, we have 
\begin{equation}\label{psi1000}
-\frac{2\mu\dot\theta_\e(t)(t-\psi_\e(t))}{W_\e(\theta_\e(t),\dot\theta_\e(t))}
=
\dot\psi_\e(t)\Big(1-
\frac{\e\dot\theta_\e(t)}{W_\e(\theta_\e(t),\dot\theta_\e(t))}\Big)
=\dot\psi_\e(t) \frac{V'(\theta_\e(t))}{W_\e(\theta_\e(t),\dot\theta_\e(t))}\,.
\end{equation}

By \eqref{on-boundary2}, \eqref{adesso2000}, and \eqref{psi1000} we obtain
\begin{equation}\label{star}
\dot\psi_\e(t)=-\frac{1}{|\xi_0^s|}
\frac{\sqrt{1-W_\e^2(\theta_\e(t),\dot\theta_\e(t))}}{
W_\e(\theta_\e(t),\dot\theta_\e(t))}\dot\theta_\e(t)\,.
\end{equation}
{}From \eqref{psi1000} it follows that
$$
2\mu(t-\psi_\e(t))\dot \theta_\e(t)=
\frac{V'(\theta_\e(t))}{|\xi_0^s|}\frac{\sqrt{1-W_\e^2(\theta_\e(t),
\dot\theta_\e(t))}}{W_\e(\theta_\e(t),\dot\theta_\e(t))}
\dot \theta_\e(t)\,.
$$
On the set where $\dot\theta_\e(t) \neq 0$ we obtain
\begin{equation}\label{psi1002}
2\mu(t-\psi_\e(t))=\frac{V'(\theta_\e(t))}{|\xi_0^s|}
\frac{\sqrt{1-W_\e^2(\theta_\e(t),\dot\theta_\e(t))}}{
W_\e(\theta_\e(t),\dot\theta_\e(t))}\,.
\end{equation}
On the set where $\dot\theta_\e = 0$ \eqref{psi1002} reduces to 
$$
2\mu(t-\psi_\e(t))=
\frac{\sqrt{1-W_\e^2(\theta_\e(t),\dot\theta_\e(t))}}{|\xi_0^s|}\,.
$$
which is a consequence of \eqref{on-boundary2} and \eqref{adesso2000}, taking into account that 
$\dot\psi_\e(t)=0$ where $\dot\theta_\e(t) = 0$
(this can be easily deduced from \eqref{system-1}, \eqref{system-2}, 
and~\eqref{adesso2000}).

By differentiating \eqref{psi1002} and using \eqref{star} we finally obtain that $\theta_\e$ satisfies the equation in \eqref{neweq}.

Since $\theta_\e$ is of class $C^1$ on $[0,+\infty)$ and $\theta_\e(t)=0$ for $t\in[0,t_0)$, we deduce that $\theta_\e(t_0)=0$ and $\dot\theta_\e(t_0)=0$.

Let us prove that 
\begin{equation}\label{num**}
\dot\theta_\e(t)>0 \qquad \hbox{for every } t>t_0\,.
\end{equation}
{} First we observe that the equation gives $\ddot\theta_\e(t_0)>0$, so $\dot\theta_\e(t)>0$ in a right neighbourhood of $t_0$. If this inequality is not satisfied for every $t>t_0$, let $\tau$ be the first time greater than $t_0$ such that $\dot\theta_\e(\tau)=0$. From the equation we get $\ddot\theta_\e(\tau)=F_\e(\theta_\e(\tau),0)=A_0(\theta_\e(\tau))>0$, which leads to a contradiction with the definition of $\tau$, and concludes the proof 
of~\eqref{num**}.

Let 
$$
\bar\theta:=\lim_{t\to+\infty}\theta_\e(t)\,.
$$
We want to prove that $\bar\theta=+\infty$. If not,
by \eqref{adesso2000} and by the monotonicity of $\theta_\e$, we have that 
$$
0\le\e\dot\theta_\e(t)\le -V'(\theta_\e(t))\le -V'(\bar\theta)\,, \qquad
\theta_0\le\theta_\e(t)\le\bar\theta\,.
$$
This implies that the pair $(\theta_\e(t),\dot\theta_\e(t))$ is bounded on $[t_0,+\infty)$. It follows from  equation \eqref{neweq} that $\dot\theta_\e$ is globally Lipschitz continuous on $[t_0,+\infty)$ with some constant $L>0$. 

Let us prove that 
\begin{equation}\label{limtdot}
\lim_{t\to+\infty}\dot\theta_\e(t)=0\,.
\end{equation}
If not there exist a sequence $t_k\to+\infty$ and a constant $c>0$ such that $\dot\theta_\e(t_k)\geq c$ for every $k$. It follows from the Lipschitz continuity  of $\dot\theta_\e$ that 
$$
\dot\theta_\e(t)\geq \frac{c}{2}
$$
on the intervals $(t_k-\frac{c}{2L}, t_k+\frac{c}{2L})$. This obviously contradicts the assumption that $\bar\theta<+\infty$ and concludes the proof of \eqref{limtdot}. 

Let us consider a sequence $\tau_k\to+\infty$ such that $\ddot\theta_\e(\tau_k)\to 0$. 
It follows from \eqref{neweq} that 
$F_\e(\theta_\e(\tau_k),\dot\theta_\e(\tau_k) )\to 0$. 
Passing to the limit we deduce $F_\e(\bar\theta,0)=0$, which is impossible since 
$F_\e(\bar\theta,0)=A_0(\bar\theta)>0$.
This shows that $\bar\theta=+\infty$ and concludes the proof of the lemma.
\end{proof}

\subsection{Fast  dynamics}
In the rest of the section we shall prove that the function $\theta(t)$ given in Theorem~\ref{thm:ex1}, which determines the internal variable in the quasistatic evolution, has at most a jump time $\tau\in[t_0,+\infty)$. We shall show that the sequence $\theta_\e(t)$ converges to $\theta(t)$ uniformly on the compact sets of $[t_0,\tau)\cup(\tau,+\infty)$ and that on each of the two intervals $[t_0,\tau)$ and $(\tau,+\infty)$ $\theta(t)$ is a solution of the equation of the slow dynamics 
\begin{equation}\nonumber
B_0(\theta(t))\dot\theta(t)=A_0(\theta(t))\,.
\end{equation}
Moreover $B_0(\theta(t))>0$ for every $t\in (t_0,\tau)\cup(\tau,+\infty)$.
More precisely we shall prove that if $B_0(\theta)\geq0$ for every $\theta\in[\theta_0,+\infty)$, then there is no jump point and therefore $\theta(t)$ is the unique solution to the Cauchy problem
$$
B_0(\theta(t))\,\dot\theta(t)=A_0(\theta(t))\,, \qquad \theta(t_0)=\theta_0\,.
$$

We shall prove (see Lemma~\ref{signB}) that for some values of $\mu$ we have $B_0(\theta)\geq0$ for every $\theta\geq 0$. For other values of $\mu$ the function $B_0$ has exactly two zeros $\alpha$ and $\beta$ and $B_0(\theta)>0$ for $\theta\in(0,\alpha)$, $B_0(\theta)<0$ for $\theta\in (\alpha,\beta)$, and $B_0(\theta)>0$ for $\theta\in(\beta,+\infty)$. In this case there is a jump time $\tau$ only when $\theta_0<\beta$. More precisely if 
$\theta_0\in [\alpha,\beta)$, then $\tau=t_0$, while $\tau>t_0$ when $\theta_0\in (0,\alpha)$. In both cases the jump of $\theta$ can be determined from the relation
$$
\theta(\tau+)=\Phi(\theta(\tau-))\,,
$$
where the transition function $\Phi$ is defined, through the equation of the fast dynamics
\begin{equation}\label{fastdyn}
w'(\theta)=B(\theta, w(\theta))\,.
\end{equation}
More precisely for every $\gamma\in [\alpha,\beta)$ $\Phi(\gamma)$ is defined as the first zero greater than $\gamma$ of the solution to the Cauchy problem
\begin{equation}\label{cauchy-fast}
w'_\gamma(\theta)=B(\theta, w_\gamma(\theta))\,,\qquad w_\gamma(\gamma)=0\,.
\end{equation}
Note that $\theta(\tau-)=\theta_0$ if $\theta_0\in[\alpha,\beta)$ and hence $\tau=t_0$, while $\theta(\tau-)=\alpha$ if $\theta_0\in(0,\alpha)$.

The following lemma is concerned with the changes in sign of the function $B_0$.
\begin{lemma}\label{signB}
Let
$$
\mu_0:=\frac{79\sqrt{19}-344}{108}\sqrt{7+2\sqrt{19}}= 0.0129\dots 
\quad\hbox{and}\quad \alpha_0:=\frac{\sqrt2}{3}\sqrt{\sqrt{19}-1}=0.8639\dots \,.
$$
The following properties hold:
\begin{itemize}
\item[(a)] if $\mu> \mu_0$, then $B_0(\theta)> 0$ for every $\theta>0$;
\item[(b)] if $\mu= \mu_0$, then $B_0(\theta)> 0$
for every $\theta\neq\alpha_0$;
\item[(c)] if $\mu< \mu_0$, then there exist $\alpha$ and $\beta$, with
 $0<\alpha<\alpha_0<\beta<+\infty$, such that $B_0(\theta)> 0$
for every $\theta\in (0,\alpha)\cup(\beta,+\infty)$ and $B_0(\theta)< 0$ for every
$\theta\in (\alpha,\beta)$; moreover we have $B'_0(\alpha)<0$ and $B_0'(\beta)>0$.
\end{itemize}
\end{lemma}

\begin{proof}
The results  follow from the fact that the function
$$
\frac{V'(\theta)^2V''(\theta)}{V'(\theta)^2-1}=
\frac12\frac{\theta^2}{(4+3\theta^2)(1+\theta^2)^{3/2}}
$$
has a positive derivative in $(0,\alpha_0)$, negative derivative in $(\alpha_0,+\infty)$, tends to $0$ for $\theta\to 0$ and $\theta\to+\infty$, and takes the value $2\mu_0$ at
$\theta=\alpha_0$.
\end{proof}

The following lemma provides some properties of the solutions of the equation of the fast dynamics \eqref{fastdyn}.

\begin{lemma}\label{lm:fastdyn}
Assume $\mu<\mu_0$. For every $\gamma\in[\alpha,\beta)$ there exists a point $\Phi(\gamma)>\beta$ such that the solution $w_\gamma$ to the Cauchy problem \eqref{cauchy-fast} satisfies $w_\gamma(\theta)>0$ for every $\theta\in (\gamma,\Phi(\gamma))$, $w_\gamma(\Phi(\gamma))=0$, and $w'_\gamma(\Phi(\gamma))<0$. Moreover the function $\Phi\colon [\alpha,\beta)\to (\beta,\Phi(\alpha)]$ is continuous, decreasing, and $\Phi(\gamma)\to\beta$ as $\gamma\to \beta$.
\end{lemma}

\begin{proof}
By direct computations we find that 
\begin{equation}\label{bneg}
B(\theta, -V'(\theta))=-V''(\theta)>0\,,
\end{equation}
which implies that $w(\theta):=-V'(\theta)$ is a solution of the equation. By uniqueness we have that $w_\gamma(\theta)\leq V'(\theta)$.

For $v=0$ we have that
\begin{equation}\label{bsign}
\begin{array}{ll}
B(\theta,0)=-B_0(\theta)>0 & \text{for }\theta\in(\alpha,\beta)\,, \\
B(\theta,0)=-B_0(\theta)<0 & \text{for }\theta\in(\beta,+\infty)\,.
\end{array}
\end{equation}

Since $B(\theta,v)$ is a polynomial of the third degree in $v$, the sign of its coefficients implies that for every $\theta>0$ there exists a unique point $v_*(\theta)$ such that
$v\mapsto B(\theta,v)$ is increasing on $[0,v_*(\theta)]$ and decreasing on $[v_*(\theta),+\infty)$. It follows from \eqref{bsign} and \eqref{bneg} that 
\begin{equation}\label{btv}
B(\theta,v)>0 \quad \text{for } \theta\in(\alpha,\beta) \text{ and } v\in[0,-V'(\theta)]\,.
\end{equation}

Using the definitions  \eqref{bsmall}--\eqref{B3}, we can prove that 
\begin{eqnarray}
& \displaystyle \label{blimit}
\lim_{\theta\to +\infty} B(\theta,v)= \mu(-\tfrac32+ v+6v^2-4v^3)
\\
& \displaystyle \label{b'limit}
\lim_{\theta\to +\infty} \partial_v B(\theta,v)= \mu(1 +12v-12v^2)
\end{eqnarray}
uniformly for every $v\in[0,\frac12]$. 
It follows that there exists $\bar\theta>\beta$ such that
\begin{equation}\label{pbest}
\partial_v B(\theta,v)\geq \tfrac12 \mu
\end{equation}
for every $\theta\geq\bar\theta$ and every $v\in[0,\frac12]$.

For every $v\in(0,-V'(\alpha))$ let $w(\theta,v)$ be the solution of the Cauchy problem 
$$
w'(\theta,v)=B(\theta, w(\theta,v))\,,
\qquad
w(\alpha,v)=v\,,
$$
where prime denotes derivative with respect to $\theta$.
We want to prove that for every $v\in(0,-V'(\alpha))$ there exists $\hat\theta>\alpha$ such that $w(\hat\theta,v)=0$. 

If not, there exists $v_0\in (0,-V'(\alpha))$ such that $w(\theta,v_0)>0$ for every $\theta>\alpha$. By comparison we have also that $0<w(\theta,v)<-V'(\theta)$ for every $\theta>\alpha$ and every $v\in[v_0,-V'(\alpha))$. Note that this inequality ensures that $\theta\mapsto w(\theta,v)$ is defined on $[\alpha,+\infty)$. Taking the partial derivative with respect to $v$, we obtain
$$
(\partial_v w)'(\theta,v)=\partial_v B(\theta, w(\theta,v))\partial_v w(\theta,v)\,,
\qquad
\partial_v w(\alpha,v)=1\,.
$$
By \eqref{pbest} we have
$$
\partial_v w(\theta,v)\geq \partial_v w(\bar\theta,v)e^{\frac12 \mu(\theta-\bar\theta)}
$$
for $\theta\geq\bar\theta$. Integrating this inequality with respect to $v$ between $v_0$ and $-V'(\alpha)$, we obtain
$$
\begin{array}{c}
-V'(\theta)-w(\theta,v_0)=w(\theta,-V'(\alpha))-w(\theta,v_0)\geq
\\
\geq
(w(\bar\theta,-V'(\alpha))-w(\bar\theta,v_0))e^{\frac12 \mu(\theta-\bar\theta)}
=(-V'(\bar\theta)-w(\bar\theta,v_0))e^{\frac12 \mu(\theta-\bar\theta)}\,.
\end{array}
$$
As $-V'(\bar\theta)-w(\bar\theta,v_0)>0$ by uniqueness and $-V'(\theta)\to1/2$, we conclude that $w(\theta,v_0)\to -\infty$ as $\theta\to +\infty$, which contradicts our assumption.

By \eqref{btv} we have $w_\gamma(\theta)\geq 0$ for $\theta\in(\gamma,\beta]$ 
whenever $\gamma\in (\alpha,\beta)$.  The same conclusion can be obtained also when $\gamma=\alpha$, using the fact that 
$w'_\alpha(\alpha)=B(\alpha, 0)=-B_0(\alpha)=0$ and
$w''_\alpha(\alpha)=\partial_\theta B(\alpha,0)=-B'_0(\alpha)>0$  (this can be easily deduced from equation \eqref{cauchy-fast} and from Lemma~\ref{signB}). On the other hand by comparison we have $w_\gamma(\theta)\leq w(\theta, v)$ for every $v\in (0, -V'(\alpha))$. Since the function $w(\cdot, v)$ vanishes at some point we conclude that $w_\gamma$ vanishes too. Then we define $\Phi(\gamma)$ as the smallest $\theta>\gamma$ for which $w_\gamma(\theta)=0$. From the previous discussion it follows that $\Phi(\gamma)>\beta$, hence $w'_\gamma(\Phi(\gamma))< 0$
 by \eqref{bsign}.  
This fact, together with the continuous dependence on the initial data, implies that $\Phi$ is continuous. It follows from a comparison argument that $\alpha\le\theta_1<\theta_2<\beta$ implies that $\beta<\Phi(\theta_2)<\Phi(\theta_1)\le \alpha$.
Arguing as before, we can prove that the solution $w_\beta$ of the Cauchy problem 
\eqref{cauchy-fast} for $\gamma=\beta$ satisfies $w'_\beta(\beta)=0$ and
$w''_\beta(\beta)<0$, hence $w_\beta(\theta)<0$ for $\theta$ in a right neighbourhood of $\beta$. By the continuous dependence on the data this implies that $\Phi(\gamma)\to\beta$ as $\gamma\to\beta$.
\end{proof}

\subsection{Complete description of the quasistatic evolution}
We are now ready to prove the main result of this section.
Let $\mu_0$, $\alpha$, and $\beta$ be as in Lemma~\ref{signB}.
\begin{theorem}\label{thm:ex}
Under the hypotheses of Theorem~\ref{thm:ex1} there exists a unique approximable quasistatic evolution with boundary datum $\ww$ and initial condition $(u_0, e_0, p_0, z_0)$. It is spatially regular and satisfies  \eqref{expl} and \eqref{psi}.
\begin{itemize}
\item[(a)] If $\mu\ge \mu_0$ or if $\mu<\mu_0$ and $\theta_0\ge\beta$, then in
the interval $(t_0,+\infty)$ the function
$\theta(t)$ coincides with the unique absolutely continuous solution of the Cauchy problem
\begin{equation}\label{cauchy1}
B_0(\theta(t))\,\dot\theta(t)=A_0(\theta(t))\,, \qquad \theta(t_0)=\theta_0\,.
\end{equation}
\item[(b)] If $\mu<\mu_0$ and $\alpha\le \theta_0<\beta$, then in the interval 
$(t_0,+\infty)$ the function
$\theta(t)$ coincides with  the unique absolutely continuous solution of the Cauchy problem
\begin{equation}\label{Bcauchy}
B_0(\theta(t))\,\dot\theta(t)=A_0(\theta(t))\,, \qquad \theta(t_0)=\Phi(\theta_0)\,.
\end{equation}
\item[(c)] if $\mu< \mu_0$ and $0<\theta_0<\alpha$, then there exists a unique $\tau>t_0$
with $\theta(\tau)=\alpha$, such that on the interval $(t_0,\tau]$
$\theta(t)$ coincides with  the unique absolutely continuous solution of the Cauchy problem
\begin{equation}\label{Ccauchy}
B_0(\theta(t))\,\dot\theta(t)=A_0(\theta(t))\,, \qquad \theta(t_0)=\theta_0\,,
\end{equation}
while on the interval 
$(\tau,+\infty)$ the function
$\theta(t)$ coincides with  the unique absolutely continuous solution of the Cauchy problem
\begin{equation}\label{Phicauchy}
B_0(\theta(t))\,\dot\theta(t)=A_0(\theta(t))\,, \qquad \theta(\tau)=\Phi(\alpha)\,.
\end{equation}
\end{itemize}
\end{theorem}

\begin{remark}\label{theta-infty}
In all cases 
$$
\theta(t)\to+\infty\quad\text{and}\quad 
t-\psi(t)\to \frac{\sqrt3}{4}\frac{1}{\mu|\xi_0^s|}
$$ 
as $t\to+\infty$. By \eqref{expl} we have 
$$
\sigmaa(t)\to \frac{\sqrt3}{4}\frac{\xi_0^s}{|\xi_0^s|}
$$ 
as $t\to+\infty$. As $\theta(t)$ is increasing we deduce from \eqref{expl} and \eqref{psi} 
that $|\sigmaa(t)|$ is decreasing for $t\ge t_0$, reflecting the fact that the material softens as  plastic deformation proceeds. The asymptotic value of the stress coincides with 
the value that can be obtained in the perfectly plastic case with elastic domain
$$
K_\infty:=\bigcap_{\theta\in \R}K(-V'(\theta))\,,
$$
where $K(\zeta)$ is defined by \eqref{Kzeta}.
This reflects the fact that, with this loading program, the material behaves in the weakest
way permitted by its internal variable as time tends to infinity.

In case (a) the functions $\ee$, $\pp$, and $\zz$ are continuous with respect to $t$. In case (b) they are discontinuous only at $t=t_0$ while in case (c) they are discontinuous only at a point $\tau>t_0$ (see Figures 1, 2, 3).

\begin{figure}[ht]
\begin{center} 
  \includegraphics[width=0.56\textwidth]{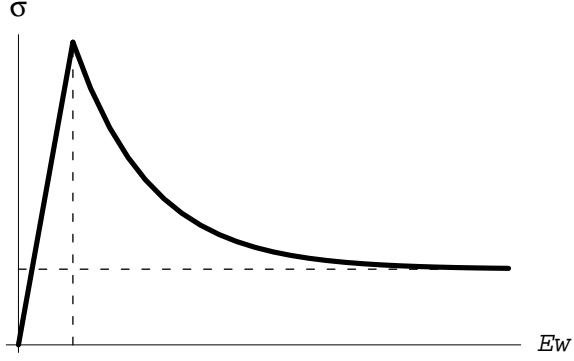}
\caption{Stress vs.\ imposed strain in case (a).}
\end{center}
\end{figure}

\begin{figure}[ht]
\begin{center} 
  \includegraphics[width=0.56\textwidth]{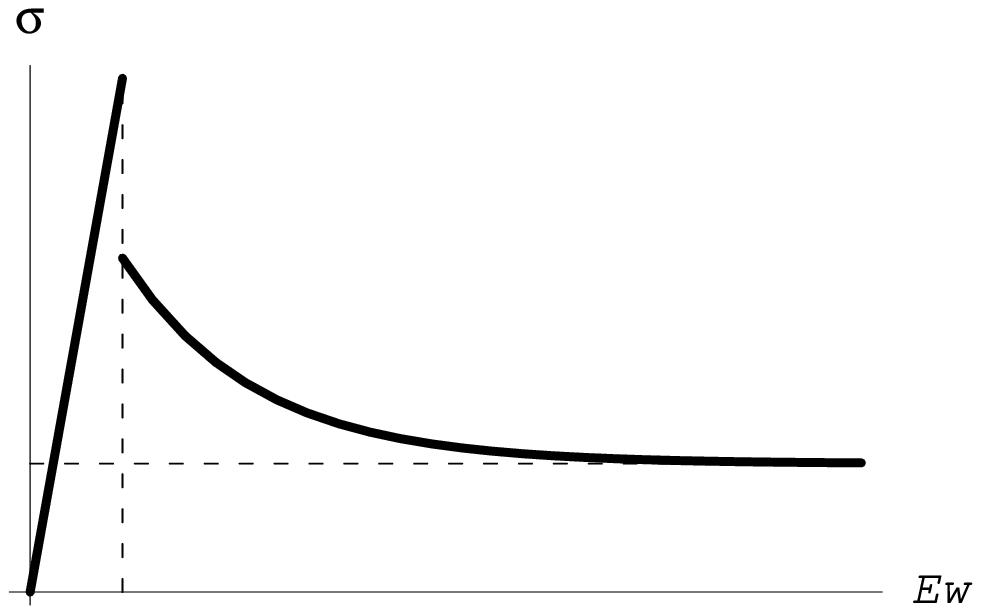}
\caption{Stress vs.\ imposed strain in case (b).}
\end{center}
\end{figure}

\begin{figure}[ht]
\begin{center} 
  \includegraphics[width=0.56\textwidth]{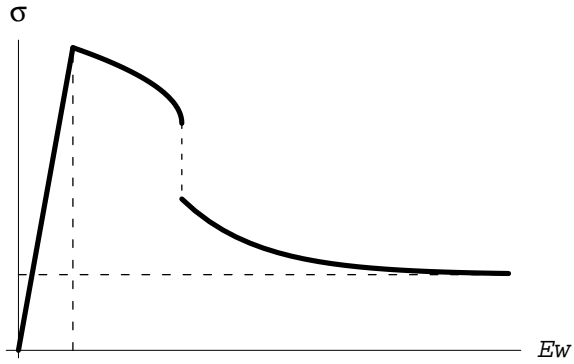}
\caption{Stress vs.\ imposed strain in case (c).}
\end{center}
\end{figure}

In case (a)  the energy inequality \eqref{ineqreduced} of Remark~\ref{RemThm55} holds with equality. In cases (b) and (c) the inequality is strict when $T$ is larger than the discontinuity time. In other words, an instantaneous dissipation occurs at the discontinuity time.
\end{remark}

\begin{proof}[Proof of Theorem~\ref{thm:ex}]
We begin by studying case (c).
Let  $\varphi\colon (0,\alpha)\cup(\beta,+\infty)\to(0,+\infty)$ be the $C^\infty$ function defined by
\begin{equation}\label{defphi}
\varphi(\theta):=\frac{A_0(\theta)}{B_0(\theta)}\,.
\end{equation}
By direct computations we see that
\begin{eqnarray}
& \displaystyle
\lim_{\theta\to+\infty}\varphi(\theta)=\frac{1}{\sqrt3}|\xi_0^s|\,,
\label{limvarphi}
\\
& \displaystyle
\inf_{\theta\geq\theta_0}\varphi(\theta)\geq \varphi_0
:=-\frac{\sqrt3}{2}|\xi_0^s|V'(\theta_0)>0\,,
\label{bdvarphi}
\\
& \displaystyle \vphantom{\frac{\sqrt3}{2}}
\lim_{\theta\to+\infty}\varphi'(\theta)=0\,.
\label{bddphi'}
\end{eqnarray}

For every $\e>0$ it is useful to introduce the 
$\e$-scaled version of the velocity, considered as a function of $\theta$, and denoted by $v_\e(\theta)$. It is characterized by the equality
\begin{equation}\label{newv}
v_\e(\theta_\e(t)):=\e \dot \theta_\e(t)\,,
\end{equation}
and satisfies the condition
\begin{equation}\label{newv'}
v_\e'(\theta_\e(t))=\e \frac{\ddot \theta_\e(t)}{\dot \theta_\e(t)}
=\e^2 \frac{\ddot \theta_\e(t)}{v_\e(\theta_\e(t))}\,,
\end{equation}
where prime denotes the derivative with respect to $\theta$. From \eqref{neweq}--\eqref{B3} it follows that 
\begin{equation}\label{neweqv}
v_\e'(\theta) = \e\frac{A_0(\theta)}{v_\e(\theta)} + 
B_\e(\theta,v_\e(\theta))\,,
\end{equation}
where 
\begin{equation}\label{defBd}
B_\e(\theta,v):=B(\theta,v)-\e A_1(\theta,v)+\e A_2(\theta,v)v\,.
\end{equation}
By \eqref{adesso2000} and \eqref{num**} we have
\begin{equation}\label{<V'}
0\le v_\e(\theta)<-V'(\theta)\qquad\hbox{for every }\theta\ge\theta_0\,.
\end{equation}
The proof will be split
into several steps which describe the behaviour of $\theta_\e(t)$ in different intervals.
Since we are considering case (c), we have $\mu<\mu_0$ and $0<\theta_0<\alpha$. 

\medskip\noindent{\it Step $\rm(c)_1$. Behaviour near $\theta_0$.} 
Let us fix $\theta_1>\theta_0$ such that $B_0(\theta)>0$ for $\theta\in[\theta_0,\theta_1]$. Let $0<\eta_1<\min\{\varphi_0 , \theta_1-\theta_0 , 1\}$, and let $\e_1\in(0,1)$.
For every $\eta\in (0,\eta_1)$
and  every $\e\in (0, \e_1)$ let
$\thone$ be the the largest element of $(\theta_0, \theta_1]$ such
that $v_\e(\theta)\le \e(\varphi(\theta)-\eta)$ for every
$\theta\in [\theta_0,\thone]$.  

In the interval $(\theta_0,\thone]$ it is convenient to write equation \eqref{neweqv} in the form 
\begin{equation}\label{cinqueeqv}
v'_\e(\theta)=\e\frac{A_0(\theta)}{v_\e(\theta)}-B_0(\theta)+
R_\e(\theta,v_\e(\theta))\,,
\end{equation}
where, by \eqref{bsmall} and \eqref{defBd},
\begin{equation}\label{defR}
R_\e(\theta,v):= -\e A_1(\theta,v)+\e A_2(\theta,v)v
+B_1(\theta)v+B_2(\theta)v^2-B_3(\theta)v^3\,.
\end{equation}
By \eqref{A1est}--\eqref{B3est}  and \eqref{<V'} there exists a constant $r_0$, 
independent of  $\e$ and $\eta$, such that the rest $R_\e(\theta,v_\e(\theta))$
can be estimated by
\begin{equation}\label{resto1}
|R_\e(\theta,v_\e(\theta))|\leq r_0 \e
\end{equation}
for all $\theta\in [\theta_0,\thone]$. 
Since 
$$
\frac{\e}{v_\e(\theta)}\ge \frac{1}{\varphi(\theta)-\eta}
$$
for every $\theta\in(\theta_0,\thone]$, recalling \eqref{defphi} we obtain 
$$
v'_\e(\theta)\geq\frac{\eta}{\varphi(\theta)-\eta}B_0(\theta)+
R_\e(\theta, v_\e(\theta))
$$
in the same interval. By \eqref{bdvarphi} and \eqref{resto1} for every
$\theta\in (\theta_0,\thone]$ we have
$$
v'_\e(\theta)\geq \frac{b_0}{\varphi_0}\eta-r_0\e\,,
$$ 
where $b_0>0$ is the minimum of $B_0(\theta)$ on the interval $[\theta_0,\theta_1]$.
This implies 
$$
v_\e(\theta)\geq \frac{b_0\eta}{2\varphi_0}(\theta-\theta_0)\,.
$$
whenever $0<\e<b_0\eta/(2r_0\varphi_0)$.

Let $\varphi_1$ be the maximum of  $\varphi(\theta)$ on the interval $[\theta_0,\theta_1]$.
As $v_\e(\thone)\le \e \varphi(\thone)\le\e \varphi_1$, we obtain
\begin{equation}\label{thone}
\thone-\theta_0 \leq  \frac{\e}{\eta}\frac{2\varphi_0\varphi_1}{b_0}<\theta_1-\theta_0\,,
\end{equation}
whenever $0<\e<\e_1(\eta):=
\min\{\eta b_0/(2r_0\varphi_0),\eta(\theta_1-\theta_0)b_0/(2\varphi_0\varphi_1), \e_1\}$.
This implies that 
\begin{equation}\label{243}
v_\e(\thone)=\e(\varphi(\thone)-\eta)\,.
\end{equation}

Let $\tone$ be the time such that $\theta_\e(\tone)=\thone$.  From \eqref{A0est},
\eqref{BBest}, \eqref{cinqueeqv}, and \eqref{resto1} we obtain
$$
\frac{\e}{v_\e(\theta)}a_0\leq v'_\e(\theta)+2\mu+r_0\e\qquad \hbox{in } (\theta_0,\thone]\,,
$$
where $a_0:=-\sqrt3 |\xi_0^s| V'(\theta_0)>0$. This implies
\begin {eqnarray*}
& \displaystyle\tone-t_0=
\e\int_{\theta_0}^{\thone}\frac{d\theta}{v_\e(\theta)}\leq 
\frac{1}{a_0}\int_{\theta_0}^{\thone}v'_\e(\theta)\, d\theta+
\frac{2\mu}{a_0}(\thone-\theta_0)+\frac{r_0}{a_0}\e(\thone-\theta_0)\leq
\nonumber\\
&\displaystyle \leq \frac{\e}{a_0}\varphi(\thone)+
\frac{\e}{\eta}\frac{4\mu\varphi_0\varphi_1}{a_0b_0}+
\frac{\e^2}{\eta}\frac{2 r_0\varphi_0\varphi_1}{a_0b_0}\,. 
\end{eqnarray*}
Hence there exists a constant $c_1$ such that 
\begin{equation}\label{tone}
\tone-t_0\leq c_1\frac{\e}{\eta}
\end{equation}
for $0<\eta<\eta_1$ and $0<\e< \e_1(\eta)$.

\medskip\noindent{\it Step $\rm(c)_2$. Behaviour between $\theta_0$ and $\alpha$.}
We want to prove that for every $\eta\in(0,\eta_1)$ there exists
$\e_2(\eta)\in(0,\e_1(\eta))$ such that
\begin{equation}\label{newfluss}
\Big|\frac{v_\e(\theta)}{\e} -\varphi(\theta)\Big|\le\eta
\end{equation}
for every $\e\in(0,\e_2(\eta))$ and every $\theta\in[\thone,\alpha-\eta]$.
Since this inequality is satisfied for $\theta=\thone$ by \eqref{243},
considering the equation satisfied by $v_\e/\e$, it is enough to prove that
$$
\frac{A_0(\theta)}{\varphi(\theta)+\eta} -B_0(\theta) 
+R_\e(\theta,\e(\varphi(\theta)+\eta))
< \e\varphi'(\theta) 
< \frac{A_0(\theta)}{\varphi(\theta)-\eta} -B_0(\theta) 
+R_\e(\theta,\e(\varphi(\theta)-\eta))\,.
$$
Arguing as in the previous step, we can show that
$R_\e(\theta,\e(\varphi(\theta)\pm\eta))\to 0$ as $\e\to 0$,
uniformly on $[\theta_0,\alpha-\eta]$. The conclusion follows from the fact that
$\e\varphi'(\theta)\to 0$ uniformly on the same interval and that 
$$
\frac{A_0(\theta)}{\varphi(\theta)+\eta} -B_0(\theta) =-\frac{\eta}{\varphi(\theta)+\eta}B_0(\theta)<0\,,
\qquad
\frac{A_0(\theta)}{\varphi(\theta)-\eta} -B_0(\theta) = \frac{\eta}{\varphi(\theta)-\eta}B_0(\theta)>0
$$
on $[\theta_0, \alpha-\eta]$. This concludes the proof of \eqref{newfluss}.

Let $\ttwo$ be the time such that $\theta_\e(\ttwo)=\alpha-\eta$. By \eqref{newfluss}
 we have
\begin{equation}\label{fluss}
|\dot\theta_\e(t) -\varphi(\theta_\e(t))|\le\eta
\end{equation}
for every $t\in[\tone,\ttwo]$.
{}From \eqref{bddphi'} it follows that for every $t\in[\tone,\ttwo]$
\begin{equation}\label{thetabdd}
0<\varphi_0-\eta\le \dot\theta_\e(t) \le \hat \varphi_0^\eta
\end{equation}
where $\hat \varphi_0^\eta$ is the maximum of $\varphi(\theta)+\eta$ on the interval
$[\theta_0,\alpha-\eta]$. 
This implies that 
$$
\frac{\alpha-\eta-\thone}{\hat \varphi_0^\eta} \le 
\ttwo - \tone \le \frac{\alpha-\eta-\thone}{\varphi_0-\eta}\,.
$$

By \eqref{tone} for every $\eta\in(0,\eta_1)$ we have $\tone\to t_0$ as $\e\to 0$.
By compactness we can assume that there exist $\Ttwo>t_0$ such that $\ttwo\to\Ttwo$ as $\e\to 0$ along a suitable sequence.
As $\ttwo < t^2_{\e,\delta}$ for $\delta < \eta$, we have $\Ttwo<\Ttwota$ for $\delta< \eta$. This implies that there exists  $\tau\in(t_0,+\infty)$ such that
$\Ttwo\to\tau$ as $\eta\to 0$.
By \eqref{thetabdd}, using a diagonal argument, we can find a subsequence,
still denoted $\theta_\e$, and a locally Lipschitz function $\theta_*\colon(t_0, \tau)\to \R$ such that for every $0<\eta<\eta_1$ and every $t_0<\hat t_1<\hat t_2<\Ttwo$ we have
\begin{eqnarray}
& \theta_\e\to \theta_* \qquad
\text{uniformly in } [\hat t_1, \hat t_2]\,,
\label{funzione}
\\
& \dot\theta_\e\wto \dot \theta_*\qquad\text{weakly$^*$ in }
L^\infty ([\hat t_1, \hat t_2])\,.
\label{derivata}
\end{eqnarray}
Passing to the limit in \eqref{fluss}, we deduce that $|\dot\theta_*(t)-\varphi(\theta_*(t))|\leq \eta$ almost everywhere in $(t_0, \Ttwo)$.
By the arbitrariness of $\eta>0$ we get that $\theta_*$ is a $C^1$ solution of the equation
\begin{equation}\label{lab}
\dot\theta_*(t)=\varphi(\theta_*(t))
\end{equation}
for $t\in(t_0, \tau)$.

Let us fix $\eta\in(0,\eta_1)$. By \eqref{tone} for every $t\in (t_0,\Ttwo)$ we have
that $t\in (\tone,\ttwo)$ for $\e$ small enough, which gives, thanks to 
\eqref{thetabdd},
$$
|\theta_\e(t)-\thone|\leq  (t-\tone)\hat \varphi_0^\eta \leq  (t-t_0)\hat \varphi_0^\eta\,.
$$
Passing to the limit as $\e\to 0$, and using \eqref{thone}, we deduce that
$$
|\theta_*(t)-\theta_0|\leq (t-t_0)\hat \varphi_0^\eta \,.
$$
This shows that $\theta_*$ satisfies $\theta_*(t_0)=\theta_0$. By uniqueness of the solution of  the Cauchy problem \eqref{Ccauchy} we deduce that $\theta_*(t)=\theta(t)$ for every $t\in(t_0,\tau)$. This implies that the full sequence $\theta_\e$ tends to $\theta$, as $\e\to 0$, uniformly on compact subsets of $(t_0,\tau)$. 

As $t^2_{\e,\eta}\to t^2_{\eta}\in (t_0,\tau)$ and $\theta_\e(t^2_{\e,\eta})=\alpha-\eta$, we deduce that  $\theta(\Ttwo)=\alpha-\eta$, and hence $\theta(\tau)=\alpha$.

\medskip\noindent{\it Step $\rm(c)_3$. Behaviour near $\alpha$.}
Assume $\eta\in (0,\eta_1)$ and 
$0<\e<\min\{\eta/(\varphi(\alpha-\eta)+\eta),\e_2(\eta)\}$.
Remark that \eqref{newfluss} implies
$v_\e(\alpha-\eta)<\eta$.
We want to find $\eta_3\in (0,\eta_1)$ such that for every $\eta\in(0,\eta_3)$ there exists $\e_3(\eta)\in(0,\e_2(\eta))$ for which
\begin{equation}\label{vd<e}
v_\e(\theta)\leq\eta
\end{equation}
for every $\e\in(0,\e_3(\eta))$ and every
$\theta\in[\alpha-\eta, \alpha]$. Let  $\hat\theta$ be the largest element of $[\alpha-\eta, \alpha]$ such that \eqref{vd<e} is true for every $\theta\in [\alpha-\eta, \hat\theta]$. 
By \eqref{defR} the rest of equation \eqref{cinqueeqv} can be estimated for $\e<\eta<\eta_1$ and $\theta\in [\alpha-\eta, \hat\theta]$ by
\begin{equation}\label{rone}
|R_\e(\theta,v_\e(\theta))|\leq r_1 \eta
\end{equation}
for a suitable constant $r_1$ independent of $\eta$ and $\e$.
Since $B_0(\theta)>0$ for $\theta\in [\alpha-\eta, \hat\theta]$, by \eqref{A0est} 
and \eqref{cinqueeqv} we have
$$
v'_\e(\theta)\leq \e\frac{\hat a_0}{v_\e(\theta)}+r_1 \eta\,,
$$
where $\hat a_0:=\mu|\xi_0^s|$.

Let $w_\e$ be the solution of the Cauchy problem
$$
w'_\e(\theta)=\e\frac{\hat a_0}{w_\e(\theta)}+r_1 \eta\,, \qquad
w(\alpha-\eta)= \hat a_0\sqrt\e\,.
$$
By \eqref{newfluss} for $\e$ small enough we have 
$v_\e(\alpha-\eta)<\e (\varphi(\alpha-\eta)+\eta) < \hat a_0 \sqrt\e$, which implies, by a comparison argument, that 
\begin{equation}\label{vdwd}
v_\e(\theta)\leq w_\e(\theta)
\end{equation}
for every $\theta\in  [\alpha-\eta, \hat\theta]$. Since $w_\e$ is increasing, we have
$w'_\e(\theta)\leq \sqrt\e +r_1\eta$, which gives 
$$
w_\e(\theta) \leq \hat a_0 \sqrt\e +
(\sqrt\e +r_1\eta)(\theta-\alpha+\eta)\,.
$$
If $\eta<\eta_3:=\min\{1/r_1,\eta_1\}$ and 
$\e<\e_3(\eta):=\min\{(1-\eta_1)^2\eta^2/(\hat a_0+\eta)^2,
{\eta/(\varphi(\alpha-\eta)+\eta)},\allowbreak\e_2(\eta)\}$, by \eqref{vdwd} we have 
$v_\e(\theta)< \eta$ for every $\theta\in  [\alpha-\eta, \hat\theta]$, which implies that
$\hat\theta=\alpha$ and concludes the proof of \eqref{vd<e}.

Let $\tthree$ be the time such that $\theta_\e(\tthree)=\alpha$.
Let us prove that $\tthree-\ttwo$ is uniformly small. To this aim, we observe that  
\eqref{A0est}, \eqref{cinqueeqv}, and \eqref{rone} give
$$
v'_\e(\theta)\geq\e\frac{a_0}{v_\e(\theta)}-B_0(\theta)-r_1\eta\quad
\hbox{for every }\theta\in (\alpha-\eta,\alpha)\,,
$$
where $a_0:=-\sqrt3 \mu |\xi_0^s|V'(\theta_0)>0$. Using \eqref{BBest}, \eqref{vd<e}, and the inequality $r_1\eta<1$, we obtain
for $0<\eta<\eta_3$ and $0<\e<\e_3(\eta)$
\begin{equation}\label{t32}
\tthree - \ttwo = \int_{\alpha-\eta}^\alpha
\frac{\e}{v_\e(\theta)}\, d\theta \leq \frac{1}{a_0}\int_{\alpha-\eta}^\alpha \Big(
v'_\e (\theta)+  B_0(\theta )+ r_1\eta\Big)\, d\theta
\leq \frac{2\mu+2}{a_0}\eta\,.
\end{equation}

Let us prove that
\begin{equation}\label{t3tau}
\lim_{\e\to 0}\tthree=\tau\,.
\end{equation}
For $0<\eta<\eta_3$ and $0<\e<\e_3(\eta)$ we have
$$
|\tthree-\tau|\le |\tthree-\ttwo| +|\ttwo-\Ttwo| +|\Ttwo-\tau|\,.
$$
Using \eqref{t32}, we deduce that
$$
\limsup_{\e\to 0}|\tthree-\tau|\le \frac{2\mu+2}{a_0}\eta +|\Ttwo-\tau|\,.
$$
Passing to the limit as $\eta\to 0$, we obtain \eqref{t3tau}.

\medskip\noindent{\it Step $(c)_4$. Behaviour between $\alpha$ and $\Phi(\alpha)$.}
Let $w_\alpha$ be the solution of the Cauchy problem \eqref{cauchy-fast} with
$\gamma=\alpha$. We want to find $\eta_4\in(0,\eta_3)$ such that for every 
$\eta\in(0,\eta_4)$
there exist $\e_4(\eta)\in(0,\e_3(\eta))$
for which
\begin{equation}\label{flux33}
 w_\alpha(\theta) -\eta\le  v_\e(\theta)\le w_\alpha(\theta)+ \eta
\end{equation}
for every $\e\in(0,\e_4(\eta))$ and every $\theta\in[\alpha,\Phi(\alpha)-\eta]$.
This will be done by comparing the equations \eqref{cauchy-fast} and \eqref{neweqv}
satisfied by $w_\alpha$ and $v_\e$, respectively.

By  \eqref{A1est}, \eqref{A2est}, and  \eqref{defBd} we have 
\begin{equation}\label{bd>}
B(\theta,v)- a_1\e\le B_\e(\theta,v)\leq B(\theta,v)+ a_2\e
\end{equation}
for $\theta\geq\theta_0$ and $v\in[0,-V'(\theta)]$, were $a_1:=4\mu|\xi_0^s|$ and  $a_2:=-\mu|\xi_0^s|/V'(\theta_0)$.
For every $\rho\in \R$ let $z_\rho$ be the solution of the the Cauchy problem
\begin{equation}\label{zrho}
z_\rho'(\theta)=\rho+B(\theta, z_\rho(\theta))\,, \qquad
z_\rho(\alpha)=\rho
\end{equation}
By comparing \eqref{neweqv} with this equation, we obtain that, if $0<\e<-\rho/a_1$, then
$v_\e(\theta)\ge z_\rho(\theta)$ for 
$\theta\ge \alpha$. By the continuous dependence on initial data 
$z_\rho\to w_\alpha$ uniformly on $[\alpha,\Phi(\alpha)]$ as $\rho\to 0$.
Therefore, for every $\eta\in (0,\eta_3)$ there exists $\rho_1(\eta)>0$ such that 
$z_{\rho}(\theta)>w_\alpha(\theta)-\eta$ for every $\rho\in[-\rho_1(\eta),0]$ and every
$\theta\in[\alpha,\Phi(\alpha)]$. It follows that, if $\eta\in (0,\eta_3)$ and $\e<\e_4^1(\eta):=\rho_1(\eta)/a_1$, then the first inequality in
\eqref{flux33} is satisfied for every $\theta\in[\alpha,\Phi(\alpha)]$.

In the previous step we have proved that $v_\e(\alpha)\le \eta$ for $\eta<\eta_3$ and
$\e<\e_3(\eta)$. Let us fix $\eta\in(0,\eta_3)$, $\rho\in (0,1)$, and $\e\in (0,\e_3(\eta))$ with $\e<\rho\eta$. If $v_\e(\alpha)\ge \e/\rho$,  we set $\thfour:=\alpha$. If $v_\e(\alpha)< \e/\rho$, let $\thfour$ be the greatest
element of $[\alpha,\beta]$ such that $v_\e(\theta)\le \e/\rho$ for every 
$\theta\in[\alpha,\thfour]$.
{}From \eqref{btv} and \eqref{bd>} it follows that $B_\e(\theta,v)\geq -a_1\e$
for $\theta\in[\alpha,\beta]$ and $v\in[0,-V'(\theta)]$.
Using \eqref{A0est}, \eqref{neweqv},  and \eqref{<V'} we deduce that
$v_\e'(\theta)\geq\e ( a_0/v_\e(\theta))-a_1\e$, with 
$ a_0:=-\sqrt3 |\xi_0^s| V'(\theta_0)>0$,
hence $v_\e(\theta) v_\e'(\theta)\geq  a_0\e-a_1\e^2/\rho$ for every 
$\theta\in [\alpha,\thfour]$. 
If $\e<\e_4^2(\eta,\rho):=\min\{\rho  a_0/(2a_1), \rho\,\eta, \e_4^1(\eta)\}$ we obtain 
$v_\e(\theta) v_\e'(\theta)\geq \e  a_0/2$ for every 
$\theta\in [\alpha,\thfour]$,
hence, by integrating,
\begin{equation}\label{vquad}
v_\e(\theta)^2\geq  a_0\e(\theta -\alpha)\quad \hbox{ for every }
\theta\in[\alpha,\thfour]\,.
\end{equation}
It follows that
\begin{equation}\label{timenew}
0\leq \thfour - \alpha \le \frac{\e}{ a_0\rho^2}\,.
\end{equation}
In particular, if $v_\e(\alpha)< \e/\rho$ and 
$\e<\e_4^3(\eta,\rho):=\min\{(\beta-\alpha) a_0\rho^2,\e_4^2(\eta,\rho)\}$, we have 
$\thfour<\beta$ and 
\begin{equation}\label{verho}
v_\e(\thfour)=\frac\e\rho
\end{equation}
and
\begin{equation}\label{veeta}
v_\e(\theta)\le \eta\qquad \hbox{for every }\theta\in [\alpha,\thfour]\,.
\end{equation}

Let $\tauf$ be the time for which
$\theta_\e(\tauf)=\thfour$.
Using \eqref{vquad} and \eqref{timenew} we deduce that
\begin{equation}\label{test}
\tauf - \tthree = \e \int_{\alpha}^{\thfour} \!\!\!
\frac{d\theta}{v_\e(\theta)}
\le \frac{\sqrt\e}{\sqrt{ a_0}}  \int_{\alpha}^{\thfour} \!\!\!
\frac{d\theta}{\sqrt{\theta-\alpha}} =
\frac{2\sqrt\e}{\sqrt{ a_0}}\sqrt{\thfour -\alpha}
\le \frac{2\e}{ a_0\rho}\le \frac{2}{a_0}\eta\,.
\end{equation}

Suppose that 
$\eta< \eta_4^1:=\min\{ a_0/(3a_1), w_\alpha(\beta)/3,\eta_3\}$, 
$\rho\in(0,1)$, and $\e<\e_4^4(\eta,\rho):=\min\{\rho \eta_4^1, \e_4^3(\eta,\rho)\}$.
Let
$\theta^4_\e$ be the greatest element of $[\alpha,\beta]$ such that
\begin{equation}\label{test2}
v_\e(\theta)\le 2 \eta_4^1 \qquad\hbox{for every }\theta\in[\alpha, \theta^4_\e]\,.
\end{equation}
By \eqref{veeta} we have $\thfour<\theta^4_\e$.
By the first inequality in \eqref{flux33} we have
$v_\e(\beta)>w_\alpha(\beta)-\eta\ge 2 \eta_4^1$, hence $\theta^4_\e<\beta$ and
\begin{equation}\label{test4}
v_\e(\theta^4_\e)= 2 \eta_4^1 \,.
\end{equation}

We want to prove that $\theta^4_\e-\alpha$ is larger than a positive constant, independent of $\e$.
Using \eqref{A0est}, \eqref{neweqv}, \eqref{bd>}, and \eqref{<V'} we deduce that
$v_\e'(\theta)\geq\e ( a_0/v_\e(\theta))-\e a_1 \ge
\e(( a_0/(2\eta_4^1))-a_1)> 0$ for every 
$\theta\in[\alpha, \theta^4_\e]$. This implies that $v_\e$ is increasing on 
$[\alpha, \theta^4_\e]$, therefore
\begin{equation}\label{test3}
v_\e(\theta)\ge \frac\e\rho\qquad \hbox{for every }\theta\in [\thfour, \theta^4_\e]\,.
\end{equation}
As $v_\e(\alpha)\le \eta<\eta_4^1$ and $v_\e(\theta^4_\e)=2\eta_4^1$, there exists a unique point $\hat \theta^4_\e$ in 
$(\alpha, \theta^4_\e)$ such that $v_\e (\hat \theta^4_\e)=\eta_4^1$. It follows that
$v_\e(\theta)\ge \eta_4^1$ for every $\theta\in[\hat\theta^4_\e, \theta^4_\e]$. By
\eqref{A0est}--\eqref{B3est} there exists a constant $b_1>0$ such that 
$B_\e(\theta,v)\le b_1$ for every $\e\in(0,1)$, $\theta\in[\theta_0,+\infty)$, and 
$v\in[0,-V'(\theta)]$. 
Let $\hat a_0:=\mu|\xi_0^s|$. Using \eqref{A0est}, \eqref{neweqv}, and \eqref{<V'} we obtain
$$
v'_\e(\theta)\le \frac{\e \hat a_0}{v_\e(\theta)}+b_1\le
\e\frac{ \hat a_0}{\eta_4^1}+b_1\le 2b_1
$$
for every $\theta\in[\hat\theta^4_\e, \theta^4_\e]$ and every 
$\e<\e^5_4(\eta,\rho):=\min\{\eta_4^1 b_1/ \hat a_0, \e_4^4(\eta,\rho)\}$. Therefore \eqref{test4} gives
$$
\eta_4^1= v_\e(\theta^4_\e)- v_\e(\hat \theta^4_\e)=
\int_{\hat \theta^4_\e}^{\theta^4_\e} v'_\e(\theta)\,d\theta\le
 2b_1(\theta^4_\e-\hat \theta^4_\e)\,,
$$
which implies
\begin{equation}\label{test21}
\theta^4_\e-\alpha\ge \theta^4_\e-\hat \theta^4_\e \ge \frac{\eta_4^1}{2b_1}\,.
\end{equation}

For every $\eta<\eta^2_4:=\min\{\beta-\alpha,\Phi(\alpha)-\beta, \eta_4^1/(2b_1),
\eta_4^1\}$ we define
$$
m_\eta:=\min \{w_\alpha(\theta): \alpha+\eta\le \theta\le \Phi(\alpha)-\eta\}>0\,.
$$
Since $m_\eta\to 0$ as $\eta\to 0$, there exists $\eta_4\in (0, \eta^2_4)$ such that
$m_\eta<\eta_3<1$ for every $\eta<\eta_4$.
Using \eqref{test21} and  the first inequality in \eqref{flux33} (with $\eta$ replaced by 
$m_\eta/2$) we obtain that 
$$
v_\e(\theta)\ge\frac{m_\eta}2\ge \frac\e\rho \qquad \text{for }
\theta\in[\theta^4_\e,\Phi(\alpha)-\eta]
$$
provided that $\eta<\eta_4$, $\rho<1$, and $\e<\e^6_4(\eta,\rho):=
\min\{\e^1_4(m_\eta/2),\rho\, m_\eta/2, \e^5_4(\eta,\rho)\}$. Using also 
\eqref{test3} we obtain
\begin{equation}\label{vest}
v_\e(\theta)\geq \frac\e\rho \qquad \text{for every }
\theta\in[\thfour,\Phi(\alpha)-\eta]\,.
\end{equation}
 
Let $\tfour$ be the time for which $\theta_\e(\tfour)=\Phi(\alpha)-\eta$.
From the previous estimate we have
\begin{equation}\label{test20}
\tfour-\tauf =\e\int_{\thfour}^{\Phi(\alpha)-\eta}\!
\frac{d\theta}{v_\e(\theta)}\le \rho (\Phi(\alpha)-\alpha)<\eta(\Phi(\alpha)-\alpha)\,,
\end{equation}
if we have also $\rho<\eta$.

By \eqref{A0est}, \eqref{neweqv},  \eqref{bd>}, \eqref{<V'}, and \eqref{vest} we obtain
$$
v'_\e(\theta)\leq(\hat a_0+a_2)\rho +B(\theta, v_\e(\theta))
 \qquad \text{for every }
\theta\in[\thfour,\Phi(\alpha)-\eta]\,,
$$
whenever $0<\eta<\eta_4$, $0<\rho<\eta$, and $0<\e<\e^6_4(\eta,\rho)$
(recall that $\e^6_4(\eta,\rho)<\rho$). Replacing $\rho$ by $\rho/({\hat a_0+a_2})$ we obtain
$$
v'_\e(\theta)\leq \rho +B(\theta, v_\e(\theta))
 \qquad \text{for every }
\theta\in[\thfour,\Phi(\alpha)-\eta]\,,
$$
whenever $0<\eta<\eta_4$, $0<\rho<\rho_2(\eta):=\min\{\eta,\eta(\hat a_0+a_2)\}$, and 
$0<\e<\e^7_4(\eta,\rho):=\min\{\e^6_4(\eta,\rho),\e^6_4(\eta,\rho/(\hat a_0+a_2))\}$.

If
$0<\e<\e^8_4(\eta,\rho):=\min\{\e^3_4(\rho,\rho),\e^7_4(\eta,\rho)\}$, by \eqref{veeta} 
we have $v_\e(\thfour)\le\rho$. By \eqref{btv} the solution $z_\rho$ of \eqref{zrho} is
increasing on $[\alpha,\beta]$, hence $z_\rho(\thfour)\ge\rho$. By comparison we have
\begin{equation}\label{vd<}
v_\e(\theta)\leq z_\rho(\theta)\qquad \hbox{for every }
\theta\in [\thfour,\Phi(\alpha)-\eta]\,.
\end{equation}
On the other hand by the continuous dependence on the data 
$z_\rho\to w_\alpha$ 
uniformly on $[\alpha,\Phi(\alpha)]$ as $\rho\to 0$.  
Therefore, for every $\eta<\eta_4$ there exists $\rho_4(\eta)\in(0, \rho_3(\eta))$ such that 
$z_\rho(\theta)\le w_\alpha(\theta)+\eta$ for every $\rho\in(0,\rho_3(\eta))$ and every $\theta\in [\alpha,\Phi(\alpha)]$.
By \eqref{vd<} we have
$v_\e(\theta)\leq w_\alpha(\theta)+\eta$ for every 
$\theta\in [\theta^4_{\e,\rho_4(\eta)},\Phi(\alpha)-\eta]$ when $\eta<\eta_4$ and 
$\e<\e_4(\eta):=
\min\{\e^8_4(\eta,\rho_4(\eta)),\eta^2\}$, which proves the 
second inequality in \eqref{flux33} on the interval $[\theta^4_{\e,\rho_4(\eta)},\Phi(\alpha)-\eta]$.
The same upper bound on the interval 
$[\alpha,\theta^4_{\e,\rho_4(\eta)}]$ follows from~\eqref{veeta}.

\medskip\noindent{\it Step $\rm(c)_5$. Behaviour near $\Phi(\alpha)$.} 
By \eqref{flux33} we obtain in particular
$$
w_\alpha(\Phi(\alpha)-\eta^2)-\eta^3\leq v_\e(\Phi(\alpha)-\eta^2)
\leq w_\alpha(\Phi(\alpha)-\eta^2)+ \eta^3
$$
for every $\eta\in(0,\eta_4)$ and every $\e\in(0,\e_4(\eta^3))$.
As $w_\alpha'(\Phi(\alpha))=-B_0(\Phi(\alpha))<0$, there exist two constants
$c_2, c_3>0$ and a constant $\eta^1_5\in (0,\eta_4)$ such that 
\begin{equation}\label{vdtra}
c_2\eta^2\leq v_\e(\Phi(\alpha)-\eta^2)\leq c_3 \eta^2
\end{equation}
for $\eta\in (0,\eta^1_5)$. Note that for $\e<\e^1_5(\eta):=\min\{c_2\eta^2/(\varphi(\Phi(\alpha)-\eta^2)+\eta),\e_4(\eta^2), \e_4(\eta^3)\}$ we have
$$
\e(\varphi(\Phi(\alpha)-\eta^2)+\eta)< v_\e(\Phi(\alpha)-\eta^2)\,.
$$

For any $\eta<\eta^1_5$ and
$\e<\e^1_5(\eta)$ let $\thfive$ be the largest element of 
$[\Phi(\alpha)-\eta^2,\Phi(\alpha)+1]$
such that 
$$
\e(\varphi(\theta)+\eta)\leq v_\e(\theta)\leq c_3\eta^2
$$
for all $\theta\in[\Phi(\alpha)-\eta^2,\thfive]$. 
Using \eqref{A1}--\eqref{B3est} the rest  $R_\e$ in equation \eqref{cinqueeqv} 
can be estimated by
\begin{equation}\label{res}
|R_\e(\theta,v_\e(\theta))|\leq r_2(\e+\eta)\,,
\end{equation}
for every $\theta\in[\Phi(\alpha)-\eta^2,\thfive]$, where $r_2$ is a constant independent 
of $\eta$ and $\e$. Since 
$$
\frac{\e}{v_\e(\theta)}\leq\frac{1}{\varphi(\theta)+\eta}
$$
for every $\theta\in[\Phi(\alpha)-\eta^2,\thfive]$, recalling \eqref{defphi} we obtain from \eqref{cinqueeqv} 
$$
v'_\e(\theta)\leq-\frac{\eta}{\varphi(\theta)+
\eta}B_0(\theta)+R_\e(\theta, v_\e(\theta))
$$
in the same interval. Let $\hat b_0>0$ be the infimum of $B_0(\theta)$ for
$\theta\geq (\Phi(\alpha)+\beta)/2$, and let $\hat \varphi_0$ be the supremum of
${\varphi(\theta)+1}$ on the same half-line. By \eqref{bdvarphi}  and  \eqref{res} in the interval 
$[\Phi(\alpha)-\eta^2,\thfive]$
we have the estimate
$$
v'_\e(\theta)\leq -\frac{\hat b_0}{\hat \varphi_0}\eta+2r_2\eta\leq
 -\frac{\hat b_0}{2\hat \varphi_0}\eta\,,
$$
provided that $\eta<\eta^2_5:=\min\{\hat b_0/(4r_2\hat \varphi_0),\eta^1_5\}$
and  $0<\e<\e^2_5:=\min\{\eta,\e^1_5(\eta)\}$. 
As $v_\e(\theta)>0$ by Lemma~\ref{lm:td} and by \eqref{newv},  we obtain
$$
\frac{\hat b_0\eta}{2\hat \varphi_0}\big(\thfive-\Phi(\alpha)+\eta^2\big)\leq 
v_\e(\Phi(\alpha)-\eta^2)-v_\e(\thfive)\leq c_3\eta^2\,,
$$
where the last inequality follows from \eqref{vdtra}. Therefore
there exists $\eta_5\in (0, \eta^2_5)$ and for every $\eta\in (0,\eta_5)$ there exists 
$\e_5(\eta)\in(0, \e_4(\eta))$ such that 
\begin{equation}\label{theta3}
0<\thfive-\Phi(\alpha)+\eta^2 \le \frac{2 c_3\hat \varphi_0}{\hat b_0}\eta<1
\end{equation}
for every $\e\in(0,\e_5(\eta))$.
As $v'_\e(\theta)<0$ on $[\Phi(\alpha)-\eta^2,\thfive]$, from the maximality of $\thfive$
we deduce that 
\begin{equation}\label{maximal}
v_\e(\thfive)=\e(\varphi(\thfive)+\eta)\,.
\end{equation}
 
Let $\tfive$ be the time for which
$\theta_\e(\tfive)= \thfive$.
Since $v_\e(\theta)\geq\e(\varphi(\theta)+\eta)$ for every 
$\theta\in[\Phi(\alpha)-\eta^2,\thfive]$, using \eqref{bdvarphi} we obtain that 
$$
\tfive - \tfourr= \e\int_{\Phi(\alpha)-\eta^2}^{\thfive}
\frac{d\theta}{v_\e(\theta)} \le \int_{\Phi(\alpha)-\eta^2}^{\thfive}
\frac{d\theta}{\varphi(\theta)+\eta} \le
\frac{1}{\varphi_0}(\thfive-\Phi(\alpha)+\eta^2)\,,
$$
which, by \eqref{theta3},  implies 
\begin{equation}\label{test30}
\tfive - \tfourr\le \frac{2c_3\hat \varphi_0}{\hat b_0 \varphi_0}\eta
\end{equation}
for every $\eta\in (0,\eta_5)$ and for every $\e\in (0,\e_5(\eta))$.

\medskip\noindent{\it Step $\rm(c)_6$. Behaviour between  $\Phi(\alpha)$ and $+\infty$.}  
Using \eqref{maximal} and arguing as in Step $\rm(c)_2$ we can prove that for every
$\eta\in(0,\eta_5)$ there exists $\e_6(\eta)\in (0,\e_5(\eta))$ such that
$$
\Big|\frac{v_\e(\theta)}{\e}-\varphi(\theta)\Big|<\eta
$$
for every $\e\in (0,\e_6(\eta))$ and every $\theta\in (\thfive, +\infty)$. 

From \eqref{theta3} it follows that, passing to a subsequence, we have 
$$
\thfive\to \theta^5_\eta\qquad\text{as }\e\to 0\qquad\hbox{and}\qquad
\theta^5_\eta\to \Phi(\alpha)\qquad\text{as }\eta\to 0\,.
$$
Using \eqref{test}, \eqref{test20}, and \eqref{test30} we get
\begin{eqnarray*}
&
\tfive-\tthree=(\tfive-\tfourr)+(\tfourr-\tau^4_{\e,\rho_3(\eta^2)})+
(\tau^4_{\e,\rho_3(\eta^2)}-\tthree)\le
\\
&\displaystyle
\le \frac{2}{a_0}\eta^2+ (\Phi(\alpha)-\alpha)\eta^2+ 
\frac{2c_3\hat \varphi_0}{\hat b_0 \varphi_0}\eta
\end{eqnarray*}
for every $\eta\in (0,\eta_5)$ and for every $\e\in (0,\e_5(\eta))$ (recall that this implies
$\e<\e_4(\eta^2)$, hence 
$\e< \e^4_4(\eta^2,\rho_3(\eta^2))$ and $\e < \e^6_4(\eta^2,\rho_3(\eta^2))$).
Passing to a subsequence we obtain
$$
\tfive\to t^5_\eta\qquad\text{as }\e\to 0
\qquad\hbox{and}\qquad
t^5_\eta\to \tau\qquad\text{as }\eta\to 0\,.
$$
Arguing as in Step $\rm(c)_2$, we deduce that
$\theta_\e(t)$ 
converges to the solution $\theta(t)$ to the Cauchy problem \eqref{Phicauchy} uniformly on compact subsets of $(\tau,+\infty)$. This concludes the proof of case~(c).

\bigskip

In case (b) the proof is divided into three steps.

\medskip\noindent{\it Step  $\rm(b)_1$. Behaviour between $\theta_0$ and $\Phi(\theta_0)$.}
It is enough to repeat the proof of Step  $\rm(c)_4$ with $\alpha$ replaced by $\theta_0$.

\medskip\noindent{\it Step $\rm(b)_2$. Behaviour near $\Phi(\theta_0)$.}
It is enough to repeat the proof of Step $\rm(c)_5$ with $\alpha$ replaced by $\theta_0$.

\medskip\noindent{\it Step  $\rm(b)_3$. Behaviour between $\Phi(\theta_0)$ and $+\infty$.}
It is enough to repeat the proof of Step $\rm(c)_6$ with $\alpha$ replaced by $\theta_0$.

\bigskip

In case (a) we first assume that $B_0(\theta)>0$ for $\theta\geq\theta_0$. This happens when $\mu>\mu_0$, or when $\mu=\mu_0$ and $\theta_0>\alpha_0$, or when $\mu<\mu_0$ and $\theta_0>\beta$. In all these cases the proof is divided into two steps.

\medskip\noindent{\it Step  $\rm(a)_1$. Behaviour near $\theta_0$.}
This is identical to Step $\rm(c)_1$.

\medskip\noindent{\it Step $\rm(a)_2$. Behaviour between $\theta_0$ and $+\infty$.}
It is enough to repeat the proof of Step $\rm(c)_6$ with $\Phi(\alpha)$ replaced by $\theta_0$.

\bigskip

In the case $\mu<\mu_0$ and $\theta_0=\beta$, or $\mu=\mu_0$ and
$\theta_0=\alpha_0$, we have only to modify the proof of Step $\rm(a)_1$. 
First of all we find 
$\eta_1\in(0,1)$ such that for every $\eta\in(0,\eta_1)$ there exists
$\e_1(\eta)\in(0,\eta)$ for which 
$$
v_\e(\theta)\le\eta
$$
for every $\e\in(0,\e_1(\eta))$ and every $\theta\in[\beta,\beta+\eta]$
(with $\beta=\alpha_0$ for $\mu=\mu_0$). This can be proved as in Step $\rm(c)_3$, replacing $\alpha-\eta$ by $\beta$. 

If $v_\e(\beta+\eta)>\e(\varphi(\beta+\eta)+\eta)$, we define $\thstar$ as the largest element of $[\beta+\eta,\beta+1]$ such that 
$$
\e(\varphi(\theta)+\eta)\leq v_\e(\theta)\leq \eta
$$
for all $\theta\in[\beta+\eta,\thstar]$. As in Step $\rm(c)_5$, we can prove that $\thstar\to\beta+\eta$
as $\e\to 0$.

If $v_\e(\beta+\eta)<\e(\varphi(\beta+\eta)-\eta)$, we define $\thstar$ as the largest element of $[\beta+\eta,\beta+1]$ such that 
$$
v_\e(\theta)\leq \e(\varphi(\theta)-\eta)
$$
for all $\theta\in[\beta+\eta,\thstar]$. As in Step $\rm(c)_1$, we can prove that $\thstar\to\beta+\eta$
as $\e\to 0$.

If $\e(\varphi(\beta+\eta)-\eta)\le v_\e(\beta+\eta)\le \e(\varphi(\beta+\eta)+\eta)$, we define $\thstar:=\beta+\eta$.

This concludes Step $\rm(a)_1$ in the case $\mu<\mu_0$ and $\theta_0=\beta$.

\bigskip

It remains to study the case $\mu=\mu_0$ and $\theta_0<\alpha_0$. The proof is split into four steps.

\medskip\noindent{\it Step $\rm(a')_1$. Behaviour near $\theta_0$.}
This is identical to Step $\rm(c)_1$.

\medskip\noindent{\it Step $\rm(a')_2$. Behaviour between $\theta_0$ and $\alpha_0$.}
This is identical to Step $\rm(c)_2$.

\medskip\noindent{\it Step $\rm(a')_3$. Behaviour near $\alpha_0$.}
The behaviour in the interval $[\alpha_0-\eta,\alpha_0]$ can be studied as in Step $\rm(c)_3$, while the behaviour in a right neighbourhood of $\alpha_0$ can be studied as in Step $\rm(a)_1$ for the case $\mu<\mu_0$ and $\theta_0=\beta$.

\medskip\noindent{\it Step $\rm(a')_4$. Behaviour between $\alpha_0$ and $+\infty$.}
It is enough to repeat the proof of Step $\rm(c)_6$ with $\beta$ replaced by $\alpha_0$.
\end{proof}

\end{section}

\begin{section}{Examples with concentrations and oscillations}\label{conc-osc}

In this section we consider solutions of the reduced problem in dimension $d=1$ considered in Section~\ref{simple-shear}
and corresponding to simple shears.
To simplify the exposition we take
\begin{equation}\label{1171}
\mu=\tfrac12 \qquad\hbox{and}\qquad
K^{\!R}:=\{(\alpha,\zeta)\in\R{\times}\R: |\alpha|+|\zeta|\le 2\} \,,
\end{equation}
so that 
\begin{equation}\label{1953}
\sigmaa^R(t)=\ee^R(t)\quad \hbox{and}\quad H^R(\xi,\theta)=2(|\xi|\lor|\theta|)\,.
\end{equation}
About the softening potential $V$ we assume that 
\begin{equation}\label{1957}
V(\theta)=V(-\theta)\qquad\hbox{for every }\theta\in \R
\end{equation}
and 
\begin{equation}\label{1353}
V'(\theta)=-1\text{ if and only if }\theta\geq1\,.
\end{equation}
Together with the general assumptions \eqref{2der} and \eqref{V'infty}, these conditions imply that $-1\le V'(\theta)\le 1$ for every $\theta\in \R$ and
\begin{equation}\label{315}
-1<V'(\theta)\le 0\quad\hbox{for }0\le \theta<1\,.
\end{equation}

We will consider the boundary datum
\begin{equation}\label{315w}
\ww^R(t,y):=ty 
\end{equation}
and  the initial condition $(u^R_0, e^R_0, p^R_0, z^R_0)$, where
\begin{equation}\label{317}
u^R_0=0\,, \quad 
e^R_0=0\,, \quad p^R_0=0
\quad\hbox{on } [-\tfrac12, \tfrac12]\,,
\end{equation}
and $z^R_0\colon [-\tfrac12, \tfrac12]\to\R$ is a suitable piecewise continuous function. 

\subsection{Strain localization}
The following theorem describes
an example where the strain has a singular component supported by a lower 
dimensional manifold.
\begin{theorem}\label{th:concentration} Besides \eqref{1171}--\eqref{317}, assume 
that $z^R_0$ is continuous and $z^R_0(0)=1>z^R_0(y)\ge 0$ for every 
$y\neq 0$. Let $(\uu^R, \ee^R, \muu^R)$ be a solution of the reduced problem
for the approximable quasistatic evolution 
with boundary datum $\ww^R$ and initial condition
$(u^R_0, e^R_0, p^R_0, z^R_0)$. Then for every 
$t\in[0,+\infty)$ we have
\begin{eqnarray}
&\displaystyle
\uu^R(t,y)=
\begin{cases}
ty&\hbox{if }\, t\in[0,1]\,,\ y\in[-\tfrac12, \tfrac12]\,,
\\
y+t-1&\hbox{if }\, t\in(1,+\infty)\,,\ y\in(0, \tfrac12]\,,
\\
y-t+1&\hbox{if }\, t\in(1,+\infty)\,,\ y\in[-\tfrac12, 0)\,,
\end{cases}
\label{pl20}
\\
\nonumber
\\
&
\ee^R(t)=t\land1\,,\qquad \muu^R_t=\delta_{(\pp^R(t),\zz^R(t))}\,,
\label{pl22}
\end{eqnarray}
 where
\begin{equation}\label{pl2}
\pp^R(t):=(t-1)^+\delta_0\,,
\quad 
\zz^R(t):=z^R_0+(t-1)^+\delta_0\,,
\end{equation} 
$\delta_0$ being the unit Dirac mass concentrated at $0$.
Moreover, for every $T>0$ we have 
\begin{equation}\label{pl21}
\D_{\!H^{\!R}}(\muu^R;0, T)=2(T-1)^+\,,
\end{equation}
and the energy inequality (ev4) of Theorem~\ref{Thm55} holds with equality.
\end{theorem}

\begin{remark}\label{shear-band}
In this example we see the phenomenon of strain localization at $y=0$ even if the boundary and initial data are smooth. If the reduced 
one-dimensional problem is used to describe the evolution of simple shears, this example exhibits the formation of a shear band on the plane $x_1=0$ for every 
time~${t> 1}$.
\end{remark}

\begin{proof}[Proof of Theorem~\ref{th:concentration}] 
By Remark~\ref{rem1d} for every $\e>0$ the solution $(\uu^R_\e,\ee^R_\e,\pp^R_\e,\zz^R_\e)$ of the reduced $\e$-regularized evolution problem with boundary datum
$\ww^R$ and initial condition
$(u^R_0, e^R_0, p^R_0, z^R_0)$ satisfies \eqref{oneD1} and \eqref{oneD2}. Moreover $\zz^R_\e(t,y)$ is continuous with respect to $(t,y)$.
We observe that
\begin{equation}\label{785}
N^\e_{K^{\!R}}(\alpha,\zeta)= 
\frac{1}{2\e}\big( (\alpha+\zeta-2)^+, (\alpha+\zeta-2)^+\big)
\end{equation}
for every $\e>0$ and every $(\alpha,\zeta)\in \R{\times}\R$ with 
$\alpha\ge 0$ and $-2\le \zeta-\alpha\le 2$.

Let $t_\e$ be the largest element of $[0,+\infty]$ such that 
\begin{equation}\label{765}
-2\le -V'(\zz^R_\e(t,y))-\sigmaa^R_\e(t)\le 2
\end{equation}
for every $t\in[0,t_\e)$ and every $y\in[-\frac12,\frac12]$.
By \eqref{oneD1}, \eqref{oneD2}, and \eqref{785} we have
\begin{eqnarray}
& \displaystyle 
\dot\ee^R_\e(t) = 1-\frac{1}{2\e}\int_{-\frac12}^{\frac12} (\ee^R_\e(t)- 
V'(\zz^R_\e(t,y))-2)^+\, dy\,,  \label{eq81}
\\
& \displaystyle 
\dot\zz^R_\e(t,y)= \frac{1}{2\e} (\ee^R_\e(t)-V'(\zz^R_\e(t,y))-2)^+\qquad\hbox{for every }
y\in[-\tfrac12,\tfrac12] \label{eq82}
\end{eqnarray}
for every $t\in[0,t_\e)$.

For $t\in[0,1]$ the functions 
\begin{equation}\label{t<1}
\ee^R_\e(t)=t\qquad\text{and}\qquad \zz^R_\e(t,y)=z_0(y)
\end{equation}
solve \eqref{eq81} and \eqref{eq82} and satisfy the initial conditions
$\ee^R_\e(0)=0$ and $\zz^R_\e(0,y)=z^R_0(y)$. Moreover they satisfy \eqref{765} with strict inequalities. Therefore $t_\e> 1$ and
\eqref{t<1} gives the solution of the reduced $\e$-regularized evolution problem for 
$t\in[0,1]$.

We observe that
\begin{equation}\label{eta83}
\zz^R_\e(t,y)\geq z^R_0(y)
\end{equation}
for every $t\in [0,t_\e)$ and every $y\in [-\tfrac12,\tfrac12]$.
This follows easily from \eqref{eq82} and from the initial condition
$\zz^R_\e(0,y)= z^R_0(y)$.

Let us prove that
\begin{equation}\label{t>0}
\ee^R_\e(t)>1
\end{equation}
for every $t\in(1,t_\e)$. Indeed, this is true for $t$ near $1$, since $\ee^R_\e(1)=1$
and $\dot\ee^R_\e(1)=1$. If \eqref{t>0} is not satisfied for every $t\in(1,t_\e)$,
we can consider the first point $\tau_\e>1$ where $\ee^R_\e(\tau_\e)=1$.
It is clear that $\dot\ee^R_\e(\tau_\e)\le0$, but from \eqref{eq81} we obtain
$\dot\ee^R_\e(\tau_\e)=1$. This contradiction proves~\eqref{t>0}.

Let us prove that
\begin{equation}\label{eta84}
\lim_{\e\to 0}\sup_{1< t<t_\e}\ee^R_\e(t)=1\,.
\end{equation}

Since the function $x\mapsto V'(z^R_0(x))$ is continuous and $z^R_0(0)=1$,
for every $\eta>0$ there exists $\delta>0$ such that
$$
-V'(z^R_0(y))\geq -V'(1)-\eta =1-\eta \,,
$$
whenever $|y|\le\delta$. Together with \eqref{eta83} this implies that
$$
-V'(z^R_\e(t,y))\geq 1-\eta \,,
$$
for $1< t<t_\e$ and $|x|\le\delta$. Therefore, \eqref{eq81} implies
$$
\dot \ee^R_\e(t) \le1 -\frac{1}{2\e}\int_{-\delta}^\delta 
(\ee^R_\e(t)- V'(\zz^R_\e(t,y))-2)^+\, dy
\le  1-\frac{\delta}{2\e} \big(\ee^R_\e(t)-1-\eta \big)^+  \,.
$$
As $\ee^R_\e(1)=1$, by a comparison argument we obtain
$$
\sup_{1<t<t_\e} \ee^R_\e(t)\leq 1+\eta +\frac{2\e}{\delta}\,.
$$
This implies 
$$
\limsup_{\e\to 0}\, \sup_{1<t<t_\e} \ee^R_\e(t) \leq 1+\eta\,,
$$
which, together with \eqref{t>0}, gives \eqref{eta84} by the arbitrariness of $\eta$.

We deduce from \eqref{t<1}, \eqref{t>0}, and  \eqref{eta84} that 
\begin{equation}\label{1299}
\ee^R_\e(t)\to t\land 1
\end{equation}
uniformly with respect to $t\ge0$. It follows from Definition~\ref{def:qsym}
that the first equality in
\eqref{pl22} is satisfied for every $t\in \Theta$, and hence for every $t\in[0,+\infty)$ 
by left-continuity.

As $z^R_0(y)\ge0$, 
by \eqref{315} and \eqref{eta83} we have $-V'(\zz^R_\e(t,y))\ge 0$. Therefore, by \eqref{eta84} there exists $\e_0>0$ such that $-3/2\le -V'(\zz^R_\e(t,y))-\sigmaa^R_\e(t)\le 1$ for every $\e\in (0,\e_0)$,  $t\in [0,t_\e)$, and  $y\in[-\tfrac12,\tfrac12]$. By continuity and by \eqref{765}  this implies that 
$t_\e=+\infty$ for every $\e\in (0,\e_0)$.

For every $\e\in (0,\e_0)$ and every $t>1$ we define
\begin{equation}\label{1300}
a_\e(t):=\sup\{|y|:\ \zz^R_\e(t,y)>z^R_0(y)\}\,.
\end{equation}
Since $\zz^R_\e(t,y)$ is increasing with respect to $t$ by \eqref{eq82}, 
we have that $a_\e(t)$
is increasing with respect to $t$. Let
\begin{equation}\label{1301}
a_\e(+\infty):=\lim_{t\to+\infty} a_\e(t)\,.
\end{equation}
Let us prove that
\begin{equation}\label{limce}
\lim_{\e\to0} a_\e(+\infty)=0\,.
\end{equation}
Let us fix $\eta>0$. As  $z^R_0(y)<1$ for $|y|\geq\eta$, by \eqref{315}
\begin{equation}\label{maxV}
\max_{|y|\geq\eta} -V'(z^R_0(y))=1-\delta_\eta 
\end{equation}
for some $\delta_\eta>0$. By \eqref{eta84} there exists $\e_\eta\in(0,\e_0)$ such that $\ee^R_\e(t)\le1+\delta_\eta$ for every $\e<\e_\eta$ and every $t\geq 1$. Therefore, if 
$|y|\geq \eta$, by \eqref{maxV} we obtain
$$
(\ee^R_\e(t)-V'(z^R_0(y))-2)^+=0
$$
for every $\e<\e_\eta$ and every $t\geq 1$. This implies $\zz^R_\e(t,y)=z^R_0(y)$ by the uniqueness of the solution of \eqref{eq82} with the initial condition $z^R_0(y)$.
We deduce that $a_\e(t)\le\eta$ for every $\e<\e_\eta$ and every $t\geq1$,
which gives \eqref{limce}.

To continue the proof of the theorem we need the following lemma 
about concentrations of positive functions or measures.

\begin{lemma}\label{lmiota}
Let $U$ be a bounded open set in $\Rn$, $d\ge 1$, and let $p_k$ be a sequence in
$M_b^+(\ol U)$ which converges to $p$
weakly$^*$ in $M_b(\ol U)$. Assume that $p$ is singular with respect to the Lebesgue measure. Then $\delta_{p_k}\wto \delta_p$ weakly$^*$ in $GY(\ol U;\R)$.
\end{lemma}

\begin{proof}
Passing to a subsequence we may assume that $\delta_{p_k}\wto \mu$ weakly$^*$ in $GY(\ol U;\R)$. Then $p=\bary(\mu)$ by \eqref{weakbar}.
As $\supp\, \delta_{p_k}\subset \ol U{\times}[0,+\infty){\times}[0,+\infty)$, we have also
\begin{equation}\label{sppos}
\supp\,\mu \subset \ol U{\times}[0,+\infty){\times}[0,+\infty)\,.
\end{equation}

Let us consider the representation of $\mu$ given by \cite[Theorem~4.3 and Remark~4.5]{DM-DeS-Mor-Mor-1}:
\begin{equation}\label{8.10}
\begin{array}{c}
\displaystyle \langle f,\mu\rangle = \int_{\ol U{\times}\R} f(x,\xi,1)\,d\mu^Y(x,\xi) +
\int_{\ol U{\times}\Sigma_\R}f(x,\xi,0)\,d\mu^\infty(x,\xi) =
\vspace{.2cm}
\\
\displaystyle = \int_{U}\Big(\int_\R
f(x,\xi,1)\,d\mu^{x,Y}(\xi)\Big)\, dx+
\int_{\ol U} \Big(\int_{\Sigma_\R} f(x,\xi,0)\,d\mu^{x,\infty}(\xi)\Big)\, d\lambda^\infty(x)\,,
\end{array}
\end{equation}
where $\mu^{x,Y}$ and $\mu^{x,\infty}$ are probability measures and 
$\lambda^\infty:=\pi_{\ol U}(\mu^\infty)$.
By \eqref{sppos} we have $\supp\,\mu^{x,Y}\subset [0,+\infty)$ and
$\supp\,\mu^{x,\infty}=\{1\}$, hence $\mu^{x,\infty}=\delta_1$, the unit Dirac mass
at~$1$.

Let 
$$
u^Y(x):=\int_0^{+\infty}\xi \, d\mu^{x,Y}(\xi)\,, \qquad
u^\infty(x):=\int_{\Sigma_\R} \xi \, d\mu^{x,Y}(\xi)=1 \,.
$$
By \cite[Remark~6.3]{DM-DeS-Mor-Mor-1} we have
$$
p=\bary(\mu)=u^Y+u^\infty\lambda^\infty=u^Y+\lambda^\infty\,.
$$
Since $p$ is singular, this implies $u^Y=0$ a.e.\ on $U$ and $\lambda^\infty=p$. 
As $\mu^{x,Y}$ is a probability measure, it follows from the definition of $u^Y$ that $\mu^{x,Y}=\delta_0$. Therefore, \eqref{8.10} yields
$$
\langle f,\mu\rangle =  \int_{U} f(x,0,1)\, dx + \int_{\ol U} f(x,1,0)\, dp(x)\,,
$$
which, by \eqref{mup}, is equivalent to $\mu=\delta_p$.
\end{proof}

{\it Proof of Theorem~\ref{th:concentration} (Continuation)}
By \eqref{dualdef} and \eqref{785} for $\e<\e_0$ we have 
$\dot\zz^R_\e(t,y)=\dot\pp^R_\e(t,y)$, so that the initial conditions \eqref{317} give
\begin{equation}\label{pl}
\zz^R_\e(t,y)=\pp^R_\e(t,y)+z^R_0(y)
\end{equation} 
for every $t\ge 0$ and every $y\in[-\tfrac12,\tfrac12]$. In particular \eqref{eta83}
 implies that
\begin{equation}\label{1302}
\pp^R_\e(t,y)\ge 0\quad\hbox{for every }y\in[-\tfrac12,\tfrac12]\,,
\end{equation}
while \eqref{1300} and  \eqref{1301}  yield
$\pp^R_\e(t,y)=0$ for $|y|\ge a_\e(+\infty)$.
We deduce that for every $t\ge0$ there exist a sequence $\e_j\to 0$ and a constant
$\varphi(t)\ge0$ such that $\pp^R_{\e_j}(t)\wto\varphi(t)\delta_0$ weakly$^*$ in 
$M_b([-\tfrac12,\tfrac12])$. Since $\ee^R_{\e_j}(t)\to t\land 1=\ee^R(t)$ and 
$(\uu^R_{\e_j}(t), \ee^R_{\e_j}(t), \pp^R_{\e_j}(t))$ satisfies condition (ev1)$\!_\e$ of Definition~\ref{def:reym}, arguing as in the proof of Theorem~\ref{Thm55}
we obtain that $\uu^R_{\e_j}(t)$ converges weakly$^*$ in 
$BV([-\tfrac12,\tfrac12])$ to a function $\uu^*(t)\in BV([-\tfrac12,\tfrac12])$ such that
$(\uu^*(t), \ee^R(t), \varphi(t)\delta_0)$ satisfies condition (ev1)$\!^R$ of
Theorem~\ref{Thm55}. 
Using this property and \eqref{315w} we deduce that $\uu^*(t)$ coincides with the right-hand side of \eqref{pl20} and $\varphi(t)=(t-1)^+$. 
As the limits do not depend on the sequence  $\e_j$, we obtain that
\begin{eqnarray}
&\label{1348}
\uu^R_\e(t)\wto \uu^*(t)\quad\hbox{weakly}^* \hbox{ in }BV([-\tfrac12,\tfrac12])\,,
\\
&\label{1349}
\pp^R_\e(t)\wto \pp^R(t) \quad\hbox{weakly}^* 
\hbox{ in }M_b([-\tfrac12,\tfrac12])\,,
\\
&\label{1649}
\zz^R_\e(t)\wto \zz^R(t)
\quad\hbox{weakly}^* \hbox{ in }M_b([-\tfrac12,\tfrac12])\,,
\end{eqnarray} 
where $\pp^R(t) $ and $\zz^R(t)$ are defined by~\eqref{pl2}.
It follows from Definition~\ref{def:qsym}
that \eqref{pl20}  is satisfied for every $t\in \Theta$, and hence for every $t\in[0,+\infty)$ 
by left-continuity.

Let $\psi\colon\![-\tfrac12,\tfrac12]{\times}\R{\times}\R\to  [-\tfrac12,\tfrac12]{\times}\R{\times}\R{\times}\R$ be defined by 
$\psi(y,\beta,\eta):=(y,\beta,\beta+\eta z^R_0(y),\eta)$.  It follows from \eqref{pl} that 
\begin{equation}\label{zzpp}
\delta_{(\pp^R_\e(t),\zz^R_\e(t))}=\psi(\delta_{\pp^R_\e(t)})\,, \qquad 
\delta_{(\pp^R(t),\zz^R(t))}=\psi(\delta_{\pp^R(t)})\,.
\end{equation} 
By \eqref{1302} and \eqref{1349} we obtain that
$\delta_{\pp^R_\e(t)}\wto \delta_{\pp^R(t)}$ weakly$^*$ in $GY([-\tfrac12,\tfrac12];\R)$ thanks to Lemma~\ref{lmiota}. By \eqref{zzpp}
we conclude that 
\begin{equation}\label{1350}
\delta_{(\pp^R_\e(t),\zz^R_\e(t))}\wto \delta_{(\pp^R(t),\zz^R(t))}
\end{equation} 
weakly$^*$ in $GY([-\tfrac12,\tfrac12];\R{\times}\R)$. 
It follows from Definition~\ref{def:qsym}
that  the second equality in
\eqref{pl22}  is satisfied for every $t\in \Theta$, and hence for every $t\in[0,+\infty)$ 
by left-continuity.

To prove \eqref{pl21} it is enough to show that that for every $t_1$ and $t_2$, with $0\le t_1<t_2$,
$$
\langle H^R(\beta_2-\beta_1,\theta_2-\theta_1),
\muu^R_{t_1t_2}(y,\beta_1,\theta_1,\beta_2,\theta_2,\eta)\rangle
=2\big((t_2-1)^+ - (t_1-1)^+\big)\,.
$$
Taking into account \eqref{1953} and the left continuity, we have to prove that
\begin{equation}\label{1352}
\langle |\beta_2-\beta_1|\lor |\theta_2-\theta_1|,
\muu^R_{t_1t_2}(y,\beta_1,\theta_1,\beta_2,\theta_2,\eta)\rangle
=(t_2-1)^+ - (t_1-1)^+
\end{equation} 
for every $t_1$, $t_2\in\Theta$, with $0\le t_1<t_2$. We know that for every 
$t_1$, $t_2\in\Theta$, with $0\le t_1<t_2$, there exists a sequence ${\e_j}\to 0$ such that
$$
\delta_{((\pp^R_{\e_j}(t_1),\zz^R_{\e_j}(t_1)), (\pp^R_{\e_j}(t_2),\zz^R_{\e_j}(t_2)))} \wto \muu^R_{t_1t_2}\,,
$$
hence
\begin{equation}\nonumber
\begin{array}{c}
\displaystyle\vphantom{\int_{-1}^1}
\langle |\beta_2-\beta_1|\lor |\theta_2-\theta_1|,
\muu^R_{t_1t_2}(y,\beta_1,\theta_1,\beta_2,\theta_2,\eta)\rangle=
\\
\displaystyle = \lim_{j\to\infty}
\int_{-\frac12}^{\frac12}|\pp^R_{\e_j}(t_2,y)-\pp^R_{\e_j}(t_1,y)|
\lor|\zz^R_{\e_j}(t_2,y)-\zz^R_{\e_j}(t_1,y)|\,dy=
\\
\displaystyle = \lim_{j\to\infty}
\int_{-\frac12}^{\frac12}\big(\pp^R_{\e_j}(t_2,y)-\pp^R_{\e_j}(t_1,y)\big)\,dy\,,
\end{array}
\end{equation}
where in the last equality we used \eqref{eq82} and~\eqref{pl}. Therefore \eqref{1349}, together with the definition of $\pp^R(t)$ given in \eqref{pl2}, yields~\eqref{1352}. 

The equality in (ev4) of  Theorem~\ref{Thm55} follows now from \eqref{pl22},
\eqref{pl2},  and~\eqref{pl21}.
\end{proof}

\subsection{Oscillation of the internal variable}
The following theorem describes an example where strain localization 
occurs together with a strong oscillation 
of the internal variable.

\begin{theorem}\label{th:oscillation}
Besides \eqref{1171}--\eqref{317}, assume that the function $z^R_0$ is odd on 
$[-\tfrac12, \tfrac12]$, continuous on  $[-\tfrac12,0)\cup(0,\tfrac12]$, and satisfies $z^R_0(0+)=1>z^R_0(y)\ge 0$ for every $y\neq 0$. 
Let $(\uu^R, \ee^R, \muu^R)$ be a solution of the reduced problem 
for the approximable quasistatic evolution with boundary datum $\ww^R$
and initial condition
$(u^R_0, e^R_0, p^R_0, z^R_0)$. Then for every 
$t\in[0,+\infty)$ we have
\begin{eqnarray}
&\displaystyle
\uu^R(t,y)=
\begin{cases}
ty&\hbox{if }\, t\in[0,1]\,,\ y\in[-\tfrac12,\tfrac12]\,,
\\
y+t-1&\hbox{if }\, t\in(1,+\infty)\,,\ y\in(0,\tfrac12]\,,
\\
y-t+1&\hbox{if }\, t\in(1,+\infty)\,,\ y\in[-\tfrac12,0)\,,
\end{cases}
\label{pl20*}
\\
\nonumber
\\
&\ee^R(t)=t\land1\,,\qquad \muu^R_t=\tfrac12 \delta_{(\pp^R(t),z^R_0+\pp^R(t))}+\tfrac12 \delta_{(\pp^R(t), z^R_0-\pp^R(t))}\,,
\label{pl22*}
\end{eqnarray}
where
\begin{equation}\label{pl2*}
\pp^R(t):=(t-1)^+\delta_0\,,
\end{equation} 
$\delta_0$ being the unit Dirac mass concentrated at~$0$.
Moreover, for every $T>0$ we have 
\begin{equation}\label{pl21*}
\D_{\!H^{\!R}}(\muu^R;0, T)=2(T-1)^+\,,
\end{equation}
and the energy inequality (ev4) of Theorem~\ref{Thm55} holds with equality.
\end{theorem}

\begin{remark}\label{noenergy}
Let $\zz^R$, $\zetaa^R\colon [0,+\infty)\to L^2([-\tfrac12,\tfrac12])$ be the functions
defined by $\zz^R(t):=z^R_0$ and $\zetaa^R(t):=-V'(z^R_0)$ for every $t\geq 0$.
It follows from \eqref{pl22*} and \eqref{pl21*} that 
$$
\zetaa^R(t)=-\pi_{\ol\Om}(\{V'\}\muu^R_t)\qquad\hbox{and}\qquad 
(\pp^R(t),\zz^R(t))=\bary(\muu^R_t)
$$
for every $t\ge 0$, so that in this case the stress constraint (ev3) of 
Theorem~\ref{Thm55} depends on the Young measure $\muu^R_t$ 
only through its barycentre.
However, the barycentres  $\pp^R$ and $\zz^R$ do not satisfy the energy 
inequality~\eqref{ineqreduced}. Indeed, 
if $T>1$, by \eqref{pl2*}  we have 
\begin{eqnarray*}
&\mu \|\ee^R(T)\|_2^2=\tfrac12\,,\quad\D_{\!H^{\!R}}(\pp^R, \zz^R;0, T)=2(T-1)\,, \\
&\V(\zz^R(0))=\V(\zz^R(T))=\V(z^R_0)\,,\quad\langle\sigmaa^R(t), D\dot\ww^R(t)\rangle=t\land 1\,,
\end{eqnarray*}
which contradicts \eqref{ineqreduced}. Note that by \eqref{pl22*}, \eqref{pl2*}, and \eqref{pl21*} we have
\begin{eqnarray*}
&\D_{\!H^{\!R}}(\muu^R;0, T)=2(T-1)\,,
\\
&
\langle \{V\}(\theta,\eta), \muu^R_{T}(y, b, \theta, \eta)\rangle=\V(z^R_0)-(T-1)\,,
\end{eqnarray*}
so that the energy inequality (ev4) of  Theorem~\ref{Thm55} holds with equality.

This example shows that, if we consider just the barycentres, which coincide with the 
weak$^*$ limits of the solutions $(\uu^R_\e,\ee^R_\e,\pp^R_\e,\zz^R_\e)$ of the reduced
$\e$-regularized evolution problems, the energy inequality cannot be obtained because,
by the lack of convexity, we neglect some important terms generated by the space
oscillations of the approximate solutions, that are captured only by the Young measure
formulation.

As in Theorem~\ref{th:concentration}, in Theorem~\ref{th:oscillation} we see the phenomenon of strain localization at $y=0$, even if the boundary and initial data are smooth. In the latter we see also,  for every time
$t\ge 1$, a strong oscillation of the inner variable $\zz_\e$ concentrated near the point where the plastic strain is localized. 
\end{remark}

\begin{proof}[Proof of Theorem~\ref{th:oscillation}] 
Let $\bar z^R_0:[-\tfrac12,\tfrac12]\to \R$ be the function defined by 
\begin{equation}\label{zetabar}
\begin{array}{l}
\bar z^R_0(y)=z^R_0(y)\quad\hbox{for }y\in (0,\tfrac12]\,,\\
 \bar z^R_0(0)=z^R_0(0+)\,,
 \\ 
 \bar z^R_0(y)=-z^R_0(y)\quad \hbox{for }y\in [-\tfrac12,0)\,.
\end{array}
\end{equation}
It turns out that $\bar z^R_0$ is even and  continuous, and satisfies 
$\bar z^R_0(0)=1>\bar z^R_0(y)$ for every $y\neq 0$.
Let   $(\uu^R_\e,\ee^R_\e,\pp^R_\e,\zz^R_\e)$ be the solution of the reduced 
$\e$-regularized evolution problem with boundary datum $\ww^R$ and initial condition
$(u^R_0, e^R_0, p^R_0, z^R_0)$ and let
 $(\bar\uu^R_\e,\bar\ee^R_\e,\bar\pp^R_\e,\bar \zz^R_\e)$ be the solution of the reduced 
$\e$-regularized evolution problem with boundary datum $\ww^R$ and initial condition
$(u^R_0, e^R_0, p^R_0, \bar z^R_0)$. By the symmetries of the problem 
we can prove that for every $t\geq 0$ and every $y\in[-\tfrac12,\tfrac12]$ we have 
\begin{equation}\label{uguali}
\begin{array}{c}
\uu^R_\e(t,y)=\bar\uu^R_\e(t,y)=-\bar\uu^R_\e(t,-y),\quad \ee^R_\e(t)= \bar\ee^R_\e(t)\,,
\smallskip
\\
\pp^R_\e(t,y)=\bar\pp^R_\e(t,y)=\bar\pp^R_\e(t,-y)\,, 
\end{array}
\end{equation}
and
\begin{equation}\label{zetabare}
\begin{array}{c}
\zz^R_\e(t,y)= \bar \zz^R_\e(t,y)=\bar \zz^R_\e(t,-y) \quad\hbox{for }y\in (0,\tfrac12]\,,
\smallskip
\\
 \zz^R_\e(t,y)=-\bar \zz^R_\e(t,y)=-\bar \zz^R_\e(t,-y)\quad \hbox{for }y\in [-\tfrac12,0)\,.
\end{array}
\end{equation}

Since $\bar z^R_0$ satisfies the hypothesis of Theorem~\ref{th:concentration}, we deduce from  \eqref{1299}, \eqref{1348},  \eqref{1350} , and \eqref{uguali} that 
\begin{eqnarray}
&\label{1348*}
\uu^R_\e(t)\wto \uu^*(t)\quad\hbox{weakly}^* \hbox{ in }BV([-\tfrac12,\tfrac12])\,,
\\
&\label{1299*}
\ee^R_\e(t)\to t\land 1\,,
\\
&\label{1349*}
\delta_{(\bar\pp^R_\e(t),\bar \zz^R_\e(t))}\wto \delta_{(\pp^R(t),\bar z^R_0+\pp^R(t))}\quad\hbox{ weakly$^*$ in $GY([-\tfrac12,\tfrac12];\R{\times}\R)$,}
\end{eqnarray} 
where $\uu^*(t)$ is defined as the right-hand side of \eqref{pl20*} and $\pp^R(t)$ is defined in~\eqref{pl2*}.

Let us prove that for every $t\ge 0$
\begin{equation}\label{1370}
\delta_{(\pp^R_\e(t),\zz^R_\e(t))}\wto 
\tfrac12\delta_{(\pp^R(t),z^R_0+\pp^R(t))}+ 
\tfrac12\delta_{(\pp^R(t),z^R_0-\pp^R(t))}\,.
\end{equation} 
To this aim, we fix $f\in C^\hom([-\tfrac12,\tfrac12]{\times}\R{\times}\R{\times}\R)$ and observe that,
by \eqref{uguali} and \eqref{zetabare}, we have
\begin{eqnarray*}
&\displaystyle
\langle f,\delta_{(\pp^R_\e(t),\zz^R_\e(t))}\rangle =
\int_{-\frac12}^{\frac12}f(y,\pp^R_\e(t,y), \zz^R_\e(t,y),1)\,dy=
\\
&\displaystyle
=\frac12 \int_{-\frac12}^{\frac12} f_1(y,\bar\pp^R_\e(t,y), \bar\zz^R_\e(t,y),1)\,dy+
\frac12 \int_{-\frac12}^{\frac12} f_2(y,\bar\pp^R_\e(t,y), \bar\zz^R_\e(t,y),1)\,dy\,,
\end{eqnarray*}
where
\begin{eqnarray}
&\displaystyle
f_1(y,\beta,\theta,\eta)=
\begin{cases}
f(y,\beta,\theta,\eta)&\hbox{if }y\ge 0\,,
\\
f(-y,\beta,\theta,\eta)&\hbox{if }y< 0\,,
\end{cases}
\label{1415}
\\
&\displaystyle  
f_2(y,\beta,\theta,\eta)=
\begin{cases}
f(-y,\beta,-\theta,\eta)&\hbox{if }y\ge 0\,,
\\
f(y,\beta,-\theta,\eta)&\hbox{if }y< 0\,.
\end{cases}
\label{1416}
\end{eqnarray}
Using \eqref{1349*} we get
$$
\lim_{\e\to 0} \langle f,\delta_{(\pp^R_\e(t),\zz^R_\e(t))}\rangle=
\tfrac12 \langle f_1,\delta_{(\pp^R(t),\bar z^R_0+\pp^R(t))}\rangle
+
\tfrac12 \langle f_2,\delta_{(\pp^R(t),\bar z^R_0+\pp^R(t))}\rangle\,.
$$
Since $\pp^R(t)=(t-1)^+\delta_0$, using \eqref{mup}, \eqref{1415}, and \eqref{1416} 
we obtain
$$
\begin{array}{c}
\tfrac12 \langle f_1,\delta_{(\pp^R(t),\bar z^R_0+\pp^R(t))}\rangle
+
\tfrac12 \langle f_2,\delta_{(\pp^R(t),\bar z^R_0+\pp^R(t))}\rangle=
\smallskip
\\
=
\tfrac12 \langle f,\delta_{(\pp^R(t), z^R_0+\pp^R(t))}\rangle
+
\tfrac12 \langle f,\delta_{(\pp^R(t), z^R_0-\pp^R(t))}\rangle\,,
\end{array}
$$
which concludes the proof of~\eqref{1370}.

It follows from \eqref{1348*},  \eqref{1299*}, \eqref{1370}, and from 
Definition~\ref{def:qsym}
that
\eqref{pl20*} and \eqref{pl22*} are satisfied for every $t\in \Theta$,
and hence for every $t\in[0,+\infty)$ 
by left-continuity.

Equality \eqref{pl21*} is proved as in Theorem~\ref{th:concentration}, as well as the equality in~(ev4).
\end{proof}

\begin{remark}\label{nonsymm}
Using Remark~\ref{rem1d} we deduce from the symmetry of $K^{\!R}$ that the
solutions of the reduced
$\e$-regularized evolution problems with $z^R_0(y)\ge0$ for every 
$y\in[-\tfrac12,\tfrac12]$ satisfy the inequality $\zz^R_\e(t,y)\ge0$ for every $t\ge0$
and every $y\in[-\tfrac12,\tfrac12]$. Therefore any solution of the reduced problem
for the approximable quasistatic evolution satisfies $\zz^R(t,y)\ge0$
for every $t\ge0$ and every $y\in[-\tfrac12,\tfrac12]$. 
In view of this property the hypotheses on $z^R_0$ in Theorem~\ref{th:oscillation} may seem artificial. This is not the case if, for instance, $K^{\!R}$ is replaced by the hexagon
$$
K_{hex}:=\{ (\alpha, \zeta)\in\R{\times}\R: |\alpha+\zeta|\le2\,,\ |\alpha-\zeta|\le2\,,\
|5\alpha-\zeta|\le7\}\,,
$$
with vertices $(\frac32,\frac12)$, $(0,2)$, $(-\frac54,\frac34)$, $(-\frac32,-\frac12)$,
$(0,-2)$, $(\frac54,-\frac34)$. Indeed,
let $\ww_T(t,y)=\varphi_T(t)y$ for every $t\ge0$ and every $y\in[-\tfrac12,\tfrac12]$, where
$$P
\varphi_T(t):=\begin{cases}
t & \text{if } 0\le t\le T\,, \\
2T-t & \text{if } t\ge T\,.
\end{cases}
$$
Let $(\uu_T,\ee_T,\pp_T,\zz_T)$ be an approximable quasistatic evolution for the problem in dimension $d=1$ corresponding to $\C\xi=\xi$ and 
$K=K_{hex}$, with boundary datum $\ww_T$ and  initial condition 
$(u_0, e_0, p_0, z_0)$. We assume  
$u_0=0$, $e_0=0$, $p_0=0$, and $z_0(y)=\theta_0$ for every 
$y\in[-\tfrac12,\tfrac12]$, where  $\zeta_0:=-V'(\theta_0)\in(-\frac12,\frac12)$.
By Remark~\ref{rm:hom} the functions $\ee_T(t)$, $\pp_T(t)$, $\zz_T(t)$ do not depend on $y$. 
We can study the evolution of the corresponding functions $\sigmaa_T(t)$ and
$\zetaa_T(t)$ as in 
Theorem~\ref{thm:ex1} and plot them as paths of 
$(\sigmaa_T(t),\zetaa_T(t))$ in the plane $(\sigma,\zeta)$ (see Figure~4).
\begin{figure}[ht]
\begin{center} 
  \includegraphics[width=0.6\textwidth]{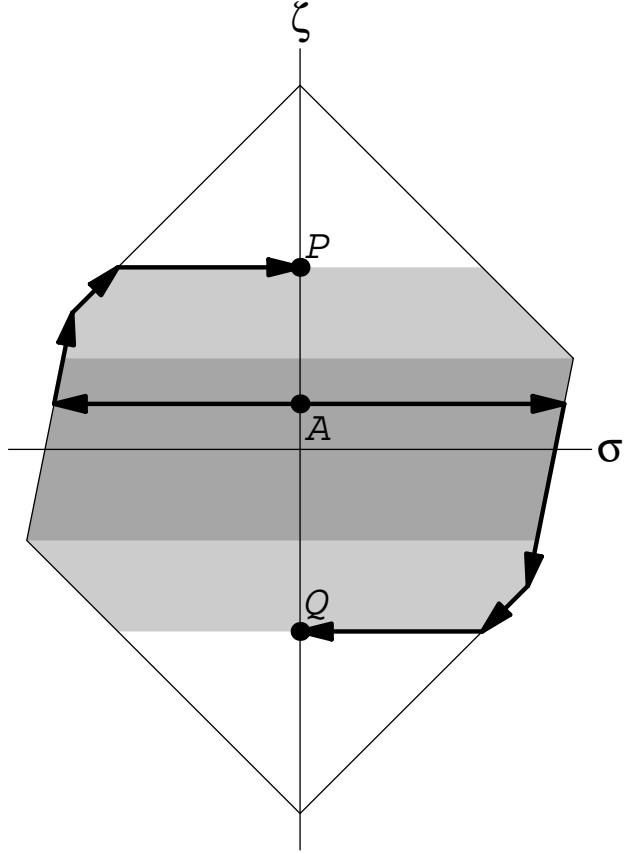}
\caption{A path $(\sigmaa_T(t),\zetaa_T(t))$ from $A:=(0,\zeta_0)$ to  $P:=(0,1)$ corresponding to an input $\varphi_T(t)$,  together with a path from $A$  to $Q:=(0,-1)$ corresponding to a different input.}
\end{center}
\end{figure}
Since $V'(\theta)=-1$ for $\theta\le-1$, it is possible to prove that there exists $T>0$ such that $(\sigmaa_T(T),\zetaa_T(T))=(-1,1)$. 
It is then easy to see that $(\sigmaa_T(t),\zetaa_T(t))=
(t-T-1,1)$ for $T\le t\le T+2$ and $(\sigmaa_T(t),\zetaa_T(t))=(1,1)$ for $t\ge T+2$.
Therefore any point at height $\zeta=1$ can be reached from any point $(0,\zeta_0)$ with $\zeta_0\in(-\frac12,\frac12)$ through an appropriate loading-unloading path (see $AP$ in Figure~4).
Similarly, replacing $\varphi_T$ by $-\varphi_T$, we can show that any point at height
$\zeta=-1$ can be reached by another loading-unloading path (see $AQ$ in Figure~4).

Note that, under the assumptions on $z_0$ considered in 
Theorem~\ref{th:oscillation}, the same oscillation phenomenon near $y=0$ occurs when $K$ is replaced by
$K_{hex}$.

Since $K_{hex}$ does not satisfy \eqref{plusminus}, the corresponding one-dimensional problem  cannot be
used to study evolutions of simple shears in elasto-plastic materials. However, they may be used in the study of uniaxial loading of cylindrical bodies, when the material exhibits a different behaviour in tension and compression
(see Remark~\ref{rem:symmetry}). 
\end{remark}

\end{section}

\bigskip

\noindent {\bf Acknowledgments.} { This work is part of the Project ``Calculus of Variations" 2004, 
supported by the Italian Ministry of Education, University, and Research and of the research project 
``Mathematical Challenges in Nanomechanics" sponsored by Istituto Nazionale di Alta Matematica (INdAM) ``F.~Severi".}

\bigskip
\bigskip

{
}

\end{document}